\documentclass[11pt]{article}
\usepackage{amssymb,amsmath,amsfonts,epsfig}
\usepackage{hyperref}
\usepackage{multirow}

\usepackage{cancel}
\usepackage{hyperref}
\usepackage{color}
\usepackage{graphicx}
\usepackage{tikz}
\usepackage{mathtools}
\usepackage{bbm}

\usepackage{float}

\usepackage{latexsym}
\usepackage{setspace}
\usepackage{caption}
\usepackage{subcaption}
 
\numberwithin{equation}{section}

\newcommand{\eprintN}[1]{{\href{http://arxiv.org/abs/#1}{[\texttt{#1 [hep-th]}]}}}
\newcommand{\eprintPH}[1]{{\href{http://arxiv.org/abs/#1}{[\texttt{#1 [hep-ph]}]}}}
\newcommand{\eprintNT}[1]{{\href{http://arxiv.org/abs/#1}{[\texttt{#1 [math.NT]}]}}}

\newcommand{\eprintRA}[1]{{\href{http://arxiv.org/abs/#1}{[\texttt{#1 [math.RA]}]}}}

\newcommand{\doi}[2]{{\hypersetup{urlcolor=darkred}\href{http://dx.doi.org/#2}{#1}\hypersetup{urlcolor=blue}}}

\newcommand{\oeis}[1]{{\hypersetup{urlcolor=darkred}\href{https://oeis.org/#1}{#1}\hypersetup{urlcolor=blue}}}

\renewcommand{\thefootnote}{\fnsymbol{footnote}}
\renewcommand{\thanks}[1]{\footnote{#1}}
\newcommand{\starttext}{
\setcounter{footnote}{0}
\renewcommand{\thefootnote}{\arabic{footnote}}}
\newcommand{\bea}{\begin{eqnarray}}
\newcommand{\eea}{\end{eqnarray}}
\newcommand{\be}{\begin{eqnarray}}
\newcommand{\ee}{\end{eqnarray}}

\def\ie{\begin{equation}\begin{aligned}}
\def\fe{\end{aligned}\end{equation}}

\makeatletter
\newcommand{\setword}[2]{\phantomsection#1\def\@currentlabel{\unexpanded{#1}}\label{#2}}
\makeatother

\def\ie{\begin{equation}\begin{aligned}}
\def\fe{\end{aligned}\end{equation}}

\def\Re{{\rm Re \,}}

\def\Tr{{\rm Tr}}

\newcommand{\mb}{{\mathsf{b}}}

\begin{document}

\begin{flushright}
{\small }
\end{flushright}

\starttext

\setcounter{footnote}{0}

\vskip 0.3in

\begin{center}

\centerline{\large \bf
Resurgent Lambert series with characters} 

\vskip 0.2in
{David Broadhurst$^{1}$ and Daniele Dorigoni$^{2}$} 
\vskip 0.15in
{\small 
$^{1}$ School of Physical Sciences, The Open University, Milton Keynes MK7 6AA, UK\\ 
$^{2}$ Centre for Particle Theory \& Department of Mathematical Sciences,\\ 
Durham University, Lower Mountjoy, Stockton Road, Durham DH1 3LE, UK}
\vskip 0.1in

\begin{abstract}
\vskip 0.1in

We consider certain Lambert series as generating functions of divisor sums twisted by Dirichlet characters and compute their exact resurgent transseries expansion near $q=1^-$.
For special values of the parameters, these Lambert series are expressible in terms of iterated integrals of holomorphic Eisenstein series twisted by the same characters and the transseries representation is a direct consequence of the action of Fricke involution on such twisted Eisenstein series.
When the parameters of the Lambert series are generic the transseries representation provides for a quantum-modular version of Fricke involution which for a particular example we show being equivalent to modular resurgent structures found in topological strings observables.

\end{abstract}  
                     
\end{center}
\nopagebreak 
\newpage

\tableofcontents

\newpage

\section{Introduction} 

The central object of study in this paper is the $q$-series
\begin{equation}
\Xi_{s_1,s_2}(\chi_{r_1},\chi_{r_2}; q) \coloneqq \sum_{n_1=1}^\infty \sum_{n_2=1}^\infty 
\frac{\chi_{r_1}(n_1)}{n_1^{s_1}} \frac{\chi_{r_2}(n_2)}{n_2^{s_2}} q^{n_1 n_2} \,,\label{eq:LamGenInt}
\end{equation}
where the parameters $s_1,s_2\in \mathbb{C}$ while $\chi_{r_1},\chi_{r_2}:\mathbb{Z} \to \mathbb{C}$ denote two Dirichlet characters modulo $r_1$ and $r_2$ respectively.
The series \eqref{eq:LamGenInt} converges absolutely within the complex unit circle $|q|<1$ and we shall be interested in computing the asymptotic expansion of $\Xi_{s_1,s_2}(\chi_{r_1},\chi_{r_2}; q)= \Xi_{s_2,s_1}(\chi_{r_2},\chi_{r_1};q)$ as $q\to 1^-$. 

For the special cases where either $s_1=0$ or $s_2=0$, the $q$-series $\Xi_{s_1,s_2}(\chi_{r_1},\chi_{r_2}; q)$ can be rewritten in terms of generalized Lambert series. For generic $s_1,s_2\in \mathbb{C}$ we can instead think of \eqref{eq:LamGenInt} as a generating function of divisor sums twisted by the two Dirichlet characters $\chi_{r_1},\chi_{r_2}$. Importantly, for particular values of the parameters $s_1,s_2$ related to the choice of characters $\chi_{r_1},\chi_{r_2}$, these generating functions are directly related to holomorphic Eisenstein series twisted by the same characters or iterated integrals thereof. 

Our main results, described in more detail below, are:
\begin{itemize}
\item[(i)] We present in~\eqref{eq:TSgen} the complete transseries asymptotic expansion, including exponentially suppressed terms, for $\Xi_{s_1,s_2}(\chi_{r_1},\chi_{r_2}; q)$ as $q\to 1^-$. Here the parameters $s_1,s_2\in \mathbb{C}$ are kept generic while $\chi_{r_1},\chi_{r_2}$ are primitive characters. 
In~\eqref{eq:TSgenImp}, we present the analogous expression for the case where the characters are imprimitive.
\item[(ii)] For the special cases with either $s_1=0$ or $s_2=0$ and where one of the characters is the trivial character,~i.e. $\chi_r(n) = \chi_{1,1}(n) = 1$ for all $n\in \mathbb{Z}^* $, we obtain in \eqref{eq:TSQM} the complete transseries asymptotic expansion, including exponentially suppressed terms, as $q\to 1^-$. Our expressions reduce to the results of \cite{Dorigoni:2020oon} when restricting both characters to the trivial case.
\item[(iii)] When the parameters $s_1$ and $s_2$ satisfy a particular integrality condition, namely $s_1,s_2\in \mathbb{Z}$ with $s_1+s_2 \equiv \kappa_1 +\kappa_2 +1 \text{ mod }2$, where $\kappa_i\in\{0,1\}$ encodes the parity of the primitive character $\chi_{r_i}$ (i.e. $\chi_{r_i}(-1) = (-1)^{\kappa_i}$), we have that the $q$-series \eqref{eq:LamGenInt} can be expressed in terms of iterated integrals of holomorphic Eisenstein series twisted by the same characters. In Section \ref{sec:ModularPrim}, we show that the complete transseries expansions mentioned above neatly encode the behaviour of  \eqref{eq:LamGenInt} under Fricke involution as inherited from its action on these twisted Eisenstein series.
 
\item[(iv)] Lastly, in Section \ref{sec:Qpoch} we discuss particular examples not related to holomorphic Eisenstein series and specialise the $q$-series \eqref{eq:LamGenInt} to the cases $(s_1,s_2)=(1,0)$ and $(s_1,s_2)=(0,1)$ where one of the characters is fixed to the trivial one. For this choice of parameters and characters we show that the complete transseries expansion of $\Xi_{s_1,s_2}(\chi_{r_1},\chi_{r_2}; q)$ is equivalent to certain strong-weak resurgent structures recently found within the context of topological strings observables in \cite{Rella:2022bwn,Fantini:2024snx,Fantini:2025wap} and particular realizations of a more general \textit{modular resurgence paradigm}~\cite{Fantini:2024ihf}.
\end{itemize}

Our analysis for the complete,~i.e.~perturbative and non-perturbative, transseries asymptotic expansion of~\eqref{eq:LamGenInt} as $q\to 1^-$ relies crucially on resurgent methods~\cite{Ecalle:1981,delabaere1999resurgent}, following closely the discussion of~\cite{Dorigoni:2020oon}. The transseries completion captures terms that are exponentially suppressed in this limit and which are fundamental in understanding the behaviour of the $q$-series~\eqref{eq:LamGenInt} under the action of Fricke involution~$\tau \to -1/(r\tau)$ where the integer $r$ depends on the choice of Dirichlet characters $\chi_{r_1},\chi_{r_2}$, while the ``modular'' parameter $\tau$ is related to $q$ via $q=e^{2\pi i \tau}$.

 In particular, we are interested in something usually referred to as the `modularity gap' which is the failure of the $q$-series to transform as a modular form, or rather a vector-valued modular form, of definite weight under the action of Fricke involution.
Depending on this particular integrality condition for the parameters $s_1$ and $s_2$ mentioned above, this modularity gap is either given by a Laurent polynomial in $\tau$ or by a multi-valued function of $\tau$. 

The case where both Dirichlet characters are chosen to be trivial is rather special in that Fricke involution reduces to the standard $S$-transformation $\tau \to -1/\tau$. Furthermore, given that the $q$-series is manifestly invariant under the $T$-transformation $\tau\to\tau+1$ and since $S$ and $T$ together generate ${\rm PSL}_2(\mathbb{Z})$ when acting on $\tau$, thanks to the transseries expansion we then know the complete behaviour of the $q$-series $\Xi_{s_1,s_2}(\chi_{1,1},\chi_{1,1}; q)$ under the modular group ${\rm PSL}_2(\mathbb{Z})$.
In \cite{Dorigoni:2020oon} it was shown that when $s_1,s_2\in \mathbb{Z}$ such that $s_1+s_2  +1 \in 2\mathbb{Z}$ the modularity gap of $\Xi_{s_1,s_2}(\chi_{1,1},\chi_{1,1}; q)$ reduces to a Laurent polynomial in $\tau$. Furthermore, the transseries expansion of $\Xi_{s_1,s_2}(\chi_{1,1},\chi_{1,1}; q)$ encodes the fact that this $q$-series can be rewritten in terms of iterated integrals of standard holomorphic Eisenstein series.
This type of iterated integrals has recently received a lot of interest both in the maths literature~\cite{Brown:2014,Brown:2017qwo,Brown:2017}, in particular for its connection to periods of moduli spaces of genus-one Riemann surfaces, as well as within the context of string scattering amplitudes~\cite{Broedel:2015hia,Broedel:2018iwv,Broedel:2018izr,Broedel:2019vjc,Dorigoni:2022npe,Dorigoni:2024oft}. 

Our present interest for considering the generalization  $\Xi_{s_1,s_2}(\chi_{r_1},\chi_{r_2}; q)$ with two non-trivial Dirichlet character is threefold.
Firstly, the particular case where \eqref{eq:LamGenInt} is specialized to $(s_1,s_2) = (2,0)$ and $(s_1,s_2) = (0,2)$ with Dirichlet characters $(\chi_{r_1}, \chi_{r_2}) = (\chi_{6,5}, \chi_{1,1})$ has been encountered by Bloch and Vanhove \cite{BV} while computing a particular class of Feynman period integrals named the \textit{sunset graph} and defined as the scalar two-point self-energy at two-loop order. Here $\chi_{6,5}$ denotes the unique odd character modulo $6$ which is induced by the odd primitive character $\chi_{3,2}$.
We will show that for this case the action of the Fricke involution $\tau \to -1/(6\tau)$ on the $q$-series $\Xi_{s_1,s_2}(\chi_{r_1},\chi_{r_2}; q)$ is related to the simple exchange in parameters $(s_1,s_2)= (2,0) \leftrightarrow (s_1,s_2) = (0,2)$. 

The reason for this phenomenon stems from the fact that these Feynman amplitudes correspond, as discussed in \cite{BV}, to a family of periods associated with the universal family of elliptic curves over the modular curve $X_1(6) \cong \mathbb{H} \backslash \Gamma_1(6)$ with $\mathbb{H} = \{ \tau \in \mathbb{C}\,\vert\,{\rm Im}(\tau)>0\}$ and where $\Gamma_1(6)\subset {\rm SL}(2,\mathbb{Z})$ denotes a particular congruence subgroup.\footnote{ For  an integer $N\geq 1 $, the congruence subgroup $\Gamma_1(N) $  is a subgroup of the modular group $ {\rm SL}(2,\mathbb{Z})$ defined by $\Gamma_1(N) \coloneqq \{ \left(\begin{smallmatrix} a & b \\ c & d \end{smallmatrix}\right) \in {\rm SL}(2,\mathbb{Z})\, \vert \, a,d \equiv 1 \, ({\rm mod} \,N)\,,\, c\equiv 0 \,({\rm {mod}} \,N)\}$.}
Correspondingly, we will show that the $q$-series $\Xi_{s_1,s_2}(\chi_{r_1},\chi_{r_2}; q)$, with parameters $(s_1,s_2) = (2,0)$ and $(s_1,s_2) = (0,2)$ and characters as above, can both be represented as iterated integrals of particular modular forms of weight $3$ given by two different holomorphic Eisenstein series twisted by the same characters and which are exchanged upon Fricke involution $\tau \to -1/(6\tau)$.
A natural question which prompted our study is understanding how does this structure change when we change both the parameters $(s_1,s_2)$ and the Dirichlet characters to different ones.

This connection between Feynman periods and iterated integrals of twisted Eisenstein series can be found in other types of diagrams, such as the kite integral \cite{Adams:2017ejb,Broedel:2018rwm,Giroux:2024yxu}. Interestingly, when one expresses the sunrise and kite integrals as power series in the dimensional regularization parameter $\epsilon$ we encounter higher and higher depth iterated integrals of twisted Eisenstein series,~i.e.~nested integrals with more and more twisted Eisenstein series as integration kernels. We believe our analysis will prove to be useful in discussing other types of Feynman periods.

A second, more speculative motivation can be found within the context of string theory amplitudes. As mentioned above, the function $\Xi_{s_1,s_2}(\chi_{1,1},\chi_{1,1}; q)$ when specialized to the case where both Dirichlet characters become trivial is related to the low-energy expansion of one-loop scattering amplitudes of massless open- and closed-string states (via a procedure called \textit{single-valued map}). The same integrality condition mentioned above implies that the relevant Lambert series are expressible in terms of iterated integrals of standard holomorphic Eisenstein series.
However, one-loop string scattering amplitudes also produce more complicated periods which are expressible as iterated integrals whose kernels are not holomorphic Eisenstein series but rather Eisenstein–Kronecker series, see e.g. \cite{Broedel:2017jdo,Hidding:2022vjf}.
Very much along the same lines as \cite{Dorigoni:2020oon} for the case of trivial Dirichlet characters, it is therefore tantalizing to speculate whether \eqref{eq:LamGenInt}, or linear combinations thereof, may be related to such one-loop string amplitudes perhaps when evaluated at some rational value for their insertion points or when (some of) the external states stem from twisted sectors of orbifolds.\footnote{We thank Oliver Schlotterer for related discussions.}.

Lastly, although rather disconnected from our original motivations, a further reason for studying the asymptotic expansion of the $q$-series \eqref{eq:LamGenInt} as $q\to 1^{-}$ can be found within topological string theory compactified on particular Calabi-Yau manifolds.
In particular \cite{Rella:2022bwn} derived resurgent structures both in the weak, i.e. $\tau \to i \infty$, and strong coupling regime, i.e. $\tau \to 0$, for a particular observable called the first fermionic spectral trace in topological string theory compactified on the Calabi-Yau manifold known as local $\mathbb{P}^2$ geometry. In the same reference it was also noticed that the perturbative and non-perturbative data in the two regimes are related by an intriguing `number-theoretic duality'. This duality was later expanded upon in \cite{Fantini:2024snx} and eventually showed in \cite{Fantini:2024ihf}  to be an example of a more general symmetry which was named \textit{modular resurgent structure}. Very recently \cite{Fantini:2025wap} showed that such structures also appear in topological string theory compactified on the broader class of local $\mathbb{P}^{m,n}$ Calabi-Yau manifolds.

We will show that the analysis of~\cite{Rella:2022bwn,Fantini:2024snx,Fantini:2025wap} corresponds to specialising the $q$-series \eqref{eq:LamGenInt} to the cases $(s_1,s_2)=(1,0)$ and $(s_1,s_2)= (0,1)$ with one of the characters reduced to unity.
Unlike the cases relevant for discussing Feynman integrals, the present parameters $s_1$ and $s_2$ do not satisfy the integrality condition previously mentioned. As a consequence, the corresponding $q$-series do not originate from iterated integrals of twisted Eisenstein series and the corresponding modularity gap do not give rise to Laurent polynomials but rather to multi-valued functions of~$\tau$ which can be written as formal, factorially divergent power series. The non-perturbative transseries completion of~\eqref{eq:LamGenInt} then reflects the transformation properties of~$\Xi_{s_1,s_2}(\chi_{r_1},\chi_{r_2}; q)$ under a quantum-modular version of Fricke involution.

We find it striking that both the work of \cite{BV} on the sunset graph and that of \cite{Rella:2022bwn,Fantini:2024snx} on the first fermionic spectral trace in topological string theory on local $\mathbb{P}^2$ rely on the very same $q$-series \eqref{eq:LamGenInt} with the same Dirichlet characters $(\chi_{r_1}, \chi_{r_2}) = (\chi_{3,2}, \chi_{1,1})$ and where the key difference between the two works lies entirely in the value taken by the parameters $(s_1,s_2)$ which consequently change dramatically the corresponding resurgence structures. For \cite{BV} we have $(s_1,s_2)=(2,0)$ and $(0,2)$ thus giving rise to a perturbative asymptotic expansion as $q\to 1^-$ which terminates after finitely many terms and from which we can reconstruct the complete transseries asymptotic expansion via something called \textit{Cheshire cat resurgence}, see \cite{Dunne:2016jsr,Kozcaz:2016wvy}. On the other hand, for \cite{Rella:2022bwn,Fantini:2024snx} the parameters $(s_1,s_2)$ take the values $(s_1,s_2)=(1,0)$ and $(0,1)$ and the perturbative asymptotic expansion as $q\to 1^-$ is given by a formal, factorially divergent power series from which we compute the complete transseries asymptotic expansion using the full resurgent analysis machinery.

It is presently not known whether other instances of $q$-series $\Xi_{s_1,s_2}(\chi_{r_1},\chi_{r_2}; q)$ do play a r\^ole within topological string theory. However, we find that the beautiful story of \cite{Rella:2022bwn,Fantini:2024snx,Fantini:2025wap} provides yet another compelling motivation for analyzing the general resurgent structures of the $q$-series $\Xi_{s_1,s_2}(\chi_{r_1},\chi_{r_2}; q)$, main purpose of the present work.

\subsection*{Outline}

The paper is organized as follows.
In Section \ref{sec:Lam} we introduce the main characters of our story and explain the connections between the $q$-series presented in \eqref{eq:LamGenInt} and generalizations of Lambert series. Using resurgence analysis methods  we then proceed to compute in Section \eqref{sec:Asy}  the complete transseries expansion of $\Xi_{s_1,s_2}(\chi_{r_1},\chi_{r_2}; q)$ as $q\to 1^{-}$. 

In Section \ref{sec:ModularPrim} we show that when the parameters $s_1$ and $s_2$ satisfy the integrality condition $s_1,s_2\in \mathbb{Z}$ with $s_1+s_2 \equiv \kappa_1 +\kappa_2 +1 \text{ mod }2$, where $\chi_{r_i}(-1) = (-1)^{\kappa_i}$, the perturbative asymptotic expansion of $\Xi_{s_1,s_2}(\chi_{r_1},\chi_{r_2}; q)$ reduces to a Laurent polynomial in $\log(q)$ plus an infinite series of non-perturbative, exponentially suppressed corrections. This phenomenon is a beautiful example of a structure dubbed \textit{Cheshire cat resurgence}. Furthermore, we clarify that for these special choice of parameters the $q$-series $\Xi_{s_1,s_2}(\chi_{r_1},\chi_{r_2}; q)$ admits a representation in terms of iterated integrals of holomorphic Eisenstein series twisted by Dirichlet characters. The transseries representation is then a consequence of a particular type of Fricke involution which acts on the space of twisted holomorphic Eisenstein series. 
Section \ref{sec:Ex} contains many concrete examples of this construction of modular primitives for twisted Eisenstein series. 

Finally, in Section \ref{sec:Qpoch} we present an application of our analysis to the study of quantum dilogarithms. In this case, the aforementioned integrality condition on the parameters $s_1$ and $s_2$ does not hold and the perturbative expansion of $\Xi_{s_1,s_2}(\chi_{r_1},\chi_{r_2}; q)$ as $q\to 1^{-}$ produces an asymptotic, factorially divergent power series. Here the transseries expansion encodes the modular resurgent structures found in the study of spectral traces for topological string theory compactified on particular Calabi-Yau manifolds. 

We review some useful background material in the appendices. In Appendix \ref{sec:DirApp} we present some details on Dirichlet characters and their associated $L$-functions. The definition of holomorphic Eisenstein series twisted by Dirichlet characters, their $q$-series expansion and their behaviour under Fricke involution is presented in Appendix \ref{sec:Twist}. Lastly, in Appendix \ref{sec:AsyApp} we show how to use a method due to Don Zagier to extract the asymptotic expansion of $\Xi_{s_1,s_2}(\chi_{r_1},\chi_{r_2}; q)$ as $q\to 1^{-}$ in the case when one of the Dirichlet characters reduces to the trivial one.

\section{Lambert series with characters}
\label{sec:Lam}

The first character of our story occurs in the class of Lambert series 
\begin{equation}
\mathcal{L}_s(\chi_r ; q) \coloneqq \sum_{n= 1}^\infty  \frac{\chi_r(n)}{n^s} \frac{q^n}{1-q^n} =
\sum_{n_1=1}^\infty\sum_{n_2=1}^\infty \frac{\chi_r(n_1)}{n_1^s} q^{n_1 n_2}\,,\label{eq:LamChi} 
\end{equation}
with $s\in \mathbb{C}$ and $\chi_r:\mathbb{Z} \to \mathbb{C}$ a Dirichlet character modulo $r$.
Throughout this work we use the notation $\chi_r$ to denote a general character modulo $r$, while with $\chi_{r,\ell}$ we denote the particular character modulo $r$ specified by the Conrey labels $(r,\ell)$. We refer to Appendix \ref{sec:DirApp} for basic properties of Dirichlet characters and the precise definition of Conrey labels.

The $q$-series \eqref{eq:LamChi} converges for all $q\in \mathbb{C}$ such that $|q|<1$.
Yet  it has a natural boundary of analyticity for $|q|=1$ due to an accumulation of poles for all $q$ near rational roots of unity. 
In particular, we wish to extract the exact asymptotic expansion of \eqref{eq:LamChi} when $q\to 1^{-}$. To this end we change variables
$$
q=\exp(-2\pi y)\,,
$$
and consider the asymptotic expansion of $\mathcal{L}_s(\chi_r ; q) $ for $y\to 0^+$.

Collecting powers of $q$ in the double sum of \eqref{eq:LamChi} we obtain the alternative representation
\begin{equation}
\mathcal{L}_s(\chi_r ; q) = \sum_{n=1}^\infty \frac{\sigma'_{s,\chi_r}(n)}{n^s} q^n\,,\label{eq:qser1} 
\end{equation}
where the twisted divisor function is defined by
\begin{equation}
\sigma'_{s,\chi_r}(n) \coloneqq \sum_{d|n} \chi_r\left(\frac{n}{d}\right) d^s= \sum_{d|n} \chi_r(d) \left(\frac{n}{d}\right)^s\,, \label{eq:Twsigma2}
\end{equation}
and the sums run over the positive divisors of $n$.
\vspace{0.2cm}

\subsection{Lambert series with a principal character} 

Our work generalizes the non-perturbative analysis of \emph{untwisted} Lambert series
in \cite{Dorigoni:2020oon}, for which the identity character $\chi_{1,1}(n)=1$  gives
\begin{equation}
\mathcal{L}_s(\chi_{1,1} ; q) = \mathcal{L}_s(q) \coloneqq \sum_{n=1}^\infty  \frac{\sigma_{s}(n)}{n^s} q^n\,,\label{eq:qserDD}
\end{equation}
where $\sigma_{s}(n) = \sum_{d\vert n}d^s$ denotes the standard divisor function.

For a {\em principal} character modulo $r$, with  $\chi_{r,1}(n)=1$, when $n$ is coprime to $r$, and
$\chi_{r,1}(n)=0$, otherwise, we obtain
\begin{equation}
\mathcal{L}_s(\chi_{r,1} ; q) = 
\sum_{d\vert r}\frac{\mu(d)}{d^{s}}\mathcal{L}_s(q^d)\,,\label{eq:chi11Lam}
\end{equation}
where the M\"obius function is given by 
\begin{equation}
\mu(n) =  \begin{cases} (-1)^{\omega(n)}\,, &{\mbox{if $n$ is square-free}} \\
0\,, &{\rm{otherwise}}\,,
\end{cases}
\end{equation}
that is $\mu(n)=(-1)^{\omega(n)}$, if $n$ is the product of $\omega(n)$ distinct primes, and by $\mu(n)=0$, if $n$ is divisible by the square of a prime.

The relation \eqref{eq:chi11Lam} is a consequence of the identity
\begin{equation}
\sigma'_{s,\chi_{r,1}}(n) = \sum_{d|\gcd(r,n)} \mu(d)  \sigma_s\left(\frac{n}{d}\right)\,,\label{eq:DivId} 
\end{equation}
which we now prove. First, we apply M\"obius inversion to $\sigma_s(n)=\sum_{d|n}d^s$, obtaining
\begin{equation}
n^s=\sum_{d|n}\mu(d) \sigma_s\left(\frac{n}{d}\right)\,.\label{eq:MuInv}
\end{equation}
Next, we combine \eqref{eq:Twsigma2} with \eqref{eq:MuInv}, to obtain
\begin{equation}\sigma'_{s,\chi_{r,1}}(n) =\sum_{\stackrel{d_1|n}{\gcd(r,d_1)=1}}\left(\frac{n}{d_1}\right)^s=
\sum_{\stackrel{d_1|n}{\gcd(r,d_1)=1}}\sum_{d_2| \tfrac{n}{d_1}}\mu(d_2) \sigma_s\left(\frac{n}{d_1d_2}\right)\,.
\label{eq:proving}
\end{equation}
For each $d|n$, the coefficient of $\sigma_s(n/d)$ in the final expression of \eqref{eq:proving} is
\begin{equation}
\sum_{\stackrel{d_1|d}{\gcd(r,d_1)=1}}\mu\left(\frac{d}{d_1}\right)
=\begin{cases}\mu(d), & \text{if }d|r;\\0, & \text{otherwise,} \end{cases}\label{eq:cases}
\end{equation}
where to obtain the r.h.s. in the above identity we have also used the classical identity $\sum_{d\vert n } \mu(d) = 0$ for $n>1$.
This is also the coefficient of $\sigma_s(n/d)$ on the r.h.s. of relation \eqref{eq:DivId}, which is now proved.

We have shown that with a principal character $\chi_{r,1}$ the Lambert series \eqref{eq:LamChi} reduces to the simple linear combination \eqref{eq:chi11Lam} of the Lambert series already studied in detail in \cite{Dorigoni:2020oon}.
In particular, for {\em odd} positive integer  $s$,  and any modulus $r>1$, we have
\begin{equation}
\left(q\frac{d}{dq}\right)^s\mathcal{L}_s(\chi_{r,1};q) = \sum_{d\vert r}\mu(d) G_{s+1}(q^d)\,,\label{eq:LtoG}
\end{equation}
where we introduce the standard holomorphic Eisenstein series,
\begin{equation}
G_{2k}(q) \coloneqq -\frac{B_{2k}}{4k} + \sum_{n=1}^\infty  \sigma_{2k-1}(n) q^n\,.\label{eq:G2k}
\end{equation}
We remark that $G_{2k}(q)$ is a modular form for $k>1$ and that the Bernoulli number $B_{2k}$
in \eqref{eq:G2k} does not appear in the linear combination \eqref{eq:LtoG}, since $\sum_{d|r}\mu(d)=0$ for $r>1$. 

\subsection{Generalized Lambert series with a pair of characters} 

We now generalize the Lambert series \eqref{eq:LamChi}
to include a {\em pair} of characters, by defining
\begin{equation}
\mathcal{L}_s(\chi_{r_1},\chi_{r_2};q) \coloneqq \sum_{n_1=1}^\infty\sum_{n_2=1}^\infty
\frac{\chi_{r_1}(n_1)}{n_1^s} \chi_{r_2}(n_2) q^{n_1n_2}\,.\label{eq:L2chi}
\end{equation}
Then, by definition, we have
\begin{equation}
\mathcal{L}_s(\chi_{r};q)=\mathcal{L}_s(\chi_{r},\chi_{1,1};q)\,,\label{eq:Lam1} 
\end{equation}
when the second character is set to unity. Setting the {\em first} character to unity, 
we obtain a new generalized Lambert series with a single character, namely
\begin{equation}
\tilde{\mathcal{L}}_s(\chi_r;q)\coloneqq\mathcal{L}_s(\chi_{1,1},\chi_r;q)
=\sum_{n_1=1}^\infty\sum_{n_2=1}^\infty\chi_r(n_2)\frac{q^{n_1n_2}}{n_1^s}
=\sum_{n=1}^\infty \chi_r(n) {\rm Li}_s(q^{n})\,,\label{eq:Lam2}
\end{equation}
where ${\rm Li}_s(z)\coloneqq\sum_{n=1}^\infty z^n/n^s$, with $|z|<1$, is a polylogarithm.

We may also evaluate $\tilde{\mathcal{L}}_s(\chi_r;q)$ as the single sum
\begin{equation}
\tilde{\mathcal{L}}_s(\chi_r;q) = \sum_{n=1}^\infty \frac{\sigma_{s,\chi_r}(n)}{n^s} q^n\,,\label{eq:qser2}
\end{equation}
with an alternative divisor function $\sigma_{s,\chi_r}(n)$ defined by
\begin{equation}
\sigma_{s,\chi_r}(n)  \coloneqq \sum_{d|n} \chi_r(d) d^s\,. \label{eq:Twsigma1}
\end{equation}

In a notable case with $s=2$ and $\chi_r=\chi_{6,5}$, encountered by Bloch and Vanhove \cite{BV}
in the evaluation of a two-loop Feynman integral in quantum field theory, this expansion
over dilogarithms, multiplied  by a character modulo $6$, was referred to as an {\em elliptic} dilogarithm.

Our intention in \eqref{eq:L2chi}, with a pair of characters, is to construct a bridge that interpolates between
work in \cite{Dorigoni:2020oon} and \cite{BV}. It leads to doubly-twisted divisor sums in
\begin{equation}
\mathcal{L}_s(\chi_{r_1},\chi_{r_2};q) = 
\sum_{n=1}^\infty \frac{\sigma_s(\chi_{r_1},\chi_{r_2};n)}{n^s} q^n\,,
\end{equation}
with
\begin{equation}
\sigma_s(\chi_{r_1},\chi_{r_2};n) \coloneqq 
\sum_{d|n} \chi_{r_1}\left(\frac{n}{d}\right) \chi_{r_2}(d) d^{s}
=n^s\sigma_{-s}(\chi_{r_2},\chi_{r_1};n)\,.\label{eq:div2chi}
\end{equation}

\subsection{Generalization to series with two parameters}

With two characters, $\chi_{r_1}$ and $\chi_{r_2}$, it seems natural to have two parameters, $s_1$ and $s_2$.
Accordingly, we define the two-parameter series
\begin{equation}
\Xi_{s_1,s_2}(\chi_{r_1},\chi_{r_2}; q) \coloneqq  \sum_{n_1=1}^\infty \sum_{n_2=1}^\infty 
\frac{\chi_{r_1}(n_1)}{n_1^{s_1}} \frac{\chi_{r_2}(n_2)}{n_2^{s_2}} q^{n_1 n_2} 
= \Xi_{s_2,s_1}(\chi_{r_2},\chi_{r_1};q)\,,\label{eq:LamGen}
\end{equation}
as a further generalization of the one-parameter cases
\begin{equation}
\mathcal{L}_s(\chi_{r_1},\chi_{r_2};q) = \Xi_{s,0}(\chi_{r_1},\chi_{r_2};q) 
= \Xi_{0,s}(\chi_{r_2},\chi_{r_1};q)\,,\label{eq:GenToLam}
\end{equation}
and the single-character cases
\begin{align}
\mathcal{L}_s(\chi_{r};q) &= \Xi_{s,0}(\chi_{r},\chi_{1,1};q)
= \Xi_{0,s}(\chi_{1,1},\chi_{r};q)\,,\label{eq:GenToL}\\
\tilde{\mathcal{L}}_s(\chi_{r};q)& = \Xi_{0,s}(\chi_{r},\chi_{1,1};q) 
= \Xi_{s,0}(\chi_{1,1},\chi_{r};q)\,.\label{eq:GenToLp}
\end{align}

The doubly-twisted divisor sums in \eqref{eq:div2chi} determine the $q$-series 
\begin{equation}
\Xi_{s_1,s_2}(\chi_{r_1},\chi_{r_2};q)= 
\sum_{n=1}^\infty \frac{\sigma_{s_1-s_2}(\chi_{r_1},\chi_{r_2};n)}{n^{s_1}} q^n
=  \sum_{n=1}^\infty \frac{\sigma_{s_2-s_1}(\chi_{r_2},\chi_{r_1};n)}{n^{s_2}} q^n\,.\label{eq:LamGenQ}
\end{equation}

\section{Asymptotic expansion near $q=1^-$}
\label{sec:Asy}

Throughout this paper we will use interchangeably the variables $q$, $y$ and $\tau$,
related by $q=e^{-2\pi y}= e^{2\pi i \tau}$, with ${\rm Re}(y)>0$ and a modular parameter $\tau=iy$, 
with  ${\rm Im}(\tau)>0$. By slight abuse of notation, we will write
\begin{equation}
\Xi_{s_1,s_2}(\chi_{r_1},\chi_{r_2} ; q) = \Xi_{s_1,s_2}(\chi_{r_1},\chi_{r_2} ; y)=\Xi_{s_1,s_2}(\chi_{r_1},\chi_{r_2} ; \tau) \quad {\rm with}\quad  q=e^{-2\pi y} = e^{2\pi i \tau}\,.
\end{equation}

We now analyze the asymptotic expansion of $\Xi_{s_1,s_2}(\chi_{r_1},\chi_{r_2} ; y)$
for $y\to 0^+$. 
It is convenient to rewrite \eqref{eq:LamGen} as
\begin{equation}
\Xi_{s_1,s_2}(\chi_{r_1},\chi_{r_2} ; q) =  \sum_{n=1}^\infty \frac{\chi_{r_1}(n)}{n^{s_1}} \Phi_{s_2}(\chi_{r_2} ; q^n)= \sum_{n=1}^\infty \frac{\chi_{r_2}(n)}{n^{s_2}} \Phi_{s_1}(\chi_{r_1} ; q^n)\,\label{eq:Repeated1}
\end{equation}
with an auxiliary function
\begin{equation}
\Phi_{s}(\chi_r ; q) \coloneqq \sum_{n=1}^\infty \frac{\chi_r(n)}{n^{s}} q^n\,.\label{eq:Phi}
\end{equation}
At $s=0$, this function takes the particularly simple form 
\begin{align}
\Phi_{0}(\chi_r ; q) &\notag\coloneqq \sum_{n= 1}^\infty \chi_r(n) q^n = \sum_{k=0}^\infty q^{kr} \sum_{a\in (\mathbb{Z}/r\mathbb{Z})^\times}  \chi_r(a) q^a\\
& \phantom{:}= \frac{ \sum_{a\in (\mathbb{Z}/r\mathbb{Z})^\times}  \chi_r(a) q^a}{1-q^r}\,,\label{eq:NPPhi}
\end{align}
where we wrote $n=kr+a$ and used the $r$-periodicity of the Dirichlet character combined with the fact that $\chi_r(a) \neq 0$ only when $\gcd(a,r)=1$.

We note that the generalized Lambert series with a pair of characters \eqref{eq:L2chi} can be expressed as
\begin{equation}
\mathcal{L}_s(\chi_{r_1},\chi_{r_2};q) = \sum_{n=1}^\infty \frac{\chi_{r_1}(n)}{n_1^s} \Phi_0(\chi_{r_2};q^n)\,.\label{eq:GenLamPhi}
\end{equation}
When $\chi_{r_2}$ is set to unity, we have that
\begin{equation}
 \Phi_0(\chi_{1,1};q) = \frac{q}{1-q}\,,
 \end{equation}
 so that \eqref{eq:GenLamPhi} reduces to the twisted Lambert series \eqref{eq:LamChi}. Conversely, when $\chi_{r_1}$ is set to unity equation \eqref{eq:GenLamPhi} reduces to the generalized Lambert series \eqref{eq:Lam2}.

The analysis of the asymptotic expansion of \eqref{eq:Repeated1} as $y\to 0^+$ depends  crucially on the Dirichlet characters under consideration.
For the particular case when either of the characters reduces to the trivial one, say $\chi_{r_2} = \chi_{1,1}$, the asymptotic expansion of \eqref{eq:Repeated1} can be neatly derived using a method presented by Zagier in \cite{Zagier} and
reviewed in Appendix \ref{sec:AsyApp}. 
Here we consider the case where the characters $\chi_{r_1}$ and $\chi_{r_2}$, as well as the parameters $s_1,s_2\in \mathbb{C}$, are generic and we will comment on particular instances thereof when appropriate.

To proceed, we first make use of a Mellin-Barnes integral representation for \eqref{eq:Phi},
\begin{equation}
\Phi_{s}(\chi_r ; y)  = \int_{\gamma-i \infty}^{\gamma+i\infty} \Gamma(t) L(\chi_r, t+s) (2\pi y)^{-t} \,\frac{{\rm d}t}{2\pi i}\,,\label{eq:Mellin}
\end{equation}
where $L(\chi_r,t+s)$ is a Dirichlet $L$-function and
$\gamma\in \mathbb{R}$ is chosen such that $\gamma >0$, for $\chi_r$ non-principal, and 
$\gamma>{\rm max}(0,1-\rm{Re}(s))$, for $\chi_r$ principal. 
With this choice of integration contour, we can replace the $L$-function by its convergent series
$L(\chi_r,t+s)=\sum_{n=1}^\infty \chi_r(n) n^{-s-t}$ and then use the dominated convergence theorem 
to exchange the order of summation and contour integration, reproducing \eqref{eq:Phi}. We refer to Appendix \ref{sec:DirApp} for key definitions and properties of Dirichlet characters and their associated $L$-functions.

Using this integral representation, we derive a Mellin-Barnes integral representation for \eqref{eq:Repeated1},
\begin{equation}
\Xi_{s_1,s_2}(\chi_{r_1},\chi_{r_2};y) = \int_{\gamma-i \infty}^{\gamma+i\infty} \Gamma(t) L(\chi_{r_1}, t+s_1)L(\chi_{r_2}, t+s_2) (2\pi y)^{-t} \,\frac{{\rm d}t}{2\pi i}\,.\label{eq:Mellin2}
\end{equation}
Importantly, the contour of integration depends on the types of character considered.
To obtain the small-$y$ expansion of \eqref{eq:Repeated1} we  close the contour of integration in \eqref{eq:Mellin2} to the left semi-half plane ${\rm Re}(t)<0$ and collect the residues coming from the poles of the integrand.
Starting with the case where both $\chi_{r_1}$ and $\chi_{r_2}$ are  non-principal characters we have that their corresponding $L$-functions, $L(\chi_r,t)$, are entire in $t$. The contour of integration in \eqref{eq:Mellin2} then requires $\gamma>0$, so that the only poles in \eqref{eq:Mellin2} are the poles of the gamma function $\Gamma(t)$. 
In this case we obtain
\begin{equation}
\Xi_{s_1,s_2}(\chi_{r_1},\chi_{r_2};y)  \sim \Xi^{\rm Pert}_{s_1,s_2}(\chi_{r_1},\chi_{r_2};y)=  
\sum_{k=0}^\infty \frac{(-2\pi y)^k}{k!} L(\chi_{r_1},s_1-k) L(\chi_{r_2},s_2-k)\,,\label{eq:XiAsyPgen}
\end{equation}
where we indicate by $\sim$ that $\Xi_{s_1,s_2}$ is asymptotic in the sense of Poincar\'e to
$\Xi^{\rm Pert}_{s_1,s_2}$, a formal perturbative asymptotic expansion in powers of $y$, 
postponing consideration of resummation and exponentially suppressed terms.

We remind the reader that a function $F(y)$ is asymptotic in the sense of Poincar\'e  to the formal power series $F_P(y)= \sum_{n=0}^\infty a_n y^n $ as $y\to 0^+$, if
$ F(y) - \sum_{n=0}^N a_N y^n = O(y^{N+1})$ as $y\to 0$ for all $N\in \mathbb{N}$. Crucially, functions which differ by exponentially suppressed, non-perturbative terms such as $F(y)$ and $F(y)+ e^{-1/y}$ are asymptotic in the sense of Poincar\'e to the same formal power series $F_P(y)$ but they are manifestly different functions. The aim of the next section is precisely to derive a more refined version of the asymptotic expansion \eqref{eq:XiAsyPgen}, referred to as the transseries expansion and which incorporates such crucial non-perturbative terms.
 
When one of the characters reduces to the trivial or to a principal character, the story is slightly different.
First, given the symmetry \eqref{eq:LamGen} under exchange $s_1\leftrightarrow s_2 $ and $\chi_{r_1}\leftrightarrow \chi_{r_2}$, we assume that the character $\chi_{r_1}$ is non-principal while $\chi_{r_2}=\chi_{r_2,1}$ since as previously remarked when both characters are either principal or trivial the story is a direct consequence of the results discussed in~\cite{Dorigoni:2020oon}. 

In this case the contour of integration in \eqref{eq:Mellin2} requires $\gamma>{\rm max}(0,1-\rm{Re}(s_2))$ and the integrand has poles located at $t=-k$ with $k\in \mathbb{N}$ as well as a single pole located at $t=1-s_2$ coming from the $L$-function $L(\chi_{r_2,1},t+s_2) $.
Since every principal character $\chi_{r,1}$ is induced by the trivial character, we use the relation between $L$-functions given in \eqref{eq:Lnonprim} to derive the identity
\begin{equation}
L(\chi_{r,1},s) = \zeta(s) \prod_{p\vert r} \left(1-p^{-s}\right)\,,
\end{equation}
where the product runs over the prime factors $p$ of $r$ and $\zeta(s)=L(\chi_{1,1},s)$
is the Riemann zeta function, with $\zeta(s)=\sum_{n=1}^\infty n^{-s}$ for $\rm{Re}(s)>1$.
We then deduce that the residue at $s=1$ of $L(\chi_{r,1},s) $ is 
\begin{equation}
{\rm res}_{s=1} L(\chi_{r,1},s)  = \frac{\phi(r)}{r}\,,
\end{equation}
with $\phi(r)$ the Euler totient function.
Assuming that $s_2\in \mathbb{C}$ is generic, we obtain
\begin{align}
& \Xi_{s_1,s_2}(\chi_{r_1},\chi_{r_2,1};y) \sim  \Xi^{\rm Pert}_{s_1,s_2}(\chi_{r_1},\chi_{r_2,1};y)=\nonumber  \\
&(2\pi y)^{s_2-1}\Gamma(1-s_2)  L(\chi_{r_1},s_1+1-s_2) \frac{\phi(r_2)}{r_2}+ \sum_{k=0}^\infty \frac{(-2\pi y)^k}{k!}  L(\chi_{r_1},s_1-k)L(\chi_{r_2,1},s_2-k)\,.\label{eq:XiAsyPng}
\end{align} 

We note that when $s_2\in \mathbb{N}$, the expression above has to be understood as a limit since in this case the first term as well as the term with $k=s_2-1$ both become singular. This is due to the fact that for $s_2\in \mathbb{N}$ the integrand of \eqref{eq:Mellin2} develops a double pole for $t=1-s_2$. However, if we take in \eqref{eq:XiAsyPng} the limit  for $s_2\to m\in \mathbb{N}$ we find the correct and finite result for the residue,
\begin{align}
&\notag \lim_{s_2\to m} \left[ (2\pi y)^{s_2-1}\Gamma(1-s_2) L(\chi_{r_1},s_1{+}1{-}s_2)\frac{\phi(r_2)}{r_2}{+} \frac{(-2\pi y)^{m-1}}{(m-1)!} L(\chi_{r_1},s_1{+}1{-}m)L(\chi_{r_2,1},s_2{+}1{-}m) \right]\\
& = \left[ \frac{L'(\chi_{r_1} , s_1{+}1{-}m)}{L(\chi_{r_1} , s_1{+}1{-}m)} + \left(\sum_{p\vert r_2} \frac{\log p}{p-1}\right) + 
 \gamma+\psi(m) - \log(2\pi y) \right] \label{eq:s2Int}   \frac{(-2\pi y)^{m-1}}{(m-1)!} L(\chi_{r_1} , s_1{+}1{-}m) \frac{\phi(r_2)}{r_2}  \,,
\end{align}
where $\gamma$ is the Euler-Mascheroni constant, $\psi(x) = \Gamma'(x)/\Gamma(x)$ is the digamma function and $L'(\chi_r,s) = d L(\chi_r,s)/ds$. 

In light of the discussion above, we can write the asymptotic perturbative expansion of the function $\Xi_{s_1,s_2}(\chi_{r_1},\chi_{r_2};y)$ for general characters in the unified symmetric form
\begin{align}
& \Xi_{s_1,s_2}(\chi_{r_1},\chi_{r_2};y) \sim  \Xi^{\rm Pert}_{s_1,s_2}(\chi_{r_1},\chi_{r_2};y)=\nonumber \\
& y^{s_2-1} A_{s_1,s_2}(\chi_{r_1},\chi_{r_2}) +y^{s_1-1} A_{s_2,s_1}(\chi_{r_2},\chi_{r_1}) 
+ \sum_{k=0}^\infty \frac{(-2\pi y)^k}{k!}  L(\chi_{r_1},s_1-k)L(\chi_{r_2},s_2-k)\,,\label{eq:XiAsyP}
\end{align} 
where the coefficient $A_{s_1,s_2}(\chi_{r_1},\chi_{r_2})$ is given by
\begin{equation}
A_{s_1,s_2}(\chi_{r_1},\chi_{r_2}) = \begin{cases}
(2\pi)^{s_2-1}\Gamma(1-s_2) L(\chi_{r_1},s_1+1-s_2)\frac{\phi(r_2)}{r_2}\,,&\text{if }\chi_{r_2} = \chi_{r_2,1}\,;\\
0\,,&\text{otherwise.}\end{cases}
\end{equation}
We stress again that instances where one of the characters is principal and $s_1,s_2 \in \mathbb{N}$ must be understood as limiting cases of \eqref{eq:XiAsyP}, as explained in \eqref{eq:s2Int}.

From \eqref{eq:XiAsyP} it is immediate to derive the asymptotic expansions for \eqref{eq:L2chi},
\begin{align}
&\notag \mathcal{L}_s(\chi_{r_1},\chi_{r_2};y) = \Xi_{s,0}(\chi_{r_1},\chi_{r_2};y)\sim\\
& \mathcal{L}^{\rm Pert}_s(\chi_{r_1},\chi_{r_2};y) =  \sum_{k=0}^\infty \frac{(-2\pi y)^k}{k!} L(\chi_{r_1},s-k) L(\chi_{r_2},-k)\,,\label{eq:AsyL2chi}
\end{align}
when neither of the characters is principal.

When $\chi_{r_2}$ is the trivial character, as in \eqref{eq:LamChi}, the above equation must be replaced by
\begin{align}
 & \mathcal{L}_s(\chi_r ; y) =\mathcal{L}_s(\chi_{r},\chi_{1,1};y)= \Xi_{s,0}(\chi_r,\chi_{1,1};y) \sim \nonumber \\
 & \mathcal{L}^{\rm Pert}_s(\chi_r ; y) =  \sum_{k=-1}^\infty \frac{(-2\pi y)^k}{k!}  L(\chi_r,s-k) \zeta(-k)\,,\label{eq:AsyL}
\end{align}
where we note that in this particular case the first term in \eqref{eq:XiAsyPng} can be understood as the limit for $k\to-1$ of the $k^{th}$ summand.

When $\chi_{r_1}$ is the trivial character, as in \eqref{eq:Lam2}, we obtain
\begin{align}
&\notag \tilde{\mathcal{L}}_s(\chi_r;y)  =\mathcal{L}_s(\chi_{1,1},\chi_{r};y)= \Xi_{0,s}(\chi_{r},\chi_{1,1};y) \sim \\
 &\tilde{\mathcal{L}}^{\rm Pert}_s(\chi_r;y)=  (2\pi y)^{s-1}\Gamma(1-s) L(\chi_r,1-s)+ \sum_{k=0}^\infty \frac{(-2\pi y)^k}{k!} \zeta(s-k) L(\chi_r,-k)\,,\label{eq:LtP}
\end{align}
where again we note that for $s\in \mathbb{N}$ the sum of first term and the summand with $k=s-1$ must be understood as the finite limit \eqref{eq:s2Int}.

Given the general asymptotic expansion \eqref{eq:XiAsyP}, we immediately notice that something peculiar
happens for parameters $s_1,s_2 \in\mathbb{Z}$.
As reviewed in \eqref{eq:TZ}, similarly to the well-known trivial zeroes of the Riemann zeta function, $\zeta(-2n)$ with $n\in \mathbb{N}$ and $n>0$, we have that the $L$-function $L(\chi_r,s)$ has trivial zeroes located at
\begin{equation}
L(\chi_r, -\kappa -2n)=0 \,, \qquad \forall n \in \mathbb{N}\,,\label{eq:TrivZer}
\end{equation}
where $\chi_r(-1) = (-1)^\kappa$ and $\kappa\in \{0,1\}$, and $n>0$ when $\chi_r$ is a principal character.
Let us now assume that $s_1,s_2 \in\mathbb{Z}$.  For sufficiently large values of $k$, i.e. $k\geq {\rm max}(0,s_1) + {\rm max}(0,s_2)+1$, the $k^{th}$ summand appearing in the general asymptotic expansion \eqref{eq:XiAsyP} is proportional to the product of $L$-values:
\begin{equation}
L(\chi_{r_1}, - n_1 ) L( \chi_{r_2}, - n_2 )\,,\label{eq:LvalProd}
\end{equation}
where importantly  $n_1=k-s_1,n_2=k-s_2$ are positive integers such that  $n_1 +n_2 = 2k-s_1- s_2$.
Given the location \eqref{eq:TrivZer} for the trivial zeroes of a Dirichlet $L$-function, we must have that both $n_1 \equiv \kappa_1+1 \text{ mod }2$ as well as  $n_2 \equiv \kappa_2+1 \text{ mod }2$ for the product \eqref{eq:LvalProd} \textit{not} to vanish. Furthermore, since we also have $n_1 +n_2\equiv s_1+s_2 \text{ mod }2$ we deduce that \eqref{eq:LvalProd} and hence the $k^{th}$ summand appearing in the general asymptotic expansion \eqref{eq:XiAsyP} can be non-vanishing for arbitrarily large $k$ only when $s_1+s_2 \equiv \kappa_1 +\kappa_2 \text{ mod }2$.
Conversely, we conclude that the asymptotic expansion \eqref{eq:XiAsyP} {\em terminates} after finitely many terms when
\begin{equation}
s_1,s_2 \in\mathbb{Z}\,,\qquad s_1+s_2 \equiv \kappa_1 +\kappa_2 +1 \text{ mod }2\,,\label{eq:CondTrunc}
\end{equation}
where $\chi_{r_i}(-1) = (-1)^{\kappa_i}$ and $\kappa_i\in\{0,1\}$ denotes the even/odd parity of the given character.

For example if we consider the special case \eqref{eq:AsyL} with $s\in \mathbb{N}$ we see that the summand 
$$
\zeta(-k) L(\chi_r,s-k)
$$ 
vanishes for all $k> s+\kappa-1$, with $\chi_{r}(-1)=(-1)^\kappa$, in each of the following cases: 
\begin{itemize}
\item the character $\chi_r$ is even and $s\in \mathbb{N}$ is odd, 
\item the character $\chi_r$ is odd and $s\in \mathbb{N}$ is even. 
\end{itemize}

Focusing for now on the case where the asymptotic expansion \eqref{eq:AsyL} terminates after finitely many terms,~i.e.\ $s$ even and $\chi_r$ odd, or $s$ odd and $\chi_r$ even, we find that the corresponding Lambert series $\mathcal{L}_s(\chi_r;y)$ possesses a terminating asymptotic expansion given by the Laurent polynomial:
\begin{equation}
\mathcal{L}_s(\chi_r ; y) \sim \mathcal{L}^{{\rm Pert}}_s(\chi_r ; y) = \sum_{k=-1}^{s+\kappa-1} \frac{(-2\pi y)^k}{k!} \zeta(-k) L(\chi_r,s-k)\,.\label{eq:Asypert}
\end{equation}
Note that in the above Laurent polynomial the term with $k=s+\kappa-1$ actually vanishes for all $s\in \mathbb{N}$ with $s>0$ as a consequence of \eqref{eq:TZ}. In the special case with $s=0$ and  $\chi_r$ odd,
the above expression receives contribution from both $k=-1$ and $k=0$, giving
\begin{equation}
\mathcal{L}_0(\chi_r ; y) \sim \mathcal{L}^{{\rm Pert}}_0(\chi_r ; y) = \frac{L(\chi_r,1)}{2\pi y} - \frac{L(\chi_r,0)}{2}\,.\label{eq:AsypertS0}
\end{equation}
We remark that $\mathcal{L}_0(\chi_r;\tau)$ with $\chi_r$ an odd primitive character corresponds, modulo a constant term, 
to the twisted Eisenstein series $G_1(\chi_r;\tau)$ defined in \eqref{eq:TwEisen2}. 
We will discuss relations to Eisenstein series at a later stage in Section \ref{sec:ModularPrim}.

A corresponding result holds for the asymptotic expansion \eqref{eq:LtP}, with $s\in \mathbb{N}$, 
where the summand 
$$
\zeta(s-k) L(\chi_r,-k)
$$ 
vanishes for all $k> s+\kappa$ when $\chi_r$ and $s$ have opposite parities.
In such cases the asymptotic expansion \eqref{eq:LtP} terminates after finitely many terms, giving 
\begin{align}
 \tilde{\mathcal{L}}_s(\chi_r ; y) \sim {\tilde{\mathcal{L}}}^{{\rm Pert}}_s(\chi_r ; y)&\notag = 
\left[ \frac{L'(\chi_r ,1-s)}{L(\chi_r , 1-s)}  +  \gamma+\psi(s) - \log(2\pi y) \right]  
\frac{(-2\pi y)^{s-1}}{(s-1)!} L(\chi_r , 1-s) \\
&\phantom{=}+  \sum_{\substack{k=0 \\ k\neq s-1}}^{s+\kappa} \frac{(-2\pi y)^k}{k!} \zeta(s-k) L(\chi_r,-k)\,,\label{eq:Asypert2}
\end{align}
where we have used the limit \eqref{eq:s2Int} to obtain a regular expression.

As previously remarked, the terminating series for $\mathcal{L}_s(\chi_r ; y) \sim \mathcal{L}^{{\rm Pert}}_s(\chi_r ; y)$, in \eqref{eq:Asypert}, 
and for ${\tilde{\mathcal{L}}}_s(\chi_r ; y) \sim \mathcal{L}^{{\rm Pert}}_s(\chi_r ; y)$, in \eqref{eq:Asypert2}, 
hold  in the sense of Poincar\'e. Yet the functions $\tilde{\mathcal{L}}_s(\chi_r;y) $ and ${\mathcal{L}}_s(\chi_r;y)$
cannot have such simple expressions, as Laurent polynomials in $y$ with a possible addition of a logarithmic term. 
The fact that for, $s\in \mathbb{N}$, these perturbative expansions near $y\to0^+$ terminate after finitely many terms  
indicates that our asymptotic analysis is incomplete and that we are missing fundamental non-perturbative corrections,
which are exponentially suppressed. 

We now proceed to resum the general asymptotic expansion \eqref{eq:XiAsyP} using resurgence analysis and 
construct the missing non-perturbative corrections which will then be specialized to the particular cases 
$\mathcal{L}_s(\chi_r ; y)$ and ${\tilde{\mathcal{L}}}_s(\chi_r ; y)$.

\subsection{Transseries for primitive characters}
\label{sec:NPPrim}

In this section, we analyze the non-perturbative resummation of the formal asymptotic power series expansion \eqref{eq:XiAsyP} under the assumption that each $\chi_r$ is a primitive character with modulus $r>1$. 

We start by rewriting the asymptotic expansion \eqref{eq:XiAsyP} making use of the functional equation \eqref{eq:LFunct}, valid only for primitive characters, to arrive at
\begin{align}
&\notag \Xi_{s_1,s_2}(\chi_{r_1},\chi_{r_2} ; y)\sim  \Xi^{\rm Pert}_{s_1,s_2}(\chi_{r_1},\chi_{r_2} ; y)
= \\
&\notag{} \Xi^{{\rm LP}}_{s_1,s_2}(\chi_{r_1},\chi_{r_2};y){+}   \frac{y^{s_1+s_2-1} r_1^{s_2-\frac{1}{2} } r_2^{s_1-\frac{1}{2} }  \epsilon(\chi_{r_1})\epsilon(\chi_{r_2})}{\pi} \sum_{k=0}^\infty2(-1)^k \! \sin\left(\tfrac{(k-s_1-\kappa_1)\pi}{2} \right) \! \sin\left(\tfrac{(k-s_2-\kappa_2) \pi}{2} \right) \\
&\phantom{\sim}\times \frac{\Gamma(k+1-s_1)\Gamma(k+1-s_2) }{k!}\left(\frac{r_1 r_2 y}{2\pi}\right)^{k+1-s_1 -s_2} L(\bar{\chi}_{r_1}, k+1-s_1) L(\bar{\chi}_{r_2}, k+1-s_2)\,.\label{eq:AsyXiSt}
\end{align}
To avoid longer expressions we denote by $\Xi^{{\rm LP}}_{s_1,s_2}(\chi_{r_1},\chi_{r_2};y)$ the first two isolated monomials appearing in \eqref{eq:XiAsyP} which may be present  when at least one of the characters is a principal character. For a general primitive character we define,
\begin{equation}
\epsilon(\chi_r) \coloneqq \frac{ \sum_{n=1}^r \chi_r(n) e^{\frac{2\pi i n}{r}}}{i^\kappa \sqrt{r}} \,,\label{eq:GaussIT}
\end{equation}
which is a complex number on the unit circle, i.e. $|\epsilon(\chi_r)|=1$. If  $\chi_r$ is real, then $\epsilon(\chi_r)=1$.

For generic $(s_1,s_2)\in \mathbb{C}^2$ the asymptotic expansion \eqref{eq:AsyXiSt} is a 
factorially divergent power series. As already noted,  in cases with $(s_1,s_2)\in \mathbb{Z}$ 
and  $s_1+s_2 \equiv \kappa_1 +\kappa_2+1\text{ mod }2$, the perturbative expansion
terminates, after finitely many terms. This is accounted for  by the sine functions in \eqref{eq:AsyXiSt}.
In particular for $\chi_{r_1}=\chi_r$, a primitive non-trivial character, and $\chi_{r_2}=\chi_{1,1}$,
the trivial character, we find that when $ (s_1,s_2)= (s,0)$ or $(s_1,s_2) = (0,s)$, with $s\in \mathbb{N}$ and 
of {\em opposite} parity to the character $\chi_r$, the series {\em terminates}.
 
To obtain the non-perturbative corrections we assume that the parameters $(s_1,s_2)\in \mathbb{C}^2$ are generic and 
specialize to particular cases only at the very end.
First, we rewrite the trigonometric factors in \eqref{eq:AsyXiSt}, using the identity
\begin{align}
&  2(-1)^k\sin\!\left(\tfrac{\!(k-s_1-\kappa_1)\pi}{2}\right)\!\sin\left(\!\tfrac{ (k-s_2-\kappa_2) \pi}{2}\right)=
\sin\!\left(\!\tfrac{(s_1+s_2+\kappa_1+\kappa_2-1)\pi}{2}\right)
-(-1)^k \sin\!\left(\!\tfrac{(s_1-s_2+\kappa_1-\kappa_2-1)\pi}{2}\right)
\end{align}
and then rewrite both $L$-functions in terms of their Dirichlet series, i.e.
\begin{equation}
L(\bar{\chi}_{r_1}, k+1-s_1) = \sum_{n_1=1}^\infty \frac{\bar{\chi}_{r_1}(n_1)}{n_1^{k+1-s_1}} \,, \qquad L(\bar{\chi}_{r_2}, k+1-s_2) = \sum_{n_2=1}^\infty \frac{\bar{\chi}_{r_2}(n_2)}{n_2^{k+1-s_2}} \,.
\end{equation}
The asymptotic perturbative expansion in \eqref{eq:AsyXiSt} can then be rearranged as
\begin{align}
&\notag \Xi_{s_1,s_2}(\chi_{r_1},\chi_{r_2} ; y)\sim  \Xi^{\rm Pert}_{s_1,s_2}(\chi_{r_1},\chi_{r_2} ; y)= \\
 &\Xi^{{\rm LP}}_{s_1,s_2}(\chi_{r_1},\chi_{r_2};y){+} \frac{ y^{s_1+s_2-1} r_1^{s_2-\frac{1}{2}} r_2^{s_1-\frac{1}{2}} \epsilon(\chi_{r_1})\epsilon(\chi_{r_2})}{\pi} \!\sum_{n_1,n_2=1}^\infty \! \!\frac{\bar{\chi}_{r_1}(n_1)}{n_1^{s_2}} \frac{\bar{\chi}_{r_2}(n_2)}{n_2^{s_1}}  F_{s_1,s_2,\kappa_1,\kappa_2}\left(\! \frac{ 2\pi n_1 n_2 }{r_1 r_2y} \right),
\label{eq:XiSeries}  
\end{align}
where the auxiliary function $F_{s_1,s_2,\kappa_1,\kappa_2}(z)$ denotes the formal power series, for $|z|\gg1$,
\begin{align}
 F_{s_1,s_2,\kappa_1,\kappa_2}(z) =  & \notag z^{s_1+s_2} \sum_{k=0}^\infty \frac{\Gamma(k+1-s_1)\Gamma(k+1-s_2)} {k!} \\
&\times   \left[ \sin\left( \tfrac{(s_1+s_2+\kappa_1+\kappa_2-1)\pi}{2}\right)   -(-1)^k  \sin\left( \tfrac{(s_1-s_2+\kappa_1-\kappa_2-1)\pi}{2}\right)\right] z^{-k-1}\,.\label{eq:Faux} 
\end{align}

The task at hand  is to define an unambiguous resummation of the formal factorially divergent power series \eqref{eq:Faux} 
for real $z>0$ when the parameters $(s_1,s_2) $ are real numbers. This goal can be achieved via Borel-\'Ecalle resummation~\cite{Ecalle:1981,delabaere1999resurgent} which we now briefly present. (See also the review \cite{Dorigoni:2014hea}.)

Given a formal asymptotic series for $|z|\gg1$
\begin{equation}\label{eq:Formal}
F(z) = \sum_{k=0}^\infty  c_k\,  z^{-k-1}\,,
\end{equation}
with factorially divergent perturbative coefficients $c_k$, we can construct the associated Borel transform,
\begin{align}
B(t) = \sum_{k=0}^\infty \frac{c_k}{k!}\, t^k\,.
\end{align}
When the series defining the Borel transform has finite radius of convergence we use the identity
\begin{align}
 z^{-k-1} \Gamma(k+1) = \int_0^\infty  e^{-t}  \left( \frac{t}{z} \right)^{k }   \frac{{\rm d}t}{z}\,,
\end{align}
to obtain an analytic continuation of the formal asymptotic power series $F(z)$.
This can be achieved via the directional Borel resummation 
\begin{equation}
\mathcal{S}_\theta\left[ F\right](z) \coloneqq \int_0^\infty e^{-t} B\left(\frac{t}{z} \right) \frac{{\rm d}t}{z} = \int_0^{e^{i\theta}\infty} e^{-t\,z } B(t) {\rm d}t\,,\label{eq:LateralDef}
\end{equation}
where $\theta =-\mbox{arg}\,z$, provided that the direction of integration does not contain any singularity of the Borel transform, i.e.\ provided that the line ${\rm arg}\,t= \theta$ is not a {\em Stokes direction} for $B(t)$.

We now apply this method to the formal power series in \eqref{eq:Faux}.
First, we obtain the Borel transform of the sum over $k$,
\begin{align}
\notag B_{s_1,s_2,\kappa_1,\kappa_2}(t) 
&=\notag  \Gamma(1-s_1) \Gamma(1-s_2)\left[ \sin\left( \tfrac{(s_1+s_2+\kappa_1+\kappa_2-1)\pi}{2}\right)
 \,_2F_1(1-s_1,1-s_2,1 \vert t)\right. \\
&\phantom{=}\left.{}-  \sin\left( \tfrac{(s_1-s_2+\kappa_1-\kappa_2-1)\pi}{2}\right) \,_2F_1(1-s_1,1-s_2,1 \vert -t)\right]\,,
\label{eq:TwoSines}
\end{align}
with $\,_2F_1(a,b,c\vert t)$ denoting the standard hypergeometric function.
The directional Borel resummation for the formal series \eqref{eq:Faux} is then given by the directional Laplace transform
\begin{equation}\label{eq:DirBorel} 
\mathcal{S}_\theta\left[ F_{s_1,s_2,\kappa_1,\kappa_2}\right](z) =  z^{s_1+s_2} \int_0^{e^{i\theta} \infty} e^{ - t \, z} \,B_{s_1,s_2,\kappa_1,\kappa_2}(t) \,{\rm d}t\,.
\end{equation}

Recalling that $z=(2\pi n_1 n_2)/(r_1 r_2y)$ in the summand of \eqref{eq:XiSeries}, we see that the direction specified in
\eqref{eq:DirBorel} is $\theta = \mbox{arg}\,y= -\mbox{arg}\,z$. Yet we cannot set $\theta=0$, for real $y>0$,
since the first hypergeometric function, ${}_2F_1(1-s_1,1-s_2,1 \vert t)$, in \eqref{eq:TwoSines} has
a branch point at $t=1$. 
The best that we can do is to consider the directional Borel resummation $\mathcal{S}_\theta$ in the two distinct limits  $\theta\to0^+$ and $\theta\to0^-$, which we will denote by $\mathcal{S}_+ \coloneqq \lim_{\theta\to 0^+} \mathcal{S}_\theta$ and $\mathcal{S}_-\coloneqq \lim_{\theta\to 0^-} \mathcal{S}_\theta$.
These two lateral Borel resummations, $\mathcal{S}_\pm [F_{s_1,s_2,\kappa_1,\kappa_2}](z) $, have the same real parts, but imaginary parts that  differ in sign.
The discontinuity across the Stokes direction $\theta=0$ is called a {\em Stokes automorphism} of \eqref{eq:DirBorel}.
Its imaginary part can be computed from the discontinuity of the hypergeometric function across the cut with
a branch point at  $t=1$. For real $t>1$, we have
\begin{align}
&\notag \lim_{\epsilon\to 0^+}\left[\,_2F_1(a,b,c\vert t + i \epsilon) - \,_2F_1(a,b,c\vert t - i \epsilon)\right]=\\
&2\pi i \frac{\Gamma(c)}{\Gamma(a)\Gamma(b)} t^{1-c} (t-1)^{c-a-b} \,_2\tilde{F}_1(1-a,1-b,c-a-b+1\vert 1-t)\,,
\end{align}
with $_2\tilde{F}_1(a,b,c\vert t) \coloneqq \,_2 F_1(a,b,c\vert t)/\Gamma(c)$, for notational convenience.
We then obtain
\begin{align}
& \Big( \mathcal{S}_+ - \mathcal{S}_- \Big) \left[ F_{s_1,s_2,\kappa_1,\kappa_2}\right](z) 
\notag =
\lim_{\theta\to 0^+} \Big[\mathcal{S}_\theta\left[ F_{s_1,s_2,\kappa_1,\kappa_2}\right](z)
- \mathcal{S}_{-\theta}\left[ F_{s_1,s_2,\kappa_1,\kappa_2}\right](z) \Big]= \\
& \label{eq:StokesAuto} 2\pi i  \sin\ \left( \tfrac{(s_1+s_2+\kappa_1+\kappa_2-1)\pi}{2}\right) z^{s_1+s_2} \int_1^\infty e^{-t\,z} (t-1)^{s_1+s_2-1}\,_2\tilde{F}_1(s_1,s_2,s_1+s_2\vert 1-t)\,{\rm d}t  \,.
\end{align}

By translating the integration variable  $t\to t+1$, we immediately see 
that the Stokes automorphism is inherently non-perturbative, since
\begin{align}
&  \Big( \mathcal{S}_+ - \mathcal{S}_- \Big) \left[ F_{s_1,s_2,\kappa_1,\kappa_2}\right](z) = \notag  \\
& \label{eq:StokesAuto2}   2\pi i  \sin\left( \tfrac{(s_1+s_2+\kappa_1+\kappa_2-1)\pi}{2}\right)  z^{s_1+s_2} e^{-z} \int_0^\infty e^{-t\,z} t^{s_1+s_2-1}\,_2\tilde{F}_1(s_1,s_2,s_1+s_2\vert -t)\,{\rm d}t 
\end{align}
exhibits exponential suppression by a factor $\exp(-z)=\exp(-(2\pi n_1 n_2)/(r_1 r_2y))$,
all of whose derivatives w.r.t.\  $y$ vanish as $y\to0^+$.

To obtain a result that is real, for real values of $(s_1,s_2,z)$, we may take
an average of the two lateral resummations, referred to as a {\em median} resummation \cite{delabaere1999resurgent},
namely
\begin{align}
&\mathcal{S}_{\text{med}}\left[ F_{s_1,s_2,\kappa_1,\kappa_2}\right](z)\coloneqq \notag \mathcal{S}_{\pm}\left[F_{s_1,s_2,\kappa_1,\kappa_2} \right](z) \\
&    \mp i \pi  \sin\left( \tfrac{(s_1+s_2+\kappa_1+\kappa_2-1)\pi}{2}\right) z^{s_1+s_2}  e^{-z} \int_0^\infty e^{-t\,z} t^{s_1+s_2-1}\,_2\tilde{F}_1(s_1,s_2,s_1+s_2\vert-t)\,{\rm d}t  \,.\label{eq:SmedInter} 
\end{align}
This may be regarded as an analytic gluing of two different types of formal objects: the resummations $\mathcal{S}_{\pm}\left[F_{s_1,s_2,\kappa_1,\kappa_2} \right](z)$ of the formal power series \eqref{eq:Faux} and an infinite series of non-perturbative, exponentially suppressed corrections encoded by the Stokes automorphism. It is an example of a {\em transseries}~\cite{Edgar}, a concept that generalizes formal power series expansions.

Our story is not yet complete. We have used {\em imaginary} non-perturbative terms  to ensure that the median 
resummation \eqref{eq:SmedInter} has the appropriate reality properties. Yet we have been silent about the possibility of
non-perburative terms that are {\em real}. Here we resort to a conjectural method, attested  to by many
numerical checks, at high precision.

We know that resummation of the perturbative series jumps when we cross the Stokes direction $\theta=0$.
We know that the imaginary part of this jump depends on the fixed parameters $(s_1,s_2)$ 
and the parities $(\kappa_1,\kappa_2)\in\{0,1\}$ of the characters $(\chi_{r_1},\chi_{r_2})$.
We denote the coefficient of the non-perturbative contributions by a transseries parameter $\sigma$, which jumps
when we cross the Stokes direction $\theta=0$. We take a normalization of $\sigma$ from \eqref{eq:SmedInter}, obtaining
\begin{equation}
\mbox{Im} \,\sigma_\pm (\chi_{r_1},\chi_{r_2}; s_1,s_2) 
= \mp i  \sin\left(\frac{(s_1+s_2+\kappa_1+\kappa_2-1)\pi}{2}\right)\,,
\label{eq:TSparamIm}
\end{equation}
above and below the Stokes direction $\theta=0$. It then seems natural to {\em conjecture} that
\begin{equation}
\sigma_\pm (\chi_{r_1},\chi_{r_2}; s_1,s_2)= \exp\left(\mp \,i  \frac{(s_1+s_2+\kappa_1+\kappa_2-1)\pi}{2}\right).
\label{eq:TSexp}
\end{equation}
We remark that \eqref{eq:TSexp} generalizes a conjecture in \cite[Eq.\ 3.21]{Dorigoni:2020oon}. That work discussed
both the single parameter $s$ and the two parameters $(s_1,s_2)$ cases but there were no Dirichlet characters. In this present work, with doubly twisted
divisor sums and a pair of parameters, $(s_1,s_2)$, the parallel idea of inferring a phase factor
from a trigonometric imaginary part leads to spectacular numerical confirmation, as we now recount.

Our final result, for a pair of primitive characters, takes the form
\begin{equation}
\Xi_{s_1,s_2}(\chi_{r_1},\chi_{r_2} ; y) =    
\mathcal{S}_\pm\left[ \Xi^{{\rm Pert}}_{s_1,s_2} \right](\chi_{r_1},\chi_{r_2};y ) 
+\sigma_\pm (\chi_{r_1},\chi_{r_2}; s_1,s_2)\,\mathcal{S}_0\left[\Xi^{{\rm{NP}}}_{s_1,s_2}\right](\chi_{r_1},\chi_{r_2};y)\,,\label{eq:XiTS}
\end{equation}
with resummed perturbative contributions
\begin{align}
&\notag  \mathcal{S}_\pm\left[ \Xi^{{\rm Pert}}_{s_1,s_2} \right](\chi_{r_1},\chi_{r_2};y ) = \Xi^{{\rm LP}}_{s_1,s_2}(\chi_{r_1},\chi_{r_2};y)\\
&+  \frac{y^{s_1+s_2-1} r_1^{s_2-\frac{1}{2}}r_2^{s_1-\frac{1}{2}} \epsilon(\chi_{r_1})\epsilon(\chi_{r_2})}{\pi}  \sum_{n_1,n_2=1}^\infty \frac{\bar{\chi}_{r_1}(n_1)}{n_1^{s_2}} \frac{\bar{\chi}_{r_2}(n_2)}{n_2^{s_1}}  \mathcal{S}_\pm\left[F_{s_1,s_2,\kappa_1,\kappa_2}\right]\left( \frac{ 2\pi n_1 n_2 }{r_1 r_2y} \right)\,,
\end{align}
where $F_{s_1,s_2,\kappa_1,\kappa_2}(z)$ is  given in \eqref{eq:Faux}. The two isolated perturbative terms
\begin{equation}
\Xi^{{\rm LP}}_{s_1,s_2}(\chi_{r_1},\chi_{r_2};y)=
y^{s_2-1} A_{s_1,s_2}(\chi_{r_1},\chi_{r_2}) +y^{s_1-1} A_{s_2,s_1}(\chi_{r_2},\chi_{r_1}),
\end{equation}
are specified by \eqref{eq:XiAsyP}. The resummed non-perturbative terms are determined by
\begin{align}
& \notag \mathcal{S}_0\Big[\Xi^{{\rm{NP}}}_{s_1,s_2}\Big](\chi_{r_1},\chi_{r_2};y) =\notag  y^{s_1+s_2-1} r_1^{s_2-\frac{1}{2}} r_2^{s_1-\frac{1}{2}} \epsilon(\chi_{r_1})\epsilon(\chi_{r_2})  \sum_{n_1,n_2=1}^\infty \frac{\bar{\chi}_{r_1}(n_1)}{n_1^{s_2}} \frac{\bar{\chi}_{r_2}(n_2)}{n_2^{s_1}} \left(\frac{ 2\pi n_1 n_2 }{r_1 r_2 y}\right)^{s_1+s_2}\\
&\qquad  \times  \exp\left(- \frac{2\pi n_1 n_2}{r_1 r_2 y}\right) \int_0^\infty \exp\left(-t\,\frac{2\pi n_1n_2}{r_1 r_2 y}\right) 
\frac{t^{s_1+s_2-1}}{\Gamma(s_1+s_2)}\,
{}_2F_1(s_1,s_2,s_1+s_2\vert-t)\,{\rm d}t\,.\label{eq:XiNP1}
\end{align}

We have not proved that the transseries parameter $\sigma$ exponentiates in accordance with \eqref{eq:TSexp}.
Yet we have numerically checked its validity by evaluating both the $q$-series \eqref{eq:LamGen} and its transseries \eqref{eq:XiTS}, to high numerical precision, for many different values of $y$, of the parameters $s_1,s_2$, 
and of choices for the primitive Dirichlet characters $\chi_{r_1},\chi_{r_2}$, finding full agreement, at 100 decimal digits of numerical precision. Later, we study in detail special cases of \eqref{eq:TSexp}, which provide further strong evidence for the correctness of our conjecture, by connecting with the theory of iterated integrals of twisted Eisenstein series.

Usually when we try to construct transseries solutions to say a non-linear ODE, while the imaginary part of the transseries parameter is fixed in a fashion similar to the one presented above, the real part of the transseries parameter is fixed via some initial or boundary condition, see e.g. \cite{CostinP1} for the Painlev\'e I case. 
For the present case, we do not know which ODE or functional equation satisfied by $\Xi_{s_1,s_2}$ leads to the conjectural transseries parameter \eqref{eq:TSexp}. In a certain sense we are trying to bootstrap the full transseries entirely out of its purely perturbative expansion \eqref{eq:AsyXiSt} without relying on any underlying ODE or functional equation.
As a consequence, we are forced to start by considering $\Xi_{s_1,s_2}$ at a generic point in parameter space for $(s_1,s_2)$, and only fix $(s_1,s_2)$ to particular values (such as integers) at the very end, as we shall shortly do.

Before moving on to study impritive characters,  we develop a transseries for the non-perturbative
contribution \eqref{eq:XiNP1}.
This is achieved by simply expanding as a Gauss series the hypergometric function ${}_2F_1(s_1,s_2,s_1+s_2|-t)$ in \eqref{eq:XiNP1} and then commute the series with the integral hence performing the integral over $t$, for each power of $(-t)^n$.
 This yields the remarkable formal expression
\begin{align}
\Xi^{{\rm{NP}}}_{s_1,s_2}(\chi_{r_1},\chi_{r_2};y) =&\notag    y^{s_1+s_2-1} r_1^{s_2-\frac{1}{2}}r_2^{s_1-\frac{1}{2}} \epsilon(\chi_{r_1}) \epsilon(\chi_{r_2})\\
&\times  \sum_{n=0}^\infty(-1)^n  \frac{(s_1)_n (s_2)_n}{n!} \left(\frac{r_1 r_2 y}{2\pi}\right)^n \Xi_{s_1+n,s_2+n}\left(\bar{\chi}_{r_2},\bar{\chi}_{r_1}; \frac{1}{r_1 r_2y}\right)\,,\label{eq:XiNP2}
\end{align}
where $(s)_n \coloneqq \Gamma(s+n)/\Gamma(s)$ is a Pochhammer symbol. 
We stress that, unlike its Borel resummation \eqref{eq:XiNP1}, this transseries representation for the non-perturbative effects is a formal, 
factorially divergent series expansion in which we see $q$-series, of the form with which we began, 
but now with shifted parameters, 
$(s_1,s_2)\to (s_1+n,s_2+n)$, with characters exchanged and complex conjugated,
$(\chi_{r_1},\chi_{r_2})\to(\bar{\chi}_{r_2},\bar{\chi}_{r_1})$, and most importantly 
subject to the Fricke inversion  $y\to 1/(r_1 r_2y)$, which suppresses each term by 
the non-perturbative factor $e^{-\frac{2\pi}{r_1r_2y}}$, as $y\to0^+$.  

We emphasize that formula \eqref{eq:XiNP1} is {\em exact}, delivering the required result for 
$\Xi^{{\rm{NP}}}_{s_1,s_2}(\chi_{r_1},\chi_{r_2};y)$, at any desired numerical precision, for all $y$ with
${\rm Re}(y)>0$. In contrast, to obtain
an approximate value at small $y$ the formal transseries \eqref{eq:XiNP2} needs to be suitably truncated. While less effective, numerically, it is conceptually instructive.

Rewriting \eqref{eq:XiTS} in terms of the upper half-plane variable $\tau = i y$ and selecting the directional Borel resummation $\mathcal{S}_-$, we obtain from \eqref{eq:TSexp} and \eqref{eq:XiNP1}-\eqref{eq:XiNP2} the transseries
\begin{align}
 \Xi_{s_1,s_2}(\chi_{r_1},\chi_{r_2};\tau) 
&\notag = \mathcal{S}_-\left[\Xi_{s_1,s_2}^{\rm Pert}\right]
(\chi_{r_1},\chi_{r_2};\tau) +\mathcal{S}_0\Big[\Xi^{{\rm{NP}}}_{s_1,s_2}\Big](\chi_{r_1},\chi_{r_2};y) \\
& \notag= \mathcal{S}_-\left[\Xi_{s_1,s_2}^{\rm Pert}\right]
(\chi_{r_1},\chi_{r_2};\tau) 
+r_1^{\frac{1}{2}-s_1}r_2^{\frac{1}{2}-s_2} \epsilon(\chi_{r_1})\epsilon(\chi_{r_2}) i^{\kappa_1+\kappa_2}  \\
&\phantom{=} \times \sum_{n=0}^\infty \frac{ (s_1)_n (s_2)_n}{n!}\, \frac{(r_1 r_2 \tau)^{s_1+s_2+n-1}}{(-2\pi i)^n}
\,\Xi_{s_1+n,s_2+n}\left(\bar{\chi}_{r_2},\bar{\chi}_{r_1}; \frac{-1}{r_1 r_2\tau}\right)\,. \label{eq:TSgen}
\end{align}
For general $s_1,s_2\in \mathbb{C}$, this may be interpreted as an example of {\em quantum modularity}, 
in the sense of~\cite{ZagierQ}, for the vector-valued form
$$(\Xi_{s_1+n,s_2+n}(\chi_{r_1},\chi_{r_2};\tau), \Xi_{s_1+n,s_2+n}(\bar{\chi}_{r_2},\bar{\chi}_{r_1};\tau))_{n\in \mathbb{N}}$$
with respect to the Fricke involution $\tau \to - 1/(r_1 r_2\tau)$.

Importantly, the {\em modularity gap}, i.e.\ the difference between the original form, on the left of \eqref{eq:TSgen}, 
and those transformed on the right, by Fricke involution, is given by a resummation of the perturbative terms, in 
$\mathcal{S}_-\left[\Xi_{s_1,s_2}^{\rm Pert}\right](\chi_{r_1},\chi_{r_2};\tau) $, 
which is {\em analytic} in the upper-half $\tau$-plane.
If we specialize to the case where both characters are trivial and $(s_1,s_2)=(s,0)$,
we recover from \eqref{eq:TSgen} a much simpler example of quantum modularity
found in \cite[Eq.\ 2.25]{Dorigoni:2020oon}.
We stress that \eqref{eq:TSgen} provides the complete transseries representation for the generalized Lambert series \eqref{eq:LamGen}. Even if \eqref{eq:TSgen} has been derived by analyzing the limit in which $\tau \to i \infty$, its validity extends to all $\tau \in \mathbb{C}$ such that ${\rm Im}(\tau)>0$.
 
In Section \ref{sec:ModularPrim} we show that the transseries representation \eqref{eq:TSgen} 
simplifies dramatically for the special cases with $s_1,s_2\in \mathbb{Z}$ and 
$s_1+s_2 \equiv \kappa_1 +\kappa_2 +1 \text{ mod }2$. These are precisely the cases for which the perturbative expansion \eqref{eq:XiAsyP} contains finitely many terms and the $q$-series \eqref{eq:LamGenQ} are iterated integrals of doubly twisted Eisenstein series.
The transseries \eqref{eq:TSgen} then becomes an iterated integral version of the action of 
the Fricke involution on those doubly twisted Eisenstein series.

We conclude this subsection by specializing the transseries representation \eqref{eq:XiTS} and the non-perturbative terms \eqref{eq:XiNP2} to the case where the first character $\chi_{r_1} = \chi_r$ is a non-trivial primitive character, with $\chi_r(-1)= (-1)^\kappa$, and $\chi_{r_2} = \chi_{1,1} $ is trivial.
For these choices of character, and  $(s_1,s_2)=(s,0)$ or $(s_1,s_2)=(0,s)$,
the transseries \eqref{eq:XiTS} gives
\begin{align}
\mathcal{L}_{s}(\chi_r;y) &=  \mathcal{S}_\pm\left[ \mathcal{L}^{{\rm Pert}}_{s} \right](\chi_r;y ) +(\mp i y)^{s-1} (\mp i)^\kappa r^{-\frac{1}{2}} \epsilon(\chi_r) \tilde{\mathcal{L}}_{s}\left(\bar{\chi}_r; \frac{1}{ry}\right)\,,\label{eq:TSL}\\
\tilde{\mathcal{L}}_{s}(\chi_r;y) &=  \mathcal{S}_\pm\left[ \tilde{\mathcal{L}}^{{\rm Pert}}_{s} \right](\chi_r;y ) +(\mp i y)^{s-1} (\mp i)^\kappa r^{s-\frac{1}{2}} \epsilon(\chi_r) {\mathcal{L}}_{s}\left(\bar{\chi}_r; \frac{1}{ry}\right)\label{eq:TSLt}\,.
\end{align}
These transseries yield {\em exact} evaluations. Expressing them in terms of the upper half-plane 
modular variable $\tau = i y$ and selecting the lateral Borel resummation $\mathcal{S}_-$, we obtain
\begin{equation}
\left( \begin{matrix} 
\mathcal{L}_{s}(\chi_r;\tau) \\ 
{\tilde{\mathcal{L}}}_{s}(\chi_r;\tau ) 
\end{matrix} \right) = \left(\begin{matrix} \mathcal{S}_-\left[ \mathcal{L}^{{\rm Pert}}_{s} \right](\chi_r;\tau ) \vspace{0.1cm}\\ \mathcal{S}_-\left[ \tilde{\mathcal{L}}^{{\rm Pert}}_{s} \right](\chi_r;\tau)
\end{matrix}\right) + \frac{\epsilon(\chi_r)}{\sqrt{r}} i^\kappa  \,\tau^{s-1} \left(\begin{matrix} 0 & 1 \\ r^s & 0\end{matrix}\right)\left( \begin{matrix} 
\mathcal{L}_{s}(\bar{\chi}_r;\frac{-1}{r\tau}) \\ 
{\tilde{\mathcal{L}}}_{s}(\bar{\chi}_r;\frac{-1}{r\tau} )\label{eq:TSQM}
\end{matrix} \right)
\end{equation}
as an example of quantum modularity for the vector-valued form 
$$
(\mathcal{L}_{s}(\chi_r;\tau),\tilde{\mathcal{L}}_{s}(\chi_r;\tau),\mathcal{L}_{s}(\bar{\chi}_r;\tau),\tilde{\mathcal{L}}_{s}(\bar{\chi}_r;\tau))
$$
with modularity gaps that are analytic in the upper half-plane, ${\rm Im}(\tau)>0$.
The single non-perturbative term for $\mathcal{L}_{s}(\chi_r;\tau)$ is captured by 
$\tilde{\mathcal{L}}_{s}(\bar{\chi}_r; \frac{-1}{r\tau})$, and vice versa.

After this compelling picture of quantum modularity, with primitive characters,
we now attend to minor modifications required in cases with imprimitive characters,
such as $\chi_{6,5}$, in our motivating example from a Feynman integral in quantum field theory \cite{BV}.

 \subsection{Transseries with imprimitive characters}
\label{sec:NPImp}

As reviewed in Appendix \ref{sec:DirApp}, an imprimitive character $\chi_r$ modulo $r$ is
induced by a primitive character $\chi_D$, of lesser modulus $D|r$ with $D<r$, by the relation
\begin{equation}
\chi_r(n) = \begin{cases} \chi_D(n),&\text{if }\gcd(n,r)=1;\\
0,&\text{otherwise.}\end{cases}
\end{equation}
The $L$-function for the imprimitive character $\chi_r$ is given by 
\begin{equation}
L(\chi_r,s) = L(\chi_D,s)\sum_{d|r}\mu(d)\frac{\chi_D(d)}{d^{s}}\label{eq:Limpr}
\end{equation}
where $\mu$ is the M\"obius function. Hence only those square-free divisors $d|r$ that are 
coprime to $D$ contribute to the sum of terms. For example, 
$L(\chi_{6,5},s) = L(\chi_{3,2},s)(1+2^{-s})$, with only the divisors $d=1$ and  $d=2$ contributing 
in \eqref{eq:Limpr}. 

It is then straightforward to generalize the transseries \eqref{eq:TSgen}, to cases with $\chi_{r_i}$ induced by
$\chi_{D_i}$ and $D_i\leq r_i$, for $i\in\{1,2\}$.
The starting point is again the Mellin contour representation \eqref{eq:Mellin2} which we now rewrite by re-expressing the $L$-values $L(\chi_{r_1},t+s_1)$ and $L(\chi_{r_2},t+t_2)$ as in \eqref{eq:Limpr} thus yielding the identity
\begin{align}
 \Xi_{s_1,s_2}(\chi_{r_1},\chi_{r_2};y) = \sum_{\substack{ d_1|r_1 \\ d_2 | r_2} } \frac{\mu(d_1)\chi_{D_1}(d_1)
}{d_1^{s_1} }\cdot\frac{\mu(d_2)\chi_{D_2}(d_2)}{d_2^{s_2}}\, \Xi_{s_1,s_2}(\chi_{D_1},\chi_{D_2};d_1 d_2\,y)  \,.\label{eq:XiAsyImp}
\end{align}
Importantly, we realise that the function $\Xi_{s_1,s_2}(\chi_{r_1},\chi_{r_2};y) $ when the characters $\chi_{r_i}$ are imprimitive can be written as a finite linear combination of functions $ \Xi_{s_1,s_2}(\chi_{D_1},\chi_{D_2};  d_1 d_2\,y)$ with primitive characters $\chi_{D_i}$ and with the rescaled variable $y\to d_1 d_2 y$.

Since the inducing characters $\chi_{D_i}$ are primitive, we find that each of the summands  in \eqref{eq:XiAsyImp} can be studied separately for $y\to 0^+$ using the analysis carried out in Section \ref{sec:NPPrim}. In particular, their corresponding  transseries expansions are immediately derived from \eqref{eq:TSgen} via a simple substitution for the characters, and importantly upon rescaling $\tau\to d_1 d_2 \tau$.
Once we substitute the transseries for each summand in \eqref{eq:XiAsyImp} and rearrange the various multiplicative factors we finally arrive at 
\begin{align}
\Xi_{s_1,s_2}(\chi_{r_1},&\chi_{r_2};\tau)  \notag =
\mathcal{S}_-\left[\Xi_{s_1,s_2}^{\rm Pert}\right](\chi_{r_1},\chi_{r_2};\tau)\\
&\notag+
\sum_{\substack{ d_1|r_1 \\ d_2 | r_2} }\mu(d_1)\chi_{D_1}(d_1)(d_1D_1)^{\frac{1}{2}-s_1} 
\cdot \mu(d_2)\chi_{D_2}(d_2)(d_2D_2)^{\frac{1}{2}-s_2}\;
\frac{\epsilon(\chi_{D_1})\epsilon(\chi_{D_2})}{\sqrt{d_1 d_2}} i^{\kappa_1+\kappa_2} \\
&\label{eq:TSgenImp}\times\sum_{n=0}^\infty \frac{ (s_1)_n (s_2)_n}{n!} 
\frac{(d_1D_1d_2 D_2 \tau)^{s_1+s_2+n-1}}{(-2\pi i)^n}
 \Xi_{s_1+n,s_2+n}\left(\bar{\chi}_{D_2},\bar{\chi}_{D_1}; \frac{-1}{d_1D_1d_2D_2 \tau}\right)\,.
\end{align}

Specializing to cases where $\chi_{r_1} = \chi_r$ is an imprimitive character induced 
by $\chi_D$, while $\chi_{r_2} = \chi_{1,1}$ is trivial, we obtain, as a straightforward 
generalization of \eqref{eq:TSQM},
 \begin{align}
\left( \begin{matrix} 
\mathcal{L}_{s}(\chi_r;\tau) \\ 
{\tilde{\mathcal{L}}}_{s}(\chi_r;\tau ) 
\end{matrix} \right) {=} 
 \left(\begin{matrix} \mathcal{S}_-\left[ \mathcal{L}^{{\rm Pert}}_{s} \right](\chi_r;\tau ) \vspace{0.1cm}\\ \mathcal{S}_-\left[ \tilde{\mathcal{L}}^{{\rm Pert}}_{s} \right](\chi_r;\tau)
\end{matrix}\right) {+}\sum_{d|r}\mu(d)
\frac{\chi_D(d)\epsilon(\chi_D)}{d\sqrt{D}}\, i^\kappa  \,\tau^{s-1} 
\!\left(\begin{matrix} 0 & 1 \\ (dD)^s & 0\end{matrix}\right)\! \left( \begin{matrix} 
\mathcal{L}_{s}(\bar{\chi}_D;\frac{-1}{dD \tau}) \\ 
{\tilde{\mathcal{L}}}_{s}(\bar{\chi}_D;\frac{-1}{dD\tau})\label{eq:TSQMImp}
\end{matrix} \right) \,.
\end{align}
In Section \ref{sec:ExImp} we present examples of \eqref{eq:TSQMImp}
with terminating perturbative expansions.

\section{Modular primitives of twisted Eisenstein series}
\label{sec:ModularPrim}

In this section we analyze the transseries \eqref{eq:TSgen}, for $\Xi_{s_1,s_2}(\chi_{r_1},\chi_{r_2};\tau)$,
and its specializations with  $(s_1,s_2)=(s,0)$ in \eqref{eq:TSL} for $\mathcal{L}_s(\chi_r;\tau)$,
or $(s_1,s_2)=(0,s)$  in \eqref{eq:TSLt} for $\tilde{\mathcal{L}}_s(\chi_r;\tau)$, in the
special situations where the associated perturbative expansions contain only finitely many terms.
From \eqref{eq:CondTrunc}, we see that that this truncation occurs when $s_1$ and $s_2$ are integers satisfying
$s_1+s_2 \equiv \kappa_1 +\kappa_2 +1 \text{ mod }2$, where $\kappa_i\in\{0,1\}$ encodes
the parity of the primitive character $\chi_{r_i}$.
In these terminating cases, we show that $\Xi_{s_1,s_2}(\chi_{r_1},\chi_{r_2};\tau)$, and hence 
$\mathcal{L}_s(\chi_r;\tau)$ or $\tilde{\mathcal{L}}_s(\chi_r;\tau)$,
are related, by integrations, to holomorphic modular forms called twisted Eisenstein series.
Their transseries representations are then direct consequences of a Fricke involution 
relating  such twisted Eisenstein series.

\subsection{Twisted Eisenstein series}

The symmetry $\Xi_{s_1,s_2}(\chi_{r_1},\chi_{r_2};\tau)=\Xi_{s_2,s_1}(\chi_{r_2},\chi_{r_1};\tau)$
implies that, without loss of generality, we may assume that $s_1\geq s_2$
and hence that $m=s_1-s_2+1$ is a positive integer, whose parity is given by
$m \equiv \kappa_1 +\kappa_2 \text{ mod }2$, to ensure that the perturbative terms are finite in number
and that the transseries parameter of \eqref{eq:TSexp} reduces to a mere sign independent from the $\pm$ choice in resummation,
\begin{equation}
\sigma_\pm(\chi_{r_1},\chi_{r_2}; s_1,s_2)= \exp\left(\mp \,i  \frac{(s_1+s_2+\kappa_1+\kappa_2-1)\pi}{2}\right)
=(-1)^{s_1+\frac{\kappa_1+\kappa_2-m}{2}}\,.\label{eq:TSred}
\end{equation}

The vanishing  of ${\rm Im}(\sigma)$  is to be expected, when the perturbative terms are finite in number
and the issue of directional Borel resumation does not arise. Yet $\sigma\in\{-1,1\}$ still multiplies
the non-perturbative terms $\Xi^{\rm NP}_{s_1,s_2}(\chi_{r_1},\chi_{r_2};\tau)$, whose existence
we might have missed, had it not been for the problem of an infinite number of perturbative terms,
in a more demanding case. The mathematician C.L.\ Dodgson ({\em alias} Lewis Carrol) imagined
a comparable phenomenon, concerning the survival of the smile of a Cheshire cat.
\begin{quotation}\noindent{\em
``All right,'' said the Cat; and this time it vanished quite slowly,
beginning with the end of the tail, and ending with the grin, which
remained some time after the rest of it had gone.
``Well! I’ve often seen a cat without a grin,'' thought Alice; ``but a
grin without a cat! It’s the most curious thing I ever saw in my life!''}
\vspace{0.1cm}

(Charles Lutwidge Dodgson, Alice's Adventures in Wonderland, 1865.)
\end{quotation} 
This  survival of non-perturbative terms, when the perturbative terms are finite in number,
has been referred to as {\em Cheshire-cat resurgence} a phenomenon first observed in quantum mechanical examples~\cite{Dunne:2016jsr,Kozcaz:2016wvy} but also in quantum field theory setups, see e.g.~\cite{Dorigoni:2017smz,Dorigoni:2019kux}.
It was observed in \cite{Dorigoni:2020oon} for untwisted Lambert series
$\mathcal{L}_s(q)=\sum_{n=1}^\infty \sigma_s(n)n^{-s}q^n$ when $s$ is an odd integer, 

With a finite number of perturbative terms, we have no need of the symbol
$\mathcal{S}_-$ in the transseries \eqref{eq:TSgen}, which we rewrite as
\begin{align}
 \Xi_{s_1,s_2}(\chi_{r_1},\chi_{r_2};\tau) & \notag = \Xi_{s_1,s_2}^{\rm Pert}
(\chi_{r_1},\chi_{r_2};\tau) +r_1^{\frac{1}{2}-s_1}r_2^{\frac{1}{2}-s_2} \epsilon(\chi_{r_1})\epsilon(\chi_{r_2}) i^{\kappa_1+\kappa_2}  \\
 &\phantom{=}\times \sum_{n=0}^\infty \frac{ (s_1)_n (s_2)_n}{n!}\, \frac{(r_1 r_2 \tau)^{s_1+s_2+n-1}}{(-2\pi i)^n}
\,\Xi_{s_1+n,s_2+n}\left(\bar{\chi}_{r_2},\bar{\chi}_{r_1}; \frac{-1}{r_1 r_2\tau}\right)\,. \label{eq:TSgenRed}
\end{align}
Specializing  to the simpler cases $\mathcal{L}_s(\chi_r;\tau)$ and $\tilde{\mathcal{L}}_s(\chi_r;\tau)$
we require that $s=m-1$ where $m$ is a positive integer with the same parity as $\kappa\in\{0,1\}$, where
$\chi_r(-1)=(-1)^\kappa$. In these cases, we obtain
\begin{align}
&\mathcal{L}_{m-1}(\chi_r;\tau) = \sum_{k=-1}^{m+\kappa-2} \frac{(2\pi i \tau)^k}{k!} \zeta(-k) L(\chi_r,m-k-1) +i^\kappa \frac{\epsilon(\chi_r)}{\sqrt{r}} \tau^{m-2} \tilde{\mathcal{L}}_{m-1}\left(\bar{\chi}_r; \frac{-1}{r\tau}\right)\,,
\label{eq:TSLTr}\\
&\tilde{\mathcal{L}}_{m-1}(\chi_r;\tau) \notag = \left[ \frac{L'(\chi_r ,2-m)}{L(\chi_r , 2-m)}  
+ \gamma+\psi(m-1) - \log(-2\pi  i \tau)\right]  \frac{(2\pi i  \tau)^{m-2}}{(m-2)!} L(\chi_r , 2-m)\\
&\qquad \qquad +  \sum_{\substack{ k=0 \\ k\neq m-2} }^{m+\kappa-1} \frac{(2\pi i \tau)^k}{k!} \zeta(m-k-1) L(\chi_r,-k)+ i^\kappa  \epsilon(\chi_r) \sqrt{r}(r\tau)^{m-2} {\mathcal{L}}_{m-1}\left(\bar{\chi}_r; \frac{-1}{r\tau}\right)
\label{eq:TSLtTr}\,,
\end{align}
from the perturbative terms in \eqref{eq:Asypert} and \eqref{eq:Asypert2},
on the understanding that the term with $k=-1$ in \eqref{eq:TSLTr} is obtained as the limit $k\to-1$.

As already stressed below \eqref{eq:XiNP1}, if we had access to a suitable ODE or functional equation to which $\mathcal{L}_{s}(\chi_r;\tau) $ and ${\tilde{\mathcal{L}}}_{s}(\chi_r;\tau ) $ are solutions of, we would be able to prove the validity of the transseries expansions \eqref{eq:TSQM} since the imaginary part of the transseries parameter here conjectured would be uniquely fixed by some initial or boundary condition. Due to the lack of such ODE or functional equations, to analyze the special cases $\mathcal{L}_{m-1}(\chi_r;\tau),\tilde{\mathcal{L}}_{m-1}(\chi_r;\tau)$ with $m\in \mathbb{N}$ here presented, we have been forced to considered a deformation in parameter space $s\in \mathbb{C}$ and reach such cases as limits of  \eqref{eq:TSQM}.
Although this deformation is by no means unique, the transseries expansion \eqref{eq:TSQM} and hence \eqref{eq:TSLTr}-\eqref{eq:TSLtTr} are. These transseries expansions are an exact encoding of the $q$-series \eqref{eq:LamChi}-\eqref{eq:qser2} when expanded near the essential singularity at $y=0$ as strongly suggested by the numerical tests previously mentioned and confirmed  even further by the modularity arguments we now present.

We shall show that identities \eqref{eq:TSgenRed}, \eqref{eq:TSLTr} and \eqref{eq:TSLtTr} are a consequence of 
how Fricke involution acts on modular forms called twisted Eisenstein series, which we now define.
 
Given a positive integer $m$ and a pair of  primitive characters $(\chi_{r_1},\chi_{r_2})$, of which at least one is non-trivial,
we define a twisted Eisenstein  $q$-series by 
\begin{equation}
G_m(\chi_{r_1},\chi_{r_2};\tau) = A_m(\chi_{r_1},\chi_{r_2})+ \sum_{n=1}^\infty  \sigma_{m-1}(\chi_{r_1},\chi_{r_2};n)\,q^n\,,\label{eq:qserG2chi}
\end{equation}
where $q=e^{2\pi i \tau}$ and the doubly twisted divisor function, $ \sigma_{s}(\chi_{r_1},\chi_{r_2};n)$, is
defined in  \eqref{eq:div2chi}.
The constant term in \eqref{eq:qserG2chi} depends on the nature of the characters considered:
\begin{equation}\label{eq:G0mode}
A_m(\chi_{r_1},\chi_{r_2}) = \begin{cases}
\frac12L(\chi_{r_2},1-m)\,,&\text{if }\chi_{r_1} = \chi_{1,1}\,;\\
\tfrac12\delta_{m,1}L(\chi_{r_1},0)\,,&\text{if }\chi_{r_2} = \chi_{1,1}\,;\\
0\,,&{\rm otherwise.}\end{cases}
\end{equation}

Let $\mathcal{M}_m(\Gamma_0(r),\chi_r)$ be the vector space of holomorphic modular forms,
of weight $m$ and character $\chi_r$, for the congruence subgroup
$\Gamma_0(r) \coloneqq\{ \left(\begin{smallmatrix} a & b \\ c & d \end{smallmatrix}\right) \in {\rm SL}(2,\mathbb{Z})\,,\, c \equiv 0 \,{\rm mod} \,r\}$, with action
\begin{align}
f(\tau)\in\mathcal{M}_m(\Gamma_0(r),\chi_r) \, \Rightarrow \,  f\left( \frac{a \tau+b}{c\tau+d}\right) = \chi_r(d) (c\tau+d)^m f(\tau) \,,\quad \forall\,\gamma=\left(\begin{matrix} a & b \\ c & d \end{matrix} \right)\in \Gamma_0(r)\,.
\end{align}
If $m \equiv \kappa_1 +\kappa_2\text{ mod }2$, then $G_m(\chi_{r_1},\chi_{r_2};\tau)$ 
defined by $q$-series  \eqref{eq:qserG2chi} is an element of the vector space $\mathcal{M}_m(\Gamma_0(r_1 r_2),\chi_{r_1} \chi_{r_2})$. 
Moreover, and most importantly,  we have
\begin{equation}
G_m(\chi_{r_1},\chi_{r_2}; \tau) = i^{\kappa_1 +\kappa_2} \epsilon(\chi_{r_1})\epsilon(\chi_{r_2}) r_1^{\frac{1}{2}-m} r_2^{-\frac{1}{2}}  \tau^{-m} G_m\left(\bar{\chi}_{r_2},\bar{\chi}_{r_1}; \frac{-1}{r_1 r_2\tau}\right)\,\label{eq:FrickeG2chi}
\end{equation}
under the Fricke involution $\tau\to -1/(r_1 r_2\tau)$, which is not an element of $\Gamma_0(r_1r_2)$.

When one of the characters is trivial, we obtain the singly twisted  Eisenstein series
\begin{align}
G_m(\chi_r; \tau)  &= G_m(\chi_r,\chi_{1,1};\tau)= 
\delta_{m,1} \frac{L(\chi_{r},0)}{2} + \sum_{n=1}^\infty \sigma'_{m-1,\chi_r}(n) q^n\,, \label{eq:TwEis2}\\
E_m(\chi_r ; \tau) &= G_m(\chi_{1,1},\chi_r;\tau) =  \frac{L(\chi_r,1-m)}{2} 
+ \sum_{n=1}^\infty \sigma_{m-1,\chi_r}(n) q^n\,,\label{eq:TwEis1}
\end{align}
for positive  integers $m$  with the same parity as the primitive character $\chi_r$. In these cases, the 
divisor functions are given in \eqref{eq:Twsigma2} and \eqref{eq:Twsigma1},
as special cases of \eqref{eq:div2chi}. Here 
both $E_m(\chi_r ; q)$ and $G_m(\chi_r ; q)$ belong to the vector space 
$\mathcal{M}_k(\Gamma_0(r),\chi_r)$ of holomorphic modular forms with weight $m$ and character $\chi_r$
and are related by the Fricke involution $\tau \to -1/(r\tau)$,  which gives
\begin{equation}
G_m(\chi_r ;\tau) = i^\kappa \epsilon(\chi_r) r^{\frac{1}{2}-m} \tau^{-m} E_m\left( \bar{\chi}_r;\frac{-1}{r\tau}\right)\,.\label{eq:Fricke}
\end{equation}

In  the general case of $q$-series with $s_1=s_2+ m-1>0$ and $m$ a positive integer of the same parity as 
$\kappa_1+\kappa_2$, differentiation of \eqref{eq:LamGen} gives the crucial differential equation
\begin{equation}
\left(q \frac{   d }{dq}\right)^{s_1} \Xi_{s_1,s_2}(\chi_{r_1},\chi_{r_2} ; q)  = G^0_m(\chi_{r_1},\chi_{r_2};q)
\coloneqq G_m(\chi_{r_1},\chi_{r_2};q) - A_m(\chi_{r_1},\chi_{r_2})\,,\label{eq:IterIntGen}
\end{equation}
where $G^0_m(\chi_{r_1},\chi_{r_2};q)$ vanishes as $q\to0$. 
The corresponding single-parameter cases are
\begin{align}
\left(q \frac{   d }{dq}\right)^{m-1}\mathcal{L}_{m-1}(\chi_r ; q) & \label{eq:IterInt}  =G^0_{m}(\chi_r;q)
\coloneqq  G_{m}(\chi_r;q) -\delta_{m,1} \frac{L(\chi_{r},0)}{2}\,,\\
\left(q \frac{   d }{dq}\right)^{m-1}\tilde{\mathcal{L}}_{m-1}(\chi_r ; q) &=E^0_{m}(\chi_r;q)
: =E_{m}(\chi_r;q) -  \frac{L(\chi_r,1-m)}{2} \,.\label{eq:IterIntt}
\end{align}
These differential equations enable us to interpret the series $\Xi_{s_1,s_2}(\chi_{r_1},\chi_{r_2} ; q)$,
and its specializations, as iterated integrals of twisted  Eisenstein series, from which they inherit
the structure under Fricke inversion apparent in \eqref{eq:TSgenRed}, \eqref{eq:TSLTr}
and \eqref{eq:TSLtTr}.

\subsection{Fricke involution and transseries representation}

Iterated integration of \eqref{eq:IterIntGen}, for positive integers $s_1=s_2+m-1$, gives
\begin{equation}
\Xi_{s_1,s_2}(\chi_{r_1},\chi_{r_2};q) = \int_{0\leq q_1 \leq ...\leq q_{s_1}\leq q} {\rm d}\log q_1 \cdots {\rm d}\log q_{s_1}\, G^0_{m}(\chi_{r_1},\chi_{r_2};q_1)\,.\label{eq:XiItInt}
\end{equation}
Now let us consider the collection of all $\Xi_{s_1,s_2}(\chi_{r_1},\chi_{r_2};q)$ as we vary $s_1\in \mathbb{Z}$,
while keeping fixed the positive integer $m = s_1-s_2 +1$, with $m\equiv \kappa_1+\kappa_2\text{ mod }2$.
We discuss, in turn,  three distinct regimes: first $s_1\leq 0$, then $0< s_1 \leq m-1$, and finally $s_1 >m-1$.

With $s_1\leq0$, we interpret \eqref{eq:IterIntGen}
as saying that
$\Xi_{s_1,s_2}(\chi_{r_1},\chi_{r_2};q)$ is given by a $|s_1|$-fold derivative of $G_m^0(\chi_{r_1},\chi_{r_2};q)$. 
In this case the hypergeometric function appearing in the non-perturbative corrections \eqref{eq:XiNP1} 
reduces to a polynomial of degree $|s_1|$.
Thus the transseries \eqref{eq:TSgen} simplifies, dramatically, since the non-perturbative terms terminate
at $n=|s_1|$. Then we obtain a result equivalent to the $|s_1|$-fold derivative, w.r.t.\ $\tau$,
of the Fricke relation \eqref{eq:FrickeG2chi}. In particular, for $s_1=0$ the transseries \eqref{eq:TSgen} becomes directly
that Fricke relation. This fact becomes obvious when we compare directly the $q$-series \eqref{eq:LamGenQ} and \eqref{eq:qserG2chi}, thus yielding the identities
\begin{align}
\Xi_{0,1-m}(\chi_{r_1},\chi_{r_2};q) & = G^0_m(\chi_{r_1},\chi_{r_2};q) \,,\qquad  \Xi_{1-m,0}(\chi_{r_1},\chi_{r_2};q)  = G^0_m(\chi_{r_2},\chi_{r_1};q)\,, \\
\tilde{\mathcal{L}}_{1-m}(\chi_r;q) & = G^0_m(\chi_r;q) \,, \qquad \qquad\qquad\quad\,\,  \mathcal{L}_{1-m}(\chi_r;q) =  E^0_m(\chi_r;q)\,.
\end{align}

Moving to the second interval, $0< s_1 \leq m-1$, and hence $s_2\leq0$, we find that the hypergeometric 
function in the non-perturbative corrections \eqref{eq:XiNP1} gives a polynomial of degree $|s_2|$ and 
the transseries \eqref{eq:TSgen} simplifies to an upper-triangular structure, 
with the non-perturbative terms terminating at $n=|s_2|$.
For example, with $s_1-s_2+1=m=3$ and  $\kappa_1+\kappa_2=1$,  we obtain
\begin{align}
 \left( \begin{matrix} \Xi_{1,\text{-}1}(\chi_{r_1},\chi_{r_2};\tau) \\ \Xi_{2,0}(\chi_{r_1},\chi_{r_2};\tau)
\end{matrix} \right) &\notag =
\left(\begin{matrix} 
\Xi^{\rm Pert}_{1,\text{-}1}(\chi_{r_1},\chi_{r_2};\tau) \\ 
\Xi^{\rm Pert}_{2,0}(\chi_{r_1},\chi_{r_2};\tau)\end{matrix} \right) \\
&\phantom{=}+  i  \epsilon(\chi_{r_1})\epsilon(\chi_{r_2}) r_1^{-\frac{3}{2}} r_2^{\frac{1}{2}} 
\left(\begin{matrix}  \tau^{-1} &  \frac{r_1r_2}{2\pi i } \\ 0 &  r_1 r_2 \tau\end{matrix} \right)
\left(\begin{matrix} 
\Xi_{1,\text{-}1}(\bar{\chi}_{r_2},\bar{\chi}_{r_1};\frac{-1}{r_1 r_2\tau}) \\ 
\Xi_{2,0}(\bar{\chi}_{r_2},\bar{\chi}_{r_1};\frac{-1}{r_1 r_2\tau})\end{matrix} \right)\,,\label{eq:UpperTrEx}
\end{align}
where we highlight the fact that the final Fricke-dual vector is multiplied by an upper-triangular matrix.

This structure can be understood by combining the iterated integral \eqref{eq:XiItInt} with the Fricke inversion \eqref{eq:FrickeG2chi}. Starting from \eqref{eq:XiItInt} we  write the \textit{Eichler integral} representation
\begin{align}
\Xi_{s_1,s_2}(\chi_{r_1},\chi_{r_2};\tau) \notag &\notag =-\frac{(2\pi i)^{s_1}}{ \Gamma(s_1)}\int_\tau^{i \infty}
 (\tau-\tau_1)^{s_1-1} G^0_m(\chi_{r_1},\chi_{r_2};\tau_1) {\rm d}\tau_1\\
&= -\frac{(2\pi i)^{s_1}}{ \Gamma(s_1)}\int_\tau^{i \infty} (\tau-\tau_1)^{s_1-1} G_m(\chi_{r_1},\chi_{r_2};\tau_1) {\rm d}\tau_1 - \frac{(2\pi i \tau)^{s_1}}{\Gamma(s_1+1)} A_{m}(\chi_{r_1},\chi_{r_2})\,,\label{eq:IterInt2}
\end{align}
where $A_{s_1,s_2}(\chi_{r_1},\chi_{r_2})$, from \eqref{eq:G0mode}, is the constant term of the modular twisted 
Eisenstein series.
This is clearly correct for $s_1=1$. For $s_1>1$, the polynomial factor  $(\tau-\tau_1)^{s_1-1}$ 
in the integral over $\tau_1$ avoids the multiple integrations in~\eqref{eq:XiItInt}, as may be shown by 
binomial expansion. However, when we add to $G^0_m$ a constant, required for modularity
of $G_m$, we may encounter problems as $\tau_1\to i\infty$. In such cases, the integral should be 
understood in the sense of tangential-basepoint regularization, described in~\cite{Brown:2014}. 
For our purposes, this amounts to the 
prescription~$\int_\tau^{i\infty} \tau_1^k {\rm d}\tau_1 {=} {-} \tau^{k+1}/(k+1)$, for~$k\geq0$.

Applying  the Fricke inversion \eqref{eq:FrickeG2chi} to \eqref{eq:IterInt2}, we encounter the Fricke dual Eisenstein series 
$G_m(\bar{\chi}_{r_2},\bar{\chi}_{r_1}; \tilde{\tau}_1) $, with $\tilde{\tau}_1 = -1/(r_1 r_2 \tau_1)$. 
With $\tilde{\tau}_1$ as the new integration variable, this gives
\begin{equation}
\Xi_{s_1,s_2}(\chi_{r_1},\chi_{r_2};\tau) = - \frac{(2\pi i \tau)^{s_1}}{\Gamma(s_1+1)} A_{m}(\chi_{r_1},\chi_{r_2})+
 i^{\kappa_1 +\kappa_2} \epsilon(\chi_{r_1})\epsilon(\chi_{r_2}) r_1^{\frac{1}{2}-s_1} r_2^{\frac{1}{2}-s_2}I
\label{eq:IntFricke}
\end{equation}
where we isolate the term in \eqref{eq:IterInt2} that is proportional  to $\tau^{s_1}$, leaving us to compute the  integral
\begin{equation}
I= (-1)^{s_1+s_2} \frac{(2\pi i)^{s_1}}{ \Gamma(s_1)}
 \tilde{\tau}^{1-s_1} \int_{\tilde{\tau}}^0  \tilde{\tau}_1^{m-s_1-1} \left( \tilde{\tau}- \tilde{\tau}_1\right)^{s_1-1} G_m(\bar{\chi}_{r_2},\bar{\chi}_{r_1}; \tilde{\tau}_1){\rm d}\tilde{\tau}_1\,,
\end{equation} 
over the dual modular function $G_m(\bar{\chi}_{r_2},\bar{\chi}_{r_1};\tilde{\tau}_1)$, multiplied by a polynomial 
in $\tilde{\tau}_1$. The lower limit of integration is $\tilde{\tau}\coloneqq-1/(r_1r_2\tau)$.

We now decompose the integral into three terms:
\begin{align}
& I =\notag  (-1)^{s_1+s_2} \frac{(2\pi i)^{s_1}}{ \Gamma(s_1)}   \tilde{\tau}^{1-s_1}\left[
\int_{\tilde{\tau}}^0\tilde{\tau}_1^{m-s_1-1} (\tilde{\tau}-\tilde{\tau}_1)^{s_1-1} 
A_m(\bar{\chi}_{r_2},\bar{\chi}_{r_1})\,{\rm d}\tilde{\tau}_1\right. \\
&\left.\!-\!\int_0^{i\infty}\! \tilde{\tau}_1^{m-s_1-1} (\tilde{\tau}-\tilde{\tau}_1)^{s_1-1}
G^0_m(\bar{\chi}_{r_2},\bar{\chi}_{r_1}; \tilde{\tau}_1)\,{\rm d}\tilde{\tau}_1+\int_{\tilde{\tau}}^{i\infty}\! \tilde{\tau}_1^{m-s_1-1} (\tilde{\tau}-\tilde{\tau}_1)^{s_1-1}
G^0_m(\bar{\chi}_{r_2},\bar{\chi}_{r_1}; \tilde{\tau}_1)\,{\rm d}\tilde{\tau}_1
\right].\label{eq:IntStep}
\end{align}
The first integral comes from the constant term in $G_m(\bar{\chi}_{r_2},\bar{\chi}_{r_1}; \tilde{\tau}_1)$
and is easily performed, since the integrand is a polynomial. The second integral,
from $0$ to $i\infty$, is closely related to the period polynomial
\begin{align}
P_m(\chi_{r_1},\chi_{r_2};\tau)&\notag \coloneqq \int_0^{i\infty} (\tau-\tau_1)^{m-2} G_m^0(\chi_{r_1},\chi_{r_2};\tau_1){\rm d}\tau_1\\
&= - \sum_{\ell=0}^{m-2} \frac{(m-2)!}{(m-\ell-2)!} \frac{L(\chi_{r_1},\ell+1) L_(\chi_{r_2},\ell+2-m)}{(2\pi i )^{\ell+1}} \tau^{m-\ell-2}\,,\label{eq:PeriodPoly}
\end{align}
whose evaluation follows from the identity
\begin{equation}
\int_0^{i\infty} \tau_1^\ell  G_m^0(\chi_{r_1},\chi_{r_2};\tau_1){\rm d}\tau_1 
= \frac{\ell\,!}{(-2\pi i )^{\ell+1}}  L(\chi_{r_1},\ell+1) L(\chi_{r_2},\ell+2-m)\,.\label{eq:IntEll}
\end{equation}
For $s_1=m-1$, the integral from $0$ to $i\infty$ in \eqref{eq:IntStep} is precisely the period polnomial
$P_m(\bar{\chi}_{r_2},\bar{\chi}_{r_1}; \tilde{\tau})$.  More generally, \eqref{eq:IntEll} 
delivers it as a Laurent polynomial in $\tau=-1/(r_1r_2\tilde{\tau})$. 

We now  attend to the third, and most interesting contribution to \eqref{eq:IntStep}, obtaining
\begin{align}
& \notag (-1)^{s_1+s_2} \frac{(2\pi i)^{s_1}}{ \Gamma(s_1)}   \tilde{\tau}^{1-s_1} \int_{\tilde{\tau}}^{i\infty} 
\tilde{\tau}_1^{m-s_1-1} \left( \tilde{\tau}- \tilde{\tau}_1\right)^{s_1-1} 
G^0_m(\bar{\chi}_{r_2},\bar{\chi}_{r_1}; \tilde{\tau}_1){\rm d}\tilde{\tau}_1 \\
&\notag =\sum_{n=0}^{m-s_1-1} \frac{(s_1)_n (s_2)_n}{n!}  
\frac{(-\tilde{\tau})^{1-s_1-s_2-n}}{(-2\pi i)^n} \left[-\frac{(2\pi i)^{s_1+n}}{\Gamma(s_1+n)} \int_{\tilde{\tau}}^{i\infty} (\tilde{\tau}-\tau_1)^{s_1+n-1} G^0_m(\bar{\chi}_{r_2},\bar{\chi}_{r_1};\tau_1){\rm d}\tau_1\right]\\
&\label{eq:IntDual} =\sum_{n=0}^{m-s_1-1} \frac{(s_1)_n (s_2)_n}{n!}  
\frac{(-\tilde{\tau})^{1-s_1-s_2-n}}{(-2\pi i)^n} \Xi_{s_1+n,s_2+n}(\bar{\chi}_{r_2},\bar{\chi}_{r_1};\tilde{\tau})\,,
\end{align}
by developing $\tilde{\tau}_1^{m-s_1-1}$ as a polynomial in $(\tilde{\tau}- \tilde{\tau}_1)$, in the first step,
and then using an integral representation from \eqref{eq:IterInt2}, in the second. It is here that we find
the non-perturbative terms in the transseries for $\Xi_{s_1,s_2}(\chi_{r_1},\chi_{r_2};\tau)$
in case that $s_1-s_2=m-1$,  with $s_1>0$ and $s_2\leq0$.  
Straightforward, if tedious, algebra shows that the other power-behaved contributions that we 
have accumulated combine to give the Laurent polynomial perturbative terms $\Xi_{s_1,s_2}^{\rm Pert}(\chi_{r_1},\chi_{r_2};\tau)$ thus concluding our proof that the iterated integral formulation \eqref{eq:IterInt2} coincides with the transseries representation \eqref{eq:TSgenRed}.

The transseries  \eqref{eq:TSgenRed} exhibits an upper triangular structure, of which
\eqref{eq:UpperTrEx} is an example. If we regard a Fricke inversion of $\Xi_{s_1,s_2}(\chi_{r_1},\chi_{r_2};\tau)$
as the $s_1$-fold iterated integral \eqref{eq:XiItInt} with $0<s_1\leq m-1$, then we produce
all of the terms $\Xi_{s'_1,s'_2}(\bar{\chi}_{r_2},\bar{\chi}_{r_1};\tilde{\tau})$ with $s_1\leq s'_1 <m-1$, 
as in  example \eqref{eq:UpperTrEx}. 
If we consider the collection of all $\Xi_{s_1,s_2}(\chi_{r_1},\chi_{r_2};\tau)$ at fixed $s_1-s_2+1=m$ with 
$s_1 \in \{1,...,m-1\}$ as a $(m-1)$ dimensional vector, then under Fricke inversion this vector transforms with 
an $(m-1)\times (m-1)$ upper triangular matrix, plus a $(m-1)$-vector of cocycles, i.e.\ the 
perturbative Laurent polynomials in $\tau$. 
The story presented here is reminiscent of the analysis of \cite{Dorigoni:2020oon} 
where the authors considered special cases of iterated integrals \eqref{eq:XiItInt} for which both
characters are trivial. In these cases, the upper triangular structure is deeply connected to the 
broader theory of iterated integrals of untwisted holomorphic Eisenstein series,
developed by Francis Brown \cite{Brown:2014,Brown:2017}.

It is worth noting that examples of iterated integrals $\Xi_{s_1,s_2}(\chi_{r_1},\chi_{r_2};\tau)$ in this range of parameters do appear in string theory contexts also for non-trivial characters. If we specialise the upper triangular structure example \eqref{eq:UpperTrEx} to the case of a single odd non-trivial character $\chi_r = \chi_{3,2}$ we find that it encodes the action of Fricke inversion for the Eichler integrals of weight-$3$ twisted Eisenstein series:
\begin{align}
\Xi_{1,-1}(\chi_{1,1},\chi_{3,2};q) &\label{eq:E3Eich} = -\int_\tau^{i\infty} E_3^0(\chi_{3,2};\tau_1) {\rm d}\tau_1\,,\quad \Xi_{2,0}(\chi_{1,1},\chi_{3,2};q)  = -2\pi i \int_\tau^{i\infty} (\tau-\tau_1) E_3^0(\chi_{3,2};\tau_1) {\rm d}\tau_1 \,,\\
\Xi_{1,-1}(\chi_{3,2},\chi_{1,1};\tau) &\label{eq:G3Eich} = -\int_\tau^{i\infty} G_3^0(\chi_{3,2};\tau_1) {\rm d}\tau_1 \,,\quad  \Xi_{2,0}(\chi_{3,2},\chi_{1,1};\tau)  = -2\pi i \int_\tau^{i\infty} (\tau-\tau_1) G_3^0(\chi_{3,2};\tau_1) {\rm d}\tau_1 \,. 
\end{align}
Simple linear combinations of \eqref{eq:E3Eich} yield particular periods associated with the universal family of elliptic curves over the modular curve $X_1(3) \cong \mathbb{H} \backslash \Gamma_1(3)$ which are in turn 
deeply intertwined \cite{Bousseau:2022snm} with the study of BPS states in type IIA string theory compactified on a particular Calabi-Yau (CY) threefold called local $\mathbb{P}^2$.\footnote{We thank Veronica Fantini for bringing this connection to our attention.}
From a geometric perspective, the Eichler integrals \eqref{eq:E3Eich} become relevant when discussing such periods in the large radius limit of the underlying background CY geometry, while the Fricke related integrals \eqref{eq:G3Eich} when expanding around the conifold point. It would be extremely interesting to understand whether this connection extends further.

This closure property, under Fricke inversion, of a finite set of iterated integrals,
cannot hold for $s_1 >m-1\ge0$, and hence $s_2>0$, since we know  
that when $s_1$ and $s_2$ are both positive integers the non-perturbative part of
the transseries series continues indefinitely. Here the iterated integrals of Eisenstein integrals
lie outside the class analyzed in \cite{Brown:2014}. The origin of this obstacle is apparent
in the integrands of \eqref{eq:IntStep}, where the Eisenstein series are no longer mutiplied by
polynomials in $\tilde{\tau}_1$, if $s>m-1$. Thus we encounter the  infinite tower of 
higher and higher iterated integrals that was derived in \eqref{eq:TSgen}.
We remark  that in the untwisted case,  $\chi_{r_1}=\chi_{r_2} = \chi_{1,1}$, these ``overly'' integrated cases, i.e.\ 
$\Xi_{s_1,s_2}(\chi_{1,1},\chi_{1,1};\tau)$ with $s_1>0$ and $s_2>0$
appear, in very special linear combinations, in the study of one-loop scattering amplitudes of open-string 
states \cite{Broedel:2018izr,Dorigoni:2021jfr,Dorigoni:2021ngn}. 
We are not aware of similar constructions that give rise to the more general, overly integrated, 
doubly-twisted transseries for $\Xi_{s_1,s_2}(\chi_{r_1},\chi_{r_2};\tau)$ considered in this paper.
It would be also very interesting to understand whether single-valued versions of the iterated integrals here considered do reproduce the real-analytic modular forms constructed in \cite{DREWITT20251,Duhr:2025lvr}.

We conclude this section by remarking on the special cases
\allowdisplaybreaks{
\begin{align}
\mathcal{L}_{m-1}(\chi_r;\tau) &= \Xi_{m-1,0}(\chi_{r},\chi_{1,1};\tau)
=-\frac{(2\pi i)^{m-1}}{ \Gamma(m-1)}\int_\tau^{i \infty} (\tau-\tau_1)^{m-2} G^0_m(\chi_{r};\tau_1) {\rm d}\tau_1 \,, \\
{\tilde{\mathcal{L}}}_{m-1}(\chi_r;\tau) &=\Xi_{m-1,0}(\chi_{1,1},\chi_{r};\tau)=-\frac{(2\pi i)^{m-1}}{ \Gamma(m-1)}\int_\tau^{i \infty} (\tau-\tau_1)^{m-2} E^0_m(\chi_{r};\tau_1) {\rm d}\tau_1\,, \label{eq:IterIntLLt}
\end{align}}
when  $m>1$ has the same parity as the primitive character $\chi_r$ and hence the integrands
involve modular forms, from which constant terms were removed, 
in \eqref{eq:IterInt} and \eqref{eq:IterIntt}.
It is straightforward to specialize \eqref{eq:TSgenRed} to either case. For example we obtain
\begin{align}
 \mathcal{L}_{m-1}(\chi_r;\tau) = &\notag  - \delta_{m,1}  \frac{L(\chi_r,0)}{2}- i^{\kappa} \epsilon(\chi_{r})r^{\frac{3}{2}-m} (-1)^{m-1} 
(2\pi i)^{m-1} \times  \\
&  \left[  \frac{L(\bar{\chi}_r,1-m) }{2 (m-1)!} \tilde{\tau}+  \frac{{\tilde{\tau}}^{2-m}}{(m-2)!} P_m(\chi_{1,1},{\bar{\chi}}_r;\tilde{\tau}) \right] +i^{\kappa} \epsilon(\chi_{r}) r^{-\frac{1}{2}}   \tau^{2-m}\tilde{\mathcal{L}}_{m-1}(\bar{\chi}_r;\tilde{\tau})
\label{eq:TSIterIntL}
\end{align}
with  a period polynomial \eqref{eq:PeriodPoly}, evaluated at $\tilde{\tau} = -1/({r\tau})$,
neatly summing all but two of the perturbative terms, no matter how large the modular weight $m$ might be. Note that for $m=1$ the period polynomial vanishes identically, hence the above expression is regular for all $m\geq1$.
The final term, involving  $\tilde{\mathcal{L}}_{m-1}(\bar{\chi}_r;\tilde{\tau})$ is non-perturbative; in Lewis Carrol's terms, one might regard it as the grin of his Cheshire cat. There is no need
for directional Borel resummation, since the perturbative terms are manifestly finite in number. 
Yet a non-perturbative term survives, from the general result in \eqref{eq:TSgen}
that emerged when the formal perturbative series was factorially divergent.

In the cases with
s$_1,s_2\in \mathbb{Z}$ and $s_1-s_2+1=m\equiv \kappa_1+\kappa_2\text{ mod }2$,
considered in the section, we encountered iterated integrals of modular twisted Eisenstein series,
i.e.\ {\em modular primitives}. The terminating perturbative series
are the modularity gaps between $\Xi_{s_1,s_2}(\chi_{r_1},\chi_{r_2};\tau)$
and Fricke duals $\Xi_{s_1',s_2'}(\bar{\chi}_{r_2},\bar{\chi}_{r_1};-1/(r_1r_2\tau))$.

We re-emphasize that our general result in \eqref{eq:XiTS},
with exact non-perturbative terms in \eqref{eq:XiNP1}, holds for all 
$s_1,s_2\in \mathbb{C}$, with a resummed perturbative contribution
that is analytic in the upper-half $\tau$-plane.  
Moreover, in \eqref{eq:XiNP2}, we showed that the non-perturbative contribution
has a transseries containing Fricke duals. 
Thus the results of this section and the concrete examples of modular primitives
in the next section originate from a more general notion of a modularity gap, revealed
by \eqref{eq:TSgen}.

\section{Examples of modular primitives}
\label{sec:Ex}

In this section we provide several concrete examples of Lambert series with a character \eqref{eq:LamChi} and their Fricke duals \eqref{eq:Lam2}, focusing first only on real and primitive Dirichlet characters where the story is a little bit neater while discussing later the case of imprimitive characters. 
Presently we discuss only cases for which the parameters $s_1,s_2\in \mathbb{N}$ and the characters $\chi_{r_1},\chi_{r_2}$ are chosen as to satisfy the condition \eqref{eq:CondTrunc} which leads to a terminating perturbative expansion and a direct connection with iterated integrals of twisted Eisenstein series. In Section \ref{sec:Qpoch} we present an application of our general analysis for which the condition \eqref{eq:CondTrunc} is not met.

Since in this section we shall always assume the integrality condition \eqref{eq:CondTrunc}, for ease of presentation we can furthermore restrict our attention only towards the study of the single family of functions $\Xi_{m-1,0}(\chi_{r_1},\chi_{r_2};q) =\mathcal{L}_{m-1}(\chi_{r_1},\chi_{r_2};q)$ for different choices of characters and for different values of $m = s_1-s_2+1\in \mathbb{N}$, provided that $m \equiv \kappa_1 +\kappa_2\,{\rm mod}\,2$.
Thanks to the results of Section \ref{sec:ModularPrim}, the infinite dimensional family of functions $\Xi_{m-1 +s_2 ,s_2}(\chi_{r_1},\chi_{r_2};q)$ with $s_2\in \mathbb{Z}$ can easily be computed directly from the results we shall provide for $\Xi_{m-1,0}(\chi_{r_1},\chi_{r_2};q) =\mathcal{L}_{m-1}(\chi_{r_1},\chi_{r_2};q)$ via the action of a simple integral operator, when $s_2>0$, or differential operator, when $s_2<0$, that is
\begin{equation}
\Xi_{m-1+s_2,s_2}(\chi_{r_1},\chi_{r_2};q) = \begin{cases}
& \int_{0\leq q_1 \leq ...\leq q_{s_2} \leq q}\frac{{\rm d} q_1}{q_1} ... \frac{{\rm d}  q_{s_2}}{q_{s_2}} \mathcal{L}_{m-1}(\chi_{r_1},\chi_{r_2};q_{s_2})\,,\qquad s_2>0\,;\\
& \mathcal{L}_{m-1}(\chi_{r_1},\chi_{r_2};q)\,, \qquad\qquad\qquad\qquad\quad\qquad\qquad s_2=0\,;\\
&  \frac{d^{\tilde{s}_2}}{dq^{\tilde{s}_2}}\mathcal{L}_{m-1}(\chi_{r_1},\chi_{r_2};q)\,, \qquad\qquad\qquad\qquad\qquad\quad\,\, \tilde{s}_2=-s_2>0\,.
\end{cases}
\end{equation}

Strikingly, many of the modular forms here discussed can be represented in terms of linear combinations of eta quotients, i.e. they can be expressed as Laurent polynomials with integer coefficients in $\eta(q^{t})$ for some values of $t \in \mathbb{N}$, where $\eta(q)$ is the standard Dedekind eta function defined in \eqref{eq:eta}.
In Table \ref{tab:SingleTwist} and \ref{tab:DoubleTwist} we present all cases where the holomorphic Eisenstein series twisted by either one or two characters, respectively, can be represented as a single eta quotient.

\subsection{Single real and primitive Dirichlet characters}

We start our discussion by considering the twisted Lambert series \eqref{eq:LamChi} and its Fricke dual generalized Lambert series \eqref{eq:Lam2} for the special cases which originate from iterated integrals of twisted Eisenstein series as discussed in section \ref{sec:ModularPrim}.
Hence, we assume that the parameter $s=m-1 \in \mathbb{N}$ where the integer $m$ has the same parity as the character under consideration so that the perturbative expansions \eqref{eq:Asypert}-\eqref{eq:Asypert2} do terminate after finitely many terms and the complete transseries representations take the easier forms \eqref{eq:TSLTr}-\eqref{eq:TSLtTr}.

In this subsection $\chi_{r}$ denotes a real and primitive character modulo $r$. As reviewed in appendix \ref{sec:DirCh}, we remind the reader that real primitive characters are in a one-to-one correspondence with fundamental discriminants. If $r$ equals a positive fundamental discriminant such as $1,5,8,12,...$ (see the OEIS sequence \href{https://oeis.org/A003658}{$A003658$}) then we have an even real primitive character modulo $r$. Conversely if $r$ equals the absolute value of a negative fundamental discriminant such as $-3,-4,-7,-8,...$ (see the OEIS sequence \href{https://oeis.org/A003657}{$A003657$}), then there exists an odd real primitive character modulo $r$.

\subsubsection*{Case $r=3$ and $s\in \mathbb{N}$ even}

Let us start with the case of the real primitive character modulo $r=3$, i.e. in terms of Conrey labels we consider $\chi_{3,2}(n)$ defined as $\chi_{3,2}(n) = 1 $ if $n\equiv 1 \,\rm{mod}\,3$ and $\chi_{3,2}(n)=-1$ if $n\equiv -1\, \rm{mod}\,3$, and $0$ otherwise, while $s$ is a positive even integer.
The Lambert series \eqref{eq:LamChi} here takes the form
\begin{equation}
\mathcal{L}_s(\chi_{3,2} ; q) = \sum_{n= 1}^\infty  \frac{\chi_{3,2}(n)}{n^s} \frac{q^n}{1-q^n}= \sum_{n=1}^\infty \left[ \frac{1}{(3n-2)^s} \frac{q^{3n-2}}{1-q^{3n-2}} - \frac{1}{(3n-1)^s} \frac{q^{3n-1}}{1-q^{3n-1}} \right]\,.\label{eq:Lchi3}
\end{equation}

As already emphasized, the asymptotic tail \eqref{eq:AsyL} terminates with $k= s-1$ and it reduces to the Laurent polynomial \eqref{eq:Asypert} which for the present case takes the general form:
\begin{equation}
\mathcal{L}_s(\chi_{3,2} ; y) \sim \mathcal{L}^{{\rm Pert}}_s(\chi_{3,2} ; y) = \sum_{k=-1}^{s-1} \frac{(-2\pi y)^k}{k!} \zeta(-k) L(\chi_{3,2},s-k)\,.\label{eq:PertChi3}
\end{equation}
For example we have,
\begin{align}
\mathcal{L}^{\rm{Pert}}_2(\chi_{3,2} ; y) &=\frac{L(\chi_{3,2},3)}{2 \pi  y}-\frac{L(\chi_3,2)}{2}+ \frac{1}{12}   L(\chi_{3,2},1) (2\pi y)\,,\\
\mathcal{L}^{\rm{Pert}}_4(\chi_{3,2} ; y) &=\frac{L\left(\chi _{3,2},5\right)}{2 \pi  y} -\frac{1}{2} L\left(\chi _{3,2},4\right)+ \frac{1}{12}  L\left(\chi _{3,2},3\right)(2\pi  y) -\frac{1}{720}  L\left(\chi _{3,2},1\right)(2\pi y)^3\,,
\end{align}
where the $L$-values can be easily computed from
\begin{equation}
L(\chi_{3,2},s) = 3^{-s} \Big[ \zeta \big(s,\frac{1}{3}\big)-\zeta\big(s,\frac{2}{3}\big)\Big]\,,
\end{equation}
so that
\begin{align}
L(\chi_{3,2},1) &\notag = \frac{\pi }{3 \sqrt{3}}\,,\qquad  L(\chi_{3,2},2) =\frac{1}{9} \left(\psi ^{(1)}\left(\tfrac{1}{3}\right)-\psi ^{(1)}\left(\tfrac{2}{3}\right)\right)\,,\qquad L(\chi_{3,2},3) = \frac{4 \pi ^3}{81 \sqrt{3}}\,,\\
L(\chi_{3,2},4) &= \frac{1}{486} \left(\psi ^{(3)}\left(\tfrac{1}{3}\right)-\psi ^{(3)}\left(\tfrac{2}{3}\right)\right)\,,\qquad L(\chi_{3,2},5) =\frac{4 \pi ^5}{729 \sqrt{3}}\,,
\end{align}
where $\psi^{(n)}(x)$ denotes the $n^{th}$ derivative of the digamma function $\psi(x) = \Gamma'(x)/\Gamma(x)$.

The non-perturbative corrections are encoded in \eqref{eq:TSQM}. For a real primitive character we have $\epsilon(\chi_r)=1$ so that $\epsilon(\chi_{3,2})=1$ and  since the character $\chi_{3,2}$ is odd we find
\begin{align}
\mathcal{L}^{\rm{NP}}_s(\chi_{3,2} ; q) &= e^{i \pi \frac{s}{2}} \frac{y^{s-1}}{\sqrt{3}} \tilde{\mathcal{L}}_{s}\left({\chi}_{3,2}; \frac{1}{ry}\right) = e^{ i\pi \frac{s}{2} } \frac{ y^{s-1}} {\sqrt{3} } \sum_{m=1}^\infty \frac{1}{m^s} \Phi_0(\chi_{3,2}; e^{\frac{-2\pi m }{3 y}})\\
& = e^{ i\pi \frac{s}{2} } \frac{ y^{s-1}} {\sqrt{3} } \sum_{m=1}^\infty \frac{1}{m^s} \frac{\tilde{q}^m}{1+\tilde{q}^m+\tilde{q}^{2m}}\,,
\end{align}
with $\tilde{q} = e^{-\frac{2\pi}{3y}}$ and where we used \eqref{eq:GenToLam} combined with the definition \eqref{eq:NPPhi} for $\Phi_0(\chi_{3,2};x)$ here explicitly given by
\begin{align}
\Phi_0(\chi_{3,2},x)= & \frac{ \sum_{a\in (\mathbb{Z}/3\mathbb{Z})^\times}  \chi_{3,2}(a) x^a}{1-x^3} = \frac{x-x^2}{1-x^3} = \frac{x}{1+x+x^2}\,.\label{eq:Phi032}
\end{align}

For concreteness, for the cases $s=2$ and $s=4$ we have the Fricke relations
\begin{align}
 \sum_{m= 1}^\infty  \frac{\chi_{3,2}(m)}{m^2} \frac{q^m}{1-q^m} =&\notag  \frac{2 \pi ^2}{81 \sqrt{3} y}+\frac{\psi ^{(1)}\left(\tfrac{2}{3}\right)-\psi ^{(1)}\left(\tfrac{1}{3}\right)}{18} +\frac{\pi ^2 y}{18 \sqrt{3}}\\
&\label{eq:mod32s2} -\frac{y}{\sqrt{3}}  \sum_{m= 1}^\infty \frac{1}{m^2} \frac{\tilde{q}^m}{1+\tilde{q}^m+\tilde{q}^{2m}}\,,\\
 \sum_{m= 1}^\infty  \frac{\chi_{3,2}(m)}{m^4} \frac{q^m}{1-q^m}  =&\notag \frac{2 \pi ^4}{729 \sqrt{3} y}+\frac{\psi ^{(3)}\left(\frac{2}{3}\right)-\psi ^{(3)}\left(\frac{1}{3}\right)}{972} +\frac{2 \pi ^4 y}{243 \sqrt{3}}-\frac{\pi ^4 y^3}{270 \sqrt{3}}\\
&\label{eq:mod32s4}  +\frac{ y^{3}} {\sqrt{3} } \sum_{m= 1}^\infty \frac{1}{m^4} \frac{\tilde{q}^m}{1+\tilde{q}^m+\tilde{q}^{2m}}\,,
\end{align}
where again $q=e^{-2\pi y}$ and $\tilde{q} = e^{-\frac{2\pi }{3y}}$.

For $\chi_{3,2}$ and $s=2$ and $s=4$, the differential equations \eqref{eq:IterInt}-\eqref{eq:IterInt2} take the form
\begin{align}
\Big(q \frac{d}{dq}\Big)^2\mathcal{L}_2(\chi_{3,2}; q)  &=G_3(\chi_{3,2};q) = \frac{\eta(q^3)^9}{ \eta(q)^3}\label{eq:L2chi3}\,,\\
\Big(q \frac{d}{dq}\Big)^4 \mathcal{L}_4(\chi_{3,2} ; q) &=G_5(\chi_{3,2};q) = \frac{ a(q)^2 c(q)^3}{27}= \eta(q)^3\eta(q^3)^7\Big( 1+9\frac{\eta(q^9)^3}{\eta(q)^3}\Big)^2 \,,\label{eq:L4chi3}
\end{align}
where $a(q)$ and $c(q)$ are cubic AGM theta functions \cite{Borwein}
\begin{align}
a(q) =  \frac{\big[\eta(q)^3+9 \eta(q^9)^3\big]}{\eta(q^3)}\,,\qquad c(q) = 3   \frac{ \eta(q^3)^3 }{\eta(q)}
\end{align} 
The modular transformations \eqref{eq:mod32s2}-\eqref{eq:mod32s4} are a consequence of the fact that $\mathcal{L}_2(\chi_{3,2}; q) $ is an iterated integral of the $\Gamma_0(3)$ weight-$3$ modular form $G_3(\chi_{3,2};q)$ while $\mathcal{L}_4(\chi_{3,2} ; q)$ is an iterated integral of the $\Gamma_0(3)$ weight-$5$ modular form $G_5(\chi_{3,2};q) $ both of them with Dirichlet character $\chi_{3,2}$.

\subsubsection*{Case $r=4$ and $s\in \mathbb{N}$ even}

We now consider the odd real primitive character modulo $4$ denoted in Conrey labels by $\chi_{4,3}$ with $\chi_{4,3}(n) = 1$ if $n\equiv 1 \,\rm{mod}\,4$ and $\chi_{4,3}(n) = -1$ if $n\equiv -1 \,\rm{mod}\,4$  and $0$ otherwise and analyze the Lambert series
\begin{equation}
\mathcal{L}_s(\chi_{4,3} ; q) = \sum_{n= 1}^\infty  \frac{\chi_{4,3}(n)}{n^s} \frac{q^n}{1-q^n}= \sum_{n=1}^\infty \left[ \frac{1}{(4n-3)^s} \frac{q^{4n-3}}{1-q^{4n-3}} - \frac{1}{(4n-1)^s} \frac{q^{4n-1}}{1-q^{4n-1}} \right]\,,\label{eq:Lchi4}
\end{equation}
with $s\in \mathbb{N}$ an even integer.
The perturbative expansion \eqref{eq:Asypert} contains finitely many terms
\begin{equation}
\mathcal{L}_s(\chi_{4,3} ; y) \sim \mathcal{L}^{{\rm Pert}}_s(\chi_{4,3} ; y) = \sum_{k=-1}^{s-1} \frac{(-2\pi y)^k}{k!} \zeta(-k) L(\chi_{4,3},s-k)\,,\label{eq:PertChi4}
\end{equation}
which can be computed from
\begin{equation}
L(\chi_{4,3},s) = 4^{-s} \left(\zeta \left(s,\frac{1}{4}\right)-\zeta \left(s,\frac{3}{4}\right)\right)\,.
\end{equation}
In particular we have
\begin{align}
L(\chi_{4,3},1) &\notag = \frac{\pi}{4}\,,\qquad  L(\chi_{4,3},2) =C \,,\qquad L(\chi_{4,3},3) = \frac{\pi^3}{32}\,,\\
L(\chi_{4,3},4) &= \frac{\psi ^{(3)}\left(\frac{1}{4}\right)-\psi ^{(3)}\left(\frac{3}{4}\right)}{1536}\,,\qquad L(\chi_{4,3},5) =\frac{5 \pi ^5}{1536}\,,
\end{align}
with $C$ Catalan's constant.

The non-perturbative corrections are encoded in \eqref{eq:TSQM} using $\epsilon(\chi_{4,3})=1$ and $\kappa=1$, 
\begin{align}
\mathcal{L}^{\rm{NP}}_s(\chi_{4,3} ; q) &= e^{i \pi \frac{s}{2}} \frac{y^{s-1}}{2} \tilde{\mathcal{L}}_{s}\left({\chi}_{4,2}; \frac{1}{4y}\right) = e^{ i\pi \frac{s}{2} } \frac{ y^{s-1}} {2 } \sum_{m=1}^\infty \frac{1}{m^s} \Phi_0(\chi_{4,3}; e^{\frac{-2\pi m }{4 y}})\\
& = e^{ i\pi \frac{s}{2} } \frac{ y^{s-1}} {2 } \sum_{m=1}^\infty \frac{1}{m^s} \frac{\tilde{q}^m}{1+\tilde{q}^{2m}}\,,
\end{align}
where $\tilde{q} = e^{-\frac{2\pi}{4y}}$ and we used
\begin{align}
\Phi_0(\chi_{4,3};x)= & \frac{ \sum_{a\in (\mathbb{Z}/4\mathbb{Z})^\times}  \chi_{4,3}(a) x^a}{1-x^4} = \frac{x-x^3}{1-x^4} = \frac{x}{1+x^2}\,.\label{eq:Phi043}
\end{align}

Specialising to $s=2$ and $s=4$ we have
\begin{align}
\sum_{m= 1}^\infty  \frac{\chi_{4,3}(m)}{m^2} \frac{q^m}{1-q^m} = &\notag \frac{L(\chi_{4,3},3)}{2 \pi  y}-\frac{L(\chi_{4,3},2)}{2}+ \frac{1}{12}   L(\chi_{4,3},1) (2\pi y)\\
&- \frac{ y}{2}\sum_{m= 1}^\infty \frac{1}{m^2}\frac{\tilde{q}^m}{ 1 +\tilde{q}^{2m}} \,, \\
\sum_{m= 1}^\infty  \frac{\chi_{4,3}(m)}{m^4} \frac{q^m}{1-q^m}=&\notag\frac{L\left(\chi _{4,3},5\right)}{2 \pi  y} -\frac{1}{2} L\left(\chi _{4,3},4\right)+ \frac{1}{12}  L\left(\chi _{4,3},3\right)(2\pi  y) -\frac{1}{720}  L\left(\chi _{4,3},1\right)(2\pi y)^3\\
&+ \frac{ y^3}{2}\sum_{m= 1}^\infty \frac{1}{m^4} \frac{\tilde{q}^m}{ 1+\tilde{q}^{2m}}\,.
\end{align}

From the differential equations \eqref{eq:IterInt}-\eqref{eq:IterInt2} we find,
\begin{align}
\Big(q \frac{d}{dq}\Big)^2 \mathcal{L}_2(\chi_{4,3} ; q) &\label{eq:L2chi4}=G_3(\chi_{4,3};q) =  \frac{\eta(q^2)^6\, \eta(q^4)^4}{ \eta(q)^4 }\,,\\
\Big(q \frac{d}{dq}\Big)^4 \mathcal{L}_4(\chi_{4,3} ; q) &\label{eq:L4chi4}=G_5(\chi_{4,3};q) =  \eta(q^2)^2 \, \eta(q^4)^4  \Big[ \eta(q)^4 + 20  \frac{\eta(q^4)^8 }{ \eta(q)^4}\Big]\,,
\end{align}
so that we can interpret $\mathcal{L}_2(\chi_{4,3} ; q)$ and $\mathcal{L}_4(\chi_{4,3} ; q)$ as respectively iterated integrals of $\Gamma_0(4)$ modular forms of weight $3$ and $5$ respectively with Dirichlet character $\chi_{4,3}$.

\subsubsection*{Case $r=5$ and $s\in \mathbb{N}$ odd.}

 For completeness we also consider the case of the primitive real even character $\chi_{5,4}$ in Conrey labels, i.e.
\begin{center}
\begin{tabular}{ |c|c|c|c|c|} 
\hline
 $n\,\rm{mod}\,5$ & 1 & 2 & 3 & 4 \\
\hline
\multirow{1}{4em}{$\chi_{5,4}(n)$} &1 & -1 &-1& 1\\
\hline
\end{tabular}
\end{center}

Note that since now the character is even to obtain a perturbative expansion which terminates after finitely many terms we must chose $s\in\mathbb{N}$ odd so that the perturbative expansion \eqref{eq:Asypert} reduces to:
\begin{equation}
\mathcal{L}_s(\chi_{5,4} ; y) \sim \mathcal{L}^{{\rm Pert}}_s(\chi_{5,4} ; y) = \sum_{k=-1}^{s-2} \frac{(-2\pi y)^k}{k!} \zeta(-k) L(\chi_{5,4},s-k)\,,\label{eq:PertChi5}
\end{equation}
where the $L$-values can be computed from
\begin{equation}
L(\chi_{5,4},s) = 5^{-s} \left(\zeta \left(s,\frac{1}{5}\right)-\zeta \left(s,\frac{2}{5}\right)-\zeta \left(s,\frac{3}{5}\right)+\zeta \left(s,\frac{4}{5}\right)\right)\,.
\end{equation} 

Again, for the non-perturbative corrections we use \eqref{eq:TSQM} specialized to $\epsilon(\chi_{5,4})=1$ and $\kappa=0$, 
\begin{align}
\mathcal{L}^{\rm{NP}}_s(\chi_{5,4} ; q) &= e^{i \pi \frac{s-1}{2}} \frac{y^{s-1}}{5} \tilde{\mathcal{L}}_{s}\left({\chi}_{5,4}; \frac{1}{5y}\right) = e^{ i\pi \frac{s-1}{2} } \frac{ y^{s-1}} {\sqrt{5} } \sum_{m=1}^\infty \frac{1}{m^s} \Phi_0(\chi_{5,4}; e^{\frac{-2\pi m }{5 y}})\\
& = e^{ i\pi \frac{s-1}{2} } \frac{ y^{s-1}} {2 } \sum_{m=1}^\infty \frac{1}{m^s} \frac{\tilde{q}^m(\tilde{q}^m+1)(\tilde{q}^{m}-1)^2}{1-\tilde{q}^{5m}}\,,
\end{align}
where $\tilde{q} = e^{-\frac{2\pi}{5y}}$ and we used
\begin{align}
\Phi_0(\chi_{5,2};x)= & \frac{ \sum_{a\in (\mathbb{Z}/5\mathbb{Z})^\times}  \chi_{5,2}(a) x^a}{1-x^5} 
= \frac{x(1+x)(1-x)^2}{(1-x^5)} = \frac{x(1-x^2)}{1+x+x^2+x^3+x^4}\,.\label{eq:Phi052}
\end{align}

Specializing to $s=1$ and $s=3$ we have
\begin{align}
\mathcal{L}_1(\chi_{5,4} ; q) = &\notag \frac{L(\chi_{5,4},2)}{2 \pi  y}+\frac{ 1} {\sqrt{5} } \sum_{m= 1}^\infty \frac{1}{m}\frac{\tilde{q}^m(\tilde{q}^m+1)(\tilde{q}^{m}-1)^2}{1-\tilde{q}^{5m}} \,, \\
\mathcal{L}_3(\chi_{5,4} ; q) = &\notag \frac{L(\chi_{5,4},4)}{2 \pi  y}-\frac{L(\chi_{5,4},3)}{2}+ \frac{1}{12}   L(\chi_{5,4},2) (2\pi y)\\
&-\frac{ y^{2}} {\sqrt{5} } \sum_{m= 1}^\infty \frac{1}{m^3} \frac{\tilde{q}^m(\tilde{q}^m+1)(\tilde{q}^{m}-1)^2}{1-\tilde{q}^{5m}} \,,
\end{align}

The relevant $L$-values are:
\begin{align}
L\left(\chi _{5,4},2\right) &\notag  =\frac{4 \pi ^2}{25 \sqrt{5}} \,,\qquad L\left(\chi _{5,4},3\right)=\frac{\left(-\psi ^{(2)}\left(\frac{1}{5}\right)+\psi ^{(2)}\left(\frac{2}{5}\right)+\psi ^{(2)}\left(\frac{3}{5}\right)-\psi ^{(2)}\left(\frac{4}{5}\right)\right)}{250}\,,  \\
 L\left(\chi _{5,4},4\right) &\notag  =  \frac{8 \pi ^4}{375 \sqrt{5}}\,.
\end{align}

We find that
\begin{align}
q\frac{d}{dq}  \mathcal{L}_1(\chi_{5,4} ; q) &=G_2(\chi_{5,4};q) =  \frac{\eta(q^5)^5}{\eta(q)}\,,\label{eq:L1chi5}\\
\left(q \frac{d}{dq}\right)^3  \mathcal{L}_3(\chi_{5,4} ; q) &=G_4(\chi_{5,4};q)  \,,
\end{align}
i.e. we can interpret $ \mathcal{L}_1(\chi_{5,4} ; q) $ as an iterated integral of a weight-$2$ modular form while $\mathcal{L}_3(\chi_{5,4} ; q)$ is an iterated integral of a weight-$4$ modular form both with respect to the congruence subgroup $\Gamma_0(5)$  and character $\chi_{5,4}$.

The twisted Eisenstein series $G_4(\chi_{5,4};q) $ is not quite a linear combination of eta quotients but we can show that it satisfies the Borwein-like \cite{Borwein} equation
\begin{equation}
G_4(\chi_{5,4};q)^2 
 = \big[\eta(q) \eta(q^5)\big]^8 \left[ 1 + 22 \left(\frac{\eta(q^5)}{\eta(q)} \right)^6 + 125 \left(\frac{\eta(q^5)}{\eta(q)}\right)^{12} \right]\,.
\end{equation}

\subsubsection*{Other examples of modular primitives}

Given the previous discussion it should appear clear how to extend our discussion to other examples.
We list here some interesting results regarding modular primitives to the differential equation \eqref{eq:IterInt}.
\vspace{0.3cm}

\textbf{Case $r=7$ and $s\in \mathbb{N}$ even}
\vspace{0.3cm}

For the primitive real odd character $\chi_{7,6}$, i.e.
\begin{center}
\begin{tabular}{ |c|c|c|c|c|c|c|} 
\hline
 $n\,{\rm mod}\,7$ & 1 & 2 & 3 & 4 & 5 & 6 \\
\hline
\multirow{1}{4em}{$\chi_{7,6}(n)$} &1 & 1 &-1& 1 & -1 & -1\\
\hline
\end{tabular}
\end{center}

we find the modular primitives
\begin{align}
\left(q \frac{d}{dq} \right) \mathcal{L}_2(\chi_{7,6} ; q) &=G_3(\chi_{7,6};q) = \eta (q )^3 \eta (q^7 )^3+ \frac{8 \eta (q^7 )^7}{\eta (q )}\,,\label{eq:L2chi7}\\
\left( q \frac{d}{dq}\right)^4  \mathcal{L}_4(\chi_{7,6} ; q) &=G_5(\chi_{7,6};q) =  A(q)^2 \big[\eta(q)\eta(q^7)\big]^3(1+16B(q))  \,,\label{eq:L4chi7}
\end{align}
where we defined
\begin{align}
A(q) &\coloneqq \frac{\big[ \eta(q)\eta(q^7)\big]^3+4\big[\eta(q^2)\eta(q^{14})\big]^3}{\eta(q)\eta(q^2)\eta(q^7)\eta(q^{14})} = \theta_3(q) \theta_3(q^7) + \theta_2(q) \theta_2(q^7)\,,\\
B(q) & \coloneqq \frac{\eta(q^7)^4}{\eta(q)^4}\,.
\end{align}

Here $\theta_2(q)$ and $\theta_3(q)$ denote the standard Jacobi theta functions defined by
\begin{align}
\theta_2(q) &= \sum_{n=-\infty}^{\infty} q^{ - (n+\frac{1}{2})^2} =  q^{\frac{1}{4}}\big[ 2 +2 q^2+ 2q^6 + O(q^{12})\big]\,,\\
\theta_3(q) &= \sum_{n=-\infty}^{\infty} q^{ -n^2} =  1 +2 q + 2q^4 + O(q^{9})\,.
\end{align}

\textbf{Case $r=8$}
\vspace{0.3cm}

The case of modulus $r=8$ is special since for all $m=8k$ we have two real primitive characters modulo $m$, namely the even one $\chi_r(a) = \left( \frac{m}{a}\right)$ and the odd one $\chi_r(a) = \left( \frac{-m}{a}\right)$ where $\left( \frac{D}{\bullet}\right)$ denotes the Kronecker symbol, see the discussion in appendix \ref{sec:DirCh}.

Considering for $r=8$ the odd character $\chi_{8,3}(a) = \left( \frac{-8}{a}\right)$ in Conrey labels, we find for $s=2$:
\begin{align}
\left( q \frac{d}{dq}\right)^2 \mathcal{L}_2(\chi_{8,3},q) &= G_3(\chi_{8,3},q) = \eta(q)^2\eta(q^2) \eta(q^4) \eta(q^8)^2 (1+6h(q)) \,,\\
\left( q \frac{d}{dq}\right)^2\tilde{ \mathcal{L}}_2(\chi_{8,3},q) & = E_{3}(\chi_{8,3};q) -  \frac{L(\chi_{8,3},-2)}{2} =  \frac{1}{2} \Big[3-\eta(q)^2 \eta(q^2) \eta(q^4) \eta(q^8)^2 \big(16 +\frac{3}{h(q)}\big) \Big]\,,
\end{align}
with $ \frac{L(\chi_{8,3},-2)}{2}=-3$ and where $h(q)$ denotes the Hauptmodul at level 8, see e.g.~Table 1 of \cite{AB},
\begin{equation}
h(q) = \Big(\frac{\eta(q^2)}{\eta(q^4)}\Big)^2 \Big(\frac{ \eta(q^8)}{\eta(q)}\Big)^4\,.
\end{equation}

Considering instead $r=8$ and the even character $\chi_{8,5}(a)= \left( \frac{8}{a}\right)$ in Conrey labels, we find for $s=1$:
\begin{align}
q\frac{d}{dq} \mathcal{L}_1(\chi_{8,5};q) &= G_2(\chi_{8,5};q) = \frac{\eta(q^2)^3\eta(q^4) \eta(q^8)^2}{\eta(q)^2} \,,\\
q\frac{d}{dq}\tilde{ \mathcal{L}}_1(\chi_{8,5};q) & = E_{2}(\chi_{8,5};q)  -  \frac{L(\chi_{8,5},-1)}{2} =  \Big(1-\frac{\eta(q)^2 \eta(q^2) \eta(q^4)^3}{\eta(q^8)^2}\Big) \,.
\end{align}
having used $L(\chi_{8,5},-1) = -1$.
\vspace{0.3cm}

\textbf{Case $r=12$ and $s\in\mathbb{N}$ odd}
\vspace{0.3cm}

The only primitive character modulo $12$ is the even character $\chi_{12,11}$ for which we find
\allowdisplaybreaks{
\begin{align}
q\frac{d}{dq} \mathcal{L}_1(\chi_{12,11}; q) &= G_2(\chi_{12,11};q) = \frac{\eta(q^2)^2\eta(q^3)^2\eta(q^4)\eta(q^{12})}{\eta(q)^2} \,,\\
q\frac{d}{dq}\tilde{ \mathcal{L}}_1(\chi_{12,11}; q) & =  E_2(\chi_{12,11};q)  -  \frac{L(\chi_{12,11},-1)}{2}  = 1-\frac{\eta(q) \eta(q^3) \eta(q^4)^2\eta(q^6)^2}{\eta(q^{12})^2}\,,
\end{align}}
having used  $L(\chi_{12,11},-1)=-2$.
\vspace{0.2cm}

\subsection*{Eta quotients representation}
From the many examples here presented, we see that for different values of the modular weight $m=s+1$ and different choices of primitive character $\chi_r$, the corresponding twisted Eisenstein series $G_m(\chi_r; q)$ and $E_m(\chi_r;q)$ can be expressed as Laurent polynomials with integer coefficients in $\eta(q^{t})$ for some values of $t \in \mathbb{N}$, where $\eta(q)$  is the standard Dedekind $\eta$-function,
\begin{equation}
\eta(q) \coloneqq q^{\frac{1}{24}} \prod_{n=1}^\infty (1-q^n)\,,\label{eq:eta}
\end{equation}
with $q\in\mathbb{C}$ and $|q|<1$. In Table \ref{tab:SingleTwist}, we present all of the instances where for a given level $N\in \mathbb{N}$, i.e. with respect to the congruence subgroup $\Gamma_0(N)$, modular weight $m$, and characters $\chi_{r_1}(a) = \left(\frac{D_1}{a}\right)$ and $\chi_{r_2}(a) = \left(\frac{D_2}{a}\right)$, with $D_1,D_2 \in \mathbb{Z}$ two fundamental discriminants, the corresponding twisted Eisenstein series $G_m(\chi_{r_1},\chi_{r_2};q)$  can be written as a single eta quotient, i.e. as
\begin{equation}
G_m(\chi_{r_1},\chi_{r_2};q) = \prod_{j=1}^K \eta(q^{t_j})^{r_j}\,,
\end{equation}
for some $t_j\in \mathbb{N}\,, r_j\in \mathbb{Z}$ with $j=1,...,K$ and $K\in \mathbb{N}$. Note that in Table \ref{tab:SingleTwist} we use the short-hand notation $\eta_t^r \coloneqq \eta(q^t)^r$.

\begin{table}[h]\centering
\begin{tabular}{|r|r|r|r|c|c|c|}\hline
$N$ & $m$ & $D_1$ & $D_2$ & Numerator & Denominator & OEIS Sequence\\ \hline 
3 & 3 & -3 & 1 & $\eta_{3}^{9} $ & $\eta_{1}^{3} $ & \oeis{A106402} \\
3 & 3 & 1 & -3 & $\eta_{1}^{9} $ & $\eta_{3}^{3} $ & \oeis{A109041} \\
4 & 1 & 1 & -4 & $\eta_{2}^{10} $ & $\eta_{1}^{4} \eta_{4}^{4} $ & \oeis{A004018} \\
4 & 3 & -4 & 1 & $\eta_{2}^{6} \eta_{4}^{4} $ & $\eta_{1}^{4} $ & \oeis{A050470} \\
4 & 3 & 1 & -4 & $\eta_{1}^{4} \eta_{2}^{6} $ & $\eta_{4}^{4} $ & \oeis{A120030} \\
5 & 2 & 5 & 1 & $\eta_{5}^{5} $ & $\eta_{1}^{\phantom{1}} $ & \oeis{A053723} \\
5 & 2 & 1 & 5 & $\eta_{1}^{5} $ & $\eta_{5}^{\phantom{1}} $ & \oeis{A109064} \\
8 & 1 & 1 & -8 & $\eta_{2}^{3} \eta_{4}^{3} $ & $\eta_{1}^{2} \eta_{8}^{2} $ &\oeis{A033715} \\
8 & 2 & 8 & 1 & $\eta_{2}^{3} \eta_{4}^{\phantom{1}} \eta_{8}^{2} $ & $\eta_{1}^{2} $ & \oeis{A124340} \\
8 & 2 & 1 & 8 & $\eta_{1}^{2} \eta_{2}^{\phantom{1}} \eta_{4}^{3} $ & $\eta_{8}^{2} $ & \oeis{A131999} \\
12 & 2 & 12 & 1 & $\eta_{2}^{2} \eta_{3}^{2} \eta_{4}^{\phantom{1}} \eta_{12}^{\phantom{1}} $ & $\eta_{1}^{2} $ & \oeis{A124815} \\
12 & 2 & 1 & 12 & $\eta_{1}^{\phantom{1}} \eta_{3}^{\phantom{1}} \eta_{4}^{2} \eta_{6}^{2} $ & $\eta_{12}^{2} $ & \oeis{A109039} \\
15 & 1 & 1 & -15 & $\eta_{3}^{2} \eta_{5}^{2} $ & $\eta_{1}^{\phantom{1}} \eta_{15}^{\phantom{1}} $ & \oeis{A123864} \\
20 & 1 & 1 & -20 & $\eta_{2}^{\phantom{1}} \eta_{4}^{\phantom{1}} \eta_{5}^{\phantom{1}} \eta_{10}^{\phantom{1}} $ & $\eta_{1}^{\phantom{1}} \eta_{20}^{\phantom{1}} $ & \oeis{A124233} \\
24 & 1 & 1 & -24 & $\eta_{2}^{\phantom{1}} \eta_{3}^{\phantom{1}} \eta_{8}^{\phantom{1}} \eta_{12}^{\phantom{1}} $ & $\eta_{1}^{\phantom{1}} \eta_{24}^{\phantom{1}} $ & \oeis{A000377} \\
\hline\end{tabular}
\caption{
Numerators and denominators of eta quotients that are singly twisted Eisenstein series,
where $N$ is the level, $m$ is the weight, $D_1$ and $D_2$ are fundamental discriminants specifying 
the twists, and the final column gives the associated entry in the On-line Encyclopedia of Integer Sequences (\href{https://oeis.org/}{OEIS}).}
\label{tab:SingleTwist}
\end{table}

\subsection{Single real and imprimitive Dirichlet characters}
\label{sec:ExImp}

In this subsection we present some examples of $q$-series \eqref{eq:LamChi}-\eqref{eq:Lam2} for which the Dirichlet characters are imprimitive while the parameter $s=m-1 \in \mathbb{N}$ with $m$ having the same parity as the character under consideration. Once again, with this choice of parameters we find perturbative expansions \eqref{eq:Asypert}-\eqref{eq:Asypert2} that terminate after finitely many terms.
From our discussion in section \ref{sec:NPImp} we see from the transseries representation \eqref{eq:TSQMImp} that imprimitive characters  lead to particular linear combinations of the transseries \eqref{eq:TSLTr}-\eqref{eq:TSLtTr} corresponding to modular primitives of twisted Eisenstein series. 

\subsubsection*{Case $r=6$ and $s\in \mathbb{N}$ even}

Let us start by considering the odd imprimitive character modulo $6$: $\chi_{6,5}(n) = 1$ if $n\equiv 1\,\mbox{mod}\,6$ and $\chi_{6,5}(n) = -1$ if $n\equiv -1\,\mbox{mod}\,6$ and $0$ otherwise, which is induced by the primitive character $\chi_{3,2}$ already analyzed, thus leading to the $L$-functions relation,
\begin{equation}
L(\chi_{6,5},s) = L(\chi_3,s) (1+2^{-s})\,.\label{eq:Lvalues63}
\end{equation}

We then consider the twisted Lambert series
\begin{equation}
\mathcal{L}_s(\chi_{6,5} ; q) = \sum_{n= 1}^\infty  \frac{\chi_{6,5}(n)}{n^s} \frac{q^n}{1-q^n}= \sum_{n=1}^\infty \left[ \frac{1}{(6n-5)^s} \frac{q^{6n-5}}{1-q^{6n-5}} - \frac{1}{(6n-1)^s} \frac{q^{6n-1}}{1-q^{6n-1}} \right]\,,\label{eq:Lchi6}
\end{equation}

Focusing on the non-perturbative corrections in \eqref{eq:TSQMImp} we deduce 
\begin{align}
\mathcal{L}^{\rm{NP}}_s(\chi_{6,5} ; q) &\notag = e^{ i \frac{\pi s}{2}}    \frac{y^{s-1}}{\sqrt{3}}  \sum_{d\vert 6} \mu(d)\frac{ \chi_{3,2}(d)}{d}
{\tilde{\mathcal{L}}}_{s}\left(\chi_{3,2};\frac{1}{3dy}\right)\\
&\notag= e^{ i \frac{\pi s}{2}} \frac{ y^{s-1}}{\sqrt{3}} \sum_{m= 1}^\infty  \frac{1}{m^s} \left[ \Phi_0(\chi_{3,2}; e^{- \frac{2\pi m }{3 y} })+\frac{1}{2} \Phi_0(\chi_{3,2}; e^{- \frac{2\pi m }{6 y} })\right] \\
&= e^{ i \frac{\pi s}{2}} \frac{ y^{s-1}}{2 \sqrt{3}}\sum_{m= 1}^\infty \frac{1}{m^s} \frac{\tilde{q}^m}{ 1-\tilde{q}^m +\tilde{q}^{2m}}\,,\label{eq:NPchi6}
\end{align}
where $\tilde{q} = e^{-\frac{2\pi}{6y}}$ and used the previously computed expression \eqref{eq:Phi032} for $\Phi_0(\chi_{3,2};x )$.

Putting together the perturbative truncating series  \eqref{eq:AsyL} and the non-perturbative corrections \eqref{eq:NPchi6} for the particular cases $s=2$ (which is exactly the case considered in \cite{BV}) and $s=4$ we find
\begin{align}
\sum_{m= 1}^\infty  \frac{\chi_{6,5}(m)}{m^2} \frac{q^m}{1-q^m} = &\notag \frac{L(\chi_{6,5},3)}{2 \pi  y}-\frac{L(\chi_{6,5},2)}{2}+ \frac{1}{12}   L(\chi_{6,5},1) (2\pi y)\\
&- \frac{ y}{2 \sqrt{3}}\sum_{m= 1} \frac{1}{m^2}\frac{\tilde{q}^m}{ 1-\tilde{q}^m +\tilde{q}^{2m}} \,, \\
\sum_{m= 1}^\infty  \frac{\chi_{6,5}(m)}{m^4} \frac{q^m}{1-q^m} =&\notag\frac{L\left(\chi_{6,5},5\right)}{2 \pi  y} -\frac{1}{2} L\left(\chi_{6,5},4\right)+ \frac{1}{12}  L\left(\chi _{6,5},3\right)(2\pi  y) -\frac{1}{720}  L\left(\chi _{6,5},1\right)(2\pi y)^3\\
&+ \frac{ y^3}{2 \sqrt{3}}\sum_{m= 1} \frac{1}{m^4} \frac{\tilde{q}^m}{ 1-\tilde{q}^m +\tilde{q}^{2m}}\,,
\end{align}
where the particular $L$-values can be obtained from \eqref{eq:Lvalues63} and are given by
\begin{align}
L(\chi_{6,5},1) &\notag = \frac{\pi }{2 \sqrt{3}}\,,\qquad  L(\chi_{6,5},2) =  \frac{\psi ^{(1)}\left(\frac{1}{6}\right)-\psi ^{(1)}\left(\frac{5}{6}\right)}{36}\,,\qquad L(\chi_{6,5},3) = \frac{\pi ^3}{18 \sqrt{3}}\,,\\
L(\chi_{6,5},4) &= \frac{\psi ^{(3)}\left(\frac{1}{6}\right)-\psi ^{(3)}\left(\frac{5}{6}\right)}{7776}\,,\qquad L(\chi_{6,5},5) =\frac{11 \pi ^5}{1944 \sqrt{3}}\,.
\end{align}

Similar to the modular primitives \eqref{eq:L2chi3}-\eqref{eq:L4chi3} for the primitive character $\chi_{3,2}$ previously discussed, for the present case of the imprimitive character $\chi_{6,5}$ induced by $\chi_{3,2}$ we find
\begin{align}
 \Big( q\frac{d}{dq}\Big)^2 \mathcal{L}_2(\chi_{6,5} ; q) &=G_3(\chi_{6,5};q) =  \frac{c(q)^3 + c(q^2)^3}{ 27}\,.\label{eq:L2chi6}
\end{align}
Note that from \eqref{eq:L2chi3} we have that $\mathcal{L}_2(\chi_{3,2} ; q)$ is a modular primitive of
\begin{equation}
G_3(\chi_{3,2};q) =  \frac{c(q)^3}{27} = \frac{\eta(q^3)^9}{ \eta(q)^3}\,,
\end{equation}
 so that from \eqref{eq:L2chi3} we deduce
\begin{equation}
 \Big( q\frac{d}{dq}\Big)^2 \mathcal{L}_2(\chi_{6,5} ; q)  =  \Big( q\frac{d}{dq}\Big)^2\left[ \mathcal{L}_2(\chi_{3,2} ; q)+ 2^{-2} \mathcal{L}_2(\chi_{3,2} ; q^2) \right]= G_3(\chi_{3,2};q)+2^{-2}G_3(\chi_{3,2};q^2)\,.
\end{equation}

Similarly, $ \mathcal{L}_4(\chi_{6,5} ; q)$  is a modular primitive of
\begin{align}
 \Big( q\frac{d}{dq}\Big)^4 \mathcal{L}_4(\chi_{6,5} ; q) &= G_5(\chi_{6,5};q) = G_5(\chi_{3,2};q) + 2^{-4}G_5(\chi_{3,2};q^2)\,,\label{eq:L4chi6}
\end{align}
with 
\begin{equation}
G_5(\chi_{3,2};q) =  \frac{ a(q)^2 c(q)^3}{27}\,.
\end{equation}
Note in particular from \eqref{eq:L4chi3} that $\mathcal{L}_4(\chi_{3,2} ; q)$ is a modular primitive of $G_5(\chi_{3,2};q) $ so that from \eqref{eq:L4chi6} we deduce
\begin{equation}
 \mathcal{L}_4(\chi_{6,5} ; q) =  \mathcal{L}_4(\chi_{3,2} ; q)+ 2^{-4} \mathcal{L}_4(\chi_{3,2} ; q^2) \,.
\end{equation}
In direct analogy with the $L$-function identity \eqref{eq:Lvalues63}, for general $s$ we have
\begin{equation}
 \mathcal{L}_s(\chi_{6,5} ; q) =  \mathcal{L}_s(\chi_{3,2} ; q)+ 2^{-s} \mathcal{L}_s(\chi_{3,2} ; q^2) \,.\label{eq:chi6chi3}
\end{equation}
Alternatively, we could have deduced this identity by noting that at the level of characters we have
\begin{equation}
\chi_{6,5}(n) = \chi_{3,2}(n) + \chi_{3,2}(n/2)\,,
\end{equation}
where we set to zero any character evaluated on non-integers, hence by substituting the above expression  in \eqref{eq:LamChi} we obtain \eqref{eq:chi6chi3}.

From this simple example, we see that for imprimitive characters we do not really find a much different story and 
we obtain finite linear combinations of the the results previously discussed for primitive characters.

\subsection{Double real and primitive Dirichlet characters}

We now wish to consider the more general examples defined in equation \eqref{eq:L2chi} of Lambert series $\mathcal{L}_s(\chi_{r_1},\chi_{r_2};q) = \Xi_{s,0}(\chi_{r_1},\chi_{r_2};q)$, where the parameter $s\in \mathbb{N}$ and the characters $\chi_{r_1},\chi_{r_2}$ satisfy the condition  \eqref{eq:CondTrunc}.
As discussed in Section \ref{sec:ModularPrim}, this class of Lambert series can be understood as being an $s$-fold iterated integral \eqref{eq:XiItInt} of the doubly twisted Eisenstein series $G_{s+1}(\chi_{r_1},\chi_{r_2};q)$.
In this case the asymptotic perturbative expansion near $q\to 1^-$ does terminate after finitely many terms and the transseries representation \eqref{eq:TSgenRed} becomes an integrated version of the Fricke inversion formula \eqref{eq:FrickeG2chi} satisfied by the doubly twisted Eisenstein series.

\subsubsection*{Case $r_1=r_2=3$ and $s\in \mathbb{N}$ odd}

Let us start with the case $\chi_{r_1}=\chi_{r_2} = \chi_{3,2}$. For the generalized Lambert series $\mathcal{L}_s(\chi_{r_1},\chi_{r_2};q) $ to satisfy the condition \eqref{eq:CondTrunc}  we must require $s\in \mathbb{N}$ odd since both characters are odd.

The generalized Lambert series \eqref{eq:L2chi}-\eqref{eq:GenLamPhi} here takes the form
\begin{equation}
\mathcal{L}_s(\chi_{3,2},\chi_{3,2} ; q) = \sum_{n= 1}^\infty  \frac{\chi_{3,2}(n)}{n^s} \Phi_0(\chi_{3,2};q^n)= \sum_{n=1}^\infty \frac{\chi_{3,2}(n)}{n^s} \frac{q^n}{1+q^n+q^{2n}}\,,\label{eq:Lchi3232chi}
\end{equation}
having used \eqref{eq:Phi032}.
Its asymptotic perturbative expansion as $y\to0^+$ can be computed directly from \eqref{eq:XiAsyPgen} and crucially it terminates after $s-1$ terms, i.e.
\begin{equation}
\mathcal{L}_s(\chi_{3,2},\chi_{3,2} ; y) \sim \mathcal{L}^{{\rm Pert}}_s(\chi_{3,2},\chi_{3,2} ; y) = \sum_{k=0}^{s-1} \frac{(-2\pi y)^k}{k!}  L(\chi_{3,2},s-k) L(\chi_{3,2},-k)\,.\label{eq:PertChi32chi}
\end{equation}

For the non-perturbative corrections we use \eqref{eq:TSgenRed} specialized to $\epsilon(\chi_{3,2})=1$ and $\kappa_1+\kappa_2=2$ and $(s_1,s_2)=(s,0)$, 
\begin{align}
\mathcal{L}^{\rm{NP}}_s(\chi_{3,2},\chi_{3,2} ; q) &=  -(3iy)^{s-1} \mathcal{L}_{s}\left({\chi}_{3,2},\chi_{3,2}; \tilde{q}\right)\,,
\end{align}
where $\tilde{q} = e^{-\frac{2\pi}{9y}}$.
We then deduce the Fricke inversion identity
\begin{align}
\sum_{n=1}^\infty \frac{\chi_{3,2}(n)}{n^s} \frac{q^n}{1+q^n+q^{2n}} =&\notag \sum_{k=0}^{s-1} \frac{(-2\pi y)^k}{k!}  L(\chi_{3,2},s-k) L(\chi_{3,2},-k)\\
&  -(3iy)^{s-1}\sum_{n=1}^\infty \frac{\chi_{3,2}(n)}{n^s} \frac{\tilde{q}^n}{1+\tilde{q}^n+\tilde{q}^{2n}}\,.
\end{align}

From \eqref{eq:IterIntGen} we find for example the modular primitives
\begin{align}
\left(q \frac{d}{dq}\right) \mathcal{L}_1(\chi_{3,2},\chi_{3,2};q) &= G_2(\chi_{3,2},\chi_{3,2};q) = \frac{[\eta(q)\eta(q^9)]^3}{\eta(q^3)^2}\,,\\
\left(q \frac{d}{dq}\right)^3 \mathcal{L}_3(\chi_{3,2},\chi_{3,2};q) &= G_4(\chi_{3,2},\chi_{3,2};q) =\frac{\eta(q)^9\eta(q^9)^3-27\eta(q)^3\eta(q^9)^9}{\eta(q^3)^4}\,,
\end{align}
where the right hand sides are manifestly modular forms of weight $2$ and $4$ respectively for the congruence subgroup $\Gamma_0(9)$ and with Dirichlet character $\chi_{3,1} =\chi_{3,2} \chi_{3,2}$.

\subsubsection*{Case $r_1=3$, $r_2=4$ and $s\in \mathbb{N}$ odd}

We consider now $\chi_{r_1} = \chi_{3,2}$ and $\chi_{r_2}=\chi_{4,3}$. For the generalized Lambert series $\mathcal{L}_s(\chi_{r_1},\chi_{r_2};q) $ to satisfy the condition \eqref{eq:CondTrunc}  we must require $s\in \mathbb{N}$ odd since both characters are odd.

The generalized Lambert series \eqref{eq:L2chi}-\eqref{eq:GenLamPhi} here takes the form
\begin{equation}
\mathcal{L}_s(\chi_{3,2},\chi_{4,2} ; q) = \sum_{n= 1}^\infty  \frac{\chi_{3,2}(n)}{n^s} \Phi_0(\chi_{4,3};q^n)= \sum_{n=1}^\infty \frac{\chi_{3,2}(n)}{n^s} \frac{q^n}{1+q^{2n}}\,,\label{eq:Lchi3242chi}
\end{equation}
having used \eqref{eq:Phi043}.
Its asymptotic perturbative expansion as $y\to0^+$ can be computed directly from \eqref{eq:XiAsyPgen} and crucially it terminates after $s-1$ terms, i.e.
\begin{equation}
\mathcal{L}_s(\chi_{3,2},\chi_{4,3} ; y) \sim \mathcal{L}^{{\rm Pert}}_s(\chi_{3,2},\chi_{4,3} ; y) = \sum_{k=0}^{s-1} \frac{(-2\pi y)^k}{k!}  L(\chi_{3,2},s-k) L(\chi_{4,3},-k)\,.\label{eq:PertChi342chi}
\end{equation}

For the non-perturbative corrections we use \eqref{eq:TSgenRed} specialized to $\epsilon(\chi_{3,2})=\epsilon(\chi_{4,3})=1$, with $\kappa_1+\kappa_2=2$, and $(s_1,s_2)=(s,0)$, 
\begin{align}
\mathcal{L}^{\rm{NP}}_s(\chi_{3,2},\chi_{4,3} ; q) &=  -\frac{2}{\sqrt{3}}  (4 i y)^{s-1} \mathcal{L}_{s}\left({\chi}_{4,3},\chi_{3,2}; \tilde{q}\right) \notag \\
& = -\frac{2}{\sqrt{3}}  (4 i y)^{s-1}  \sum_{m=1}^\infty \frac{\chi_{4,3}(m)}{m^s} \frac{\tilde{q}^m}{1+\tilde{q}^m +\tilde{q}^{2m}}\,,
\end{align}
where $\tilde{q} = e^{-\frac{2\pi}{12y}}$.
We then deduce the Fricke inversion identity
\begin{align}
\sum_{m=1}^\infty \frac{\chi_{3,2}(m)}{m^s} \frac{q^m}{1+q^{2m}} =&\notag \sum_{k=0}^{s-1} \frac{(-2\pi y)^k}{k!}  L(\chi_{3,2},s-k) L(\chi_{4,3},-k)\\
&  -\frac{2}{\sqrt{3}}  (4 i y)^{s-1}  \sum_{m=1}^\infty \frac{\chi_{4,3}(m)}{m^s} \frac{\tilde{q}^m}{1+\tilde{q}^m +\tilde{q}^{2m}}\,.
\end{align}

From \eqref{eq:IterIntGen} we find for example the modular primitives
\begin{align}
\left(q \frac{d}{dq}\right) \mathcal{L}_1(\chi_{3,2},\chi_{4,3};q) &= G_2(\chi_{3,2},\chi_{4,3};q) =\eta(q)\eta(q^3) \left[\frac{\eta(q^2)\eta(q^{12})}{\eta(q^4)}\right]^2\,,\\
\left(q \frac{d}{dq}\right)^3 \mathcal{L}_3(\chi_{3,2},\chi_{4,3};q) &\notag = G_4(\chi_{3,2},\chi_{4,3};q) \\
&=\frac{\eta(q)\eta(q^2)^{14}\eta(q^3)\eta(q^{12})^4}{\eta(q^4)^8\eta(q^6)^4}
-12\frac{\eta(q^2)^4\eta(q^3)^7\eta(q^{12})^8}{\eta(q)\eta(q^4)^4\eta(q^6)^6}\,,
\end{align}
where the right hand side are manifestly modular forms of weight $2$ and $4$ respectively, for the congruence subgroup $\Gamma_0(12)$ and Dirichlet character $\chi_{12,11} =\chi_{3,2} \chi_{4,3}$.

\subsubsection*{Case $r_1=3,r_2=5$ and $s\in \mathbb{N}$ even}

We now consider~$\chi_{r_1}{=} \chi_{3,2}$ and~$\chi_{r_2} {=} \chi_{5,4}$. For the generalized Lambert series~$\mathcal{L}_s(\chi_{r_1},\chi_{r_2};q)$ to satisfy the condition~\eqref{eq:CondTrunc} we must require~$s\in \mathbb{N}$ even since the first character is odd while the other even.

The generalized Lambert series  \eqref{eq:L2chi}-\eqref{eq:GenLamPhi} here takes the form
\begin{equation}
\mathcal{L}_s(\chi_{3,2},\chi_{5,4} ; q) = \sum_{n= 1}^\infty  \frac{\chi_{3,2}(n)}{n^s} \Phi_0(\chi_{5,4};q^n)= \sum_{n=1}^\infty \frac{\chi_{3,2}(n)}{n^s} \frac{q^n(1+q^n)(1-q^{n})^2}{1-q^{5n}}\,,\label{eq:Lchi3254chi}
\end{equation}
having used \eqref{eq:Phi052}.
Its asymptotic perturbative expansion as $y\to0^+$ can be computed directly from \eqref{eq:XiAsyPgen} and crucially it terminates after $s-1$ terms, i.e.
\begin{equation}
\mathcal{L}_s(\chi_{3,2},\chi_{5,4} ; y) \sim \mathcal{L}^{{\rm Pert}}_s(\chi_{3,2},\chi_{5,4} ; y) = \sum_{k=1}^{s-1} \frac{(-2\pi y)^k}{k!}  L(\chi_{3,2},s-k) L(\chi_{5,4},-k)\,,\label{eq:PertChi352chi}
\end{equation}
where we notice that the $k=0$ term vanishes given that it corresponds to a trivial zero for $L(\chi_{5,4},-k)$.

For the non-perturbative corrections we use \eqref{eq:TSgenRed} specialized to $\epsilon(\chi_{3,2})=\epsilon(\chi_{5,4})=1$, with $\kappa_1+\kappa_2=1$, and $(s_1,s_2)=(s,0)$, 
\allowdisplaybreaks{
\begin{align}
\mathcal{L}^{\rm{NP}}_s(\chi_{3,2},\chi_{5,4} ; q) &=  i\sqrt{\frac{5}{3}}  (5i y)^{s-1} \mathcal{L}_{s}\left({\chi}_{5,4},\chi_{3,2}; \tilde{q}\right) \notag \\*
& = i\sqrt{\frac{5}{3}}  (5 i y)^{s-1} \sum_{m=1}^\infty \frac{\chi_{5,4}(m)}{m^s} \frac{\tilde{q}^m}{1+\tilde{q}^m +\tilde{q}^{2m}}\,,
\end{align}}
where $\tilde{q} = e^{-\frac{2\pi}{15y}}$.
We then deduce the Fricke inversion identity
\begin{align}
&\notag \sum_{m=1}^\infty \frac{\chi_{3,2}(m)}{m^s} \frac{q^m(1+q^m)(1-q^{m})^2}{1-q^{5m}} =\\
& \sum_{k=1}^{s-1} \frac{(-2\pi y)^k}{k!}  L(\chi_{3,2},s-k) L(\chi_{5,4},-k)
 +i (5 i y)^{s-1} \sqrt{\frac{5}{3}} \sum_{m=1}^\infty \frac{\chi_{5,4}(m)}{m^s} \frac{\tilde{q}^m}{1+\tilde{q}^m +\tilde{q}^{2m}}\,.
\end{align}

From \eqref{eq:IterIntGen} we find for example the modular primitives
\begin{align}
\mathcal{L}_0(\chi_{3,2},\chi_{5,4};q) &= G_1(\chi_{3,2},\chi_{5,4};q) = \frac{[\eta(q)\eta(q^{15})]^2}{\eta(q^3)\eta(q^5)}\,,\\
 \left( q\frac{d}{dq}\right)^2 \mathcal{L}_2(\chi_{3,2},\chi_{5,4};q)  &\notag= G_3(\chi_{3,2},\chi_{5,4};q) \\
& =\left[\eta(q^5)\eta(q^3)\right]^3{-}5\left[\eta(q)\eta(q^{15})\right]^3
{-}24\eta(q^{15})^5\eta(q)\left[\frac{\eta(q^3)}{\eta(q^5)}\right]^2\,,
\end{align}
where the right hand side are modular forms of weight $1$ and $3$ respectively, for the congruence subgroup $\Gamma_0(15)$ and Dirichlet character $\chi_{15,14} =\chi_{3,2} \chi_{5,4}$.
\vspace{0.2cm}

As already noted when discussing singly twisted Eisenstein series, here as well we see from the many examples provided that for some values of the modular weight $m=s+1$ and some choices of primitive characters $\chi_{r_1},\chi_{r_2}$, the corresponding twisted Eisenstein series $G_m(\chi_{r_1},\chi_{r_2}; q)$ can be expressed as linear combinations with integer coefficients of eta quotients. Following the same convention as in Table \ref{tab:SingleTwist}, we present in Table \ref{tab:DoubleTwist} all of the instances where for a given level $N\in \mathbb{N}$, i.e. with respect to the congruence subgroup $\Gamma_0(N)$, modular weight $m$, and characters $\chi_{r_1}(a) = \left(\frac{D_1}{a}\right)$ and $\chi_{r_2}(a) = \left(\frac{D_2}{a}\right)$, with $D_1,D_2 \in \mathbb{Z}$ two fundamental discriminants, the corresponding twisted Eisenstein series $G_m(\chi_{r_1},\chi_{r_2};q)$  can be written as a single eta quotient. In Table \ref{tab:DoubleTwist} we use again the short-hand notation $\eta_t^r \coloneqq \eta(q^t)^r$.

Lastly, one may wish to consider cases where the generalised Lambert series $\mathcal{L}_s(\chi_{r_1},\chi_{r_2} ; q)$ has two imprimitive characters. However, given our discussion in Section \ref{sec:ExImp}, it should appear clear that such instances will produce linear combinations of previously discussed instances.
For example, similarly to the single imprimitive character case $\mathcal{L}_s(\chi_{6,5} ; q)$ studied before we can easily see that
\begin{equation}
\mathcal{L}_s(\chi_{6,5},\chi_{6,5} ; q) = \mathcal{L}_s(\chi_{3,2},\chi_{3,2} ; q) +(1+2^{-s}) \mathcal{L}_s(\chi_{3,2},\chi_{3,2} ; q^2) + 2^{-s} \mathcal{L}_s(\chi_{3,2},\chi_{3,2} ; q^4)\,,
\end{equation}
thus reducing the problem to the simpler double primitive character case $ \mathcal{L}_s(\chi_{3,2},\chi_{3,2} ; q)$ already discussed.
For this reason we shall not present such examples here, although we stress that they can all be studied using our general analysis presented in Section \ref{sec:NPImp}.

\begin{table}[h]\centering
\begin{tabular}{|r|r|r|r|c|c|c|}\hline
$N$ & $m$ & $D_1$ & $D_2$ & Numerator & Denominator & OEIS Sequence\\ \hline 
9 & 2 & -3 & -3 & $\eta_{1}^{3} \eta_{9}^{3} $ & $\eta_{3}^{2}$ & \oeis{A106401} \\
12 & 2 & -3 & -4 & $\eta_{1}^{\phantom{1}} \eta_{2}^{2} \eta_{3}^{\phantom{1}} \eta_{12}^{2} $ & $\eta_{4}^{2} $ & \oeis{A209613} \\
12 & 2 & -4 & -3 & $\eta_{1}^{2} \eta_{4}^{\phantom{1}} \eta_{6}^{2} \eta_{12}^{\phantom{1}} $ & $\eta_{3}^{2} $ & \oeis{A113421} \\
15 & 1 & -3 & 5 & $\eta_{1}^{2} \eta_{15}^{2} $ & $\eta_{3}^{\phantom{1}} \eta_{5}^{\phantom{1}} $ & \oeis{A106406} \\
16 & 2 & -4 & -4 & $\eta_{2}^{4} \eta_{8}^{4} $ & $\eta_{4}^{4} $ & \oeis{A121613} \\
20 & 1 & -4 & 5 & $\eta_{1}^{\phantom{1}} \eta_{2}^{\phantom{1}} \eta_{10}^{\phantom{1}} \eta_{20}^{\phantom{1}} $ & $\eta_{4}^{\phantom{1}} \eta_{5}^{\phantom{1}} $ & \oeis{A111949} \\
24 & 1 & -3 & 8 & $\eta_{1}^{\phantom{1}} \eta_{4}^{\phantom{1}} \eta_{6}^{\phantom{1}} \eta_{24}^{\phantom{1}} $ & $\eta_{3}^{\phantom{1}} \eta_{8}^{\phantom{1}} $ & \oeis{A115660} \\
32 & 1 & -4 & 8 & $\eta_{2}^{2} \eta_{16}^{2} $ & $\eta_{4}^{\phantom{1}} \eta_{8}^{\phantom{1}} $ & \oeis{A125095} \\
32 & 2 & -4 & -8 & $\eta_{4}^{9} \eta_{16}^{2} $ & $\eta_{2}^{2} \eta_{8}^{5} $ & \oeis{A113419} \\
32 & 2 & -8 & -4 & $\eta_{2}^{2} \eta_{8}^{9} $ & $\eta_{4}^{5} \eta_{16}^{2} $ & \oeis{A209940} \\
36 & 1 & -3 & 12 & $\eta_{1}^{\phantom{1}} \eta_{4}^{\phantom{1}} \eta_{6}^{4} \eta_{9}^{\phantom{1}} \eta_{36}^{\phantom{1}} $ & $\eta_{2}^{\phantom{1}} \eta_{3}^{2} \eta_{12}^{2} \eta_{18}^{\phantom{1}} $ & \oeis{A129448} \\
48 & 1 & -4 & 12 & $\eta_{2}^{\phantom{1}} \eta_{6}^{\phantom{1}} \eta_{8}^{\phantom{1}} \eta_{24}^{\phantom{1}} $ & $\eta_{4}^{\phantom{1}} \eta_{12}^{\phantom{1}} $ & \oeis{A129449} \\
64 & 1 & -8 & 8 & $\eta_{4}^{2} \eta_{16}^{2} $ & $\eta_{8}^{2} $ & \oeis{A134343} \\
96 & 1 & -4 & 24 & $\eta_{2}^{\phantom{1}} \eta_{8}^{4} \eta_{12}^{4} \eta_{48}^{\phantom{1}} $ & $\eta_{4}^{3} \eta_{6}^{\phantom{1}} \eta_{16}^{\phantom{1}} \eta_{24}^{3} $ & \oeis{A190615} \\
96 & 1 & -8 & 12 & $\eta_{4}^{4} \eta_{6}^{\phantom{1}} \eta_{16}^{\phantom{1}} \eta_{24}^{4} $ & $\eta_{2}^{\phantom{1}} \eta_{8}^{3} \eta_{12}^{3} \eta_{48}^{\phantom{1}} $ & \oeis{A134177} \\
576 & 1 & -24 & 24 & $\eta_{8}^{5} \eta_{12}^{2} \eta_{48}^{2} \eta_{72}^{5} $ & $\eta_{4}^{2} \eta_{16}^{2} \eta_{24}^{4} \eta_{36}^{2} \eta_{144}^{2} $ & \oeis{A208955} \\
\hline\end{tabular}
\caption{
Numerators and denominators of eta quotients that are doubly twisted Eisenstein series,
in the same format as for Table \ref{tab:SingleTwist}.}
\label{tab:DoubleTwist}
\end{table}

\section{An application to quantum dilogarithms}
\label{sec:Qpoch}

In this section we want to discuss a particular example of $q$-series \eqref{eq:LamChi} and \eqref{eq:Lam2} with parameter $s=m-1\in \mathbb{N}$ for which the integer $m$ has opposite parity to the Dirichlet character.
The examples we are about to present do not originate from iterated integrals of twisted Eisenstein series as in section \ref{sec:ModularPrim}. In particular, since the present choice of parameters does not satisfy the condition \eqref{eq:CondTrunc} we will see that the relevant asymptotic perturbative expansions do produce formal, factorially divergent power series whose non-perturbative completions rely on the quantum Fricke relations \eqref{eq:TSQM}. 

\subsection{$q$-Pochhammer and the fermionic spectral trace of local $\mathbb{P}^2$}

The present discussion has been inspired by the recent work \cite{Fantini:2024snx} whose focus has been understanding the resurgent structures of a particular observable, called fermionic spectral trace, in topological string theory compactified on the Calabi-Yau manifold known as local $\mathbb{P}^2$ geometry, the simplest example of a toric del Pezzo Calabi-Yau threefold.
The spectral traces associated with these toric Calabi-Yau manifold were first introduced and computed as observables in the context of the topological string/spectral theory correspondence in \cite{Grassi:2014zfa,Codesido:2015dia}.
The resurgent structures of the first fermionic spectral trace of local $\mathbb{P}^2$ has been solved in \cite{Rella:2022bwn} both in the weak and strong coupling regime. 
Surprisingly \cite{Rella:2022bwn} found an intriguing `number-theoretic duality' relating the perturbative and non-perturbative data in the two regimes, connecting the perturbative series and the Stokes discontinuity in the weak coupling regime with those of the strong coupling regime. Building on this work, \cite{Fantini:2024snx} clarified and expanded upon this notion of strong-weak resurgence symmetry finding deep connections between this structure and the quantum modularity properties of this observable, in particular under the action of Fricke involution, eventually leading to the notion of  \textit{modular resurgent paradigm} in \cite{Fantini:2024ihf}.  

As discussed in \cite{Rella:2022bwn,Fantini:2024snx}, the first fermionic spectral trace of local $\mathbb{P}^2$ takes the form
\begin{equation}
\mbox{Tr}(\rho_{\mathbb{P}^2}) = \frac{1}{\sqrt{3} \mb} e^{-\frac{\pi i }{36} \mb^2 +\frac{\pi i}{12} \mb^{-2} +\frac{\pi i }{4}}\, \frac{(q'^{\frac{2}{3}};q')_\infty^2}{(q'^{\frac{1}{3}};q')_\infty}\cdot  \frac{(\omega_3;\tilde{q})_\infty}{(\omega_3^{-1};\tilde{q})_\infty^2}\,.\label{eq:TrFR}
\end{equation}
where we have defined
\begin{equation}
q' = e^{2\pi i \mb^2}\,, \qquad \tilde{q}= e^{-\frac{2\pi i }{\mb^2}}\,,\qquad \omega_3 = e^{\frac{2\pi i }{3}}\,.
\end{equation}
The parameter $\mb^2$ lies in the upper-half plane, i.e. ${\rm Im}(\mb^2)>0$, and it is related to the quantum deformation parameter, $\hbar$, of the spectral theory by $\hbar = 2\pi \mb^2 /3$ which in turn is related to the inverse string coupling of the topological string via $g_s = 4\pi^2 / \hbar$.

Importantly, in \eqref{eq:TrFR} the function $(x q^\alpha;q)_\infty$, referred to as quantum dilogarithm or (infinite) $q$-Pochhammer symbol, denotes the analytic function of the two complex variables $x,q\in \mathbb{C}$ with $|q|<1$ defined by the infinite product
\begin{equation}
(x q^\alpha;q)_\infty \coloneqq \prod_{n=0}^\infty \left(1-xq^{\alpha+n}\right)\,,\label{eq:qLog}
\end{equation}
where the parameter $\alpha\in \mathbb{R}$. 

To make a direct connection with the formulae here presented we rewrite \eqref{eq:TrFR} in terms of the upper-half plane variable $\tau = \mb^2 /3= \hbar/(2\pi)$,
\begin{equation}
\mbox{Tr}(\rho_{\mathbb{P}^2}) = \frac{1}{3 \sqrt{\tau}} e^{-\frac{\pi i }{12} \tau +\frac{\pi i}{ 36} \tau^{-1} +\frac{\pi i }{4}}\, \frac{(q^{2};q^3)_\infty^2}{(q;q^3)_\infty}\cdot  \frac{(\omega_3;\tilde{q})_\infty}{(\omega_3^{-1};\tilde{q})_\infty^2}\,,\label{eq:TrFR2}
\end{equation}
where now $q=e^{2\pi i \tau}$ and $\tilde{q} = e^{-\frac{2\pi i }{3\tau}}$ with $|q|<1$ and $|\tilde{q}|<1$.
We now proceed to express the two ratios of quantum-dilogarithms in terms of the $q$-series \eqref{eq:qser1}-\eqref{eq:qser2} and \eqref{eq:qserDD}.
From the definition \eqref{eq:qLog} it is easy to check that
\begin{align}
&\log \Big[ \frac{(q^{2};q^3)_\infty^2}{(q;q^3)_\infty} \Big]  \notag\\
&\notag = \sum_{n_2=1}^\infty  2 \log(1-q^{3n_2 -1}) - \log(1 - q^{3n_2 - 2}) = \sum_{n_1,n_2=1}^\infty  \frac{1}{n_1}\left(-2 q^{n_1(3n_2-1)} + q^{n_1(3n_2-2)}\right) \\
&= \sum_{N=1}^\infty \frac{3 \sigma_{1,\chi_{3,2}}(N) - \sigma_{1,\chi_{3,1}}(N)}{2N} q^N   = \frac{3}{2} \tilde{\mathcal{L}}_1(\chi_{3,2};q) -\frac{1}{2} \left(\mathcal{L}_1(\chi_{1,1};q) -\mathcal{L}_1(\chi_{1,1};q^3)\right)\,,
\label{eq:qpoch1}
\end{align}
where in the last equation we made use of the identity \eqref{eq:chi11Lam}.

Similarly, for the second ratio of quantum-dilogarithms we find,
\begin{align}
&\notag \log \Big[ \frac{(\omega_3;\tilde{q})_\infty }{(\omega_3^{-1};\tilde{q})^2_\infty} \Big] \\
&\notag = \sum_{n_2=0}^\infty  \log(1-\omega_3 \tilde{q}^{n_2}) - 2\log(1 - \omega_3^{-1}\tilde{q}^{n_2}) = -\frac{\log(3)+ i \pi}{2}+ \sum_{n_1,n_2=1}^\infty \frac{(2\omega_3^{-n_1}-\omega_3^{n_1})}{n_1} \tilde{q}^{n_1 n_2}\\
&\notag =   -\frac{\log(3)+ i \pi}{2} - \sum_{N=1}^\infty \frac{ 3 \sqrt{3}i \sigma'_{1,\chi_{3,2}}(N)+ \sigma_{1,\chi_{3,1}}(N)}{N} \tilde{q}^N\\
&= -\frac{\log(3)+ i \pi}{2} - \frac{3\sqrt{3}}{2} i \,\mathcal{L}_1(\chi_{3,2};\tilde{q}) -\frac{1}{2} \left(\mathcal{L}_1(\chi_{1,1};\tilde{q}) -\mathcal{L}_1(\chi_{1,1};\tilde{q}^3)\right)\,,\label{eq:qpoch2}
\end{align}
where we use the trigonometric identity
\begin{equation}
2\omega_3^{-n} -\omega_3^n = 1 - \frac{3}{2} \chi_{3,1}(n) - \frac{3 \sqrt{3}  }{2}\cdot i\, \chi_{3,2}(n) \,,\qquad {\rm with} \quad \omega_3=e^{\frac{2\pi i}{3}}\,,
\end{equation}
valid for all $n\in \mathbb{N}$, combined with \eqref{eq:DivId} and \eqref{eq:chi11Lam}.

Substituting \eqref{eq:qpoch1}-\eqref{eq:qpoch2} in \eqref{eq:TrFR2} and expressing everything in terms of the variable $\tau$ using $q=e^{2\pi i \tau}$ and $\tilde{q} = e^{-\frac{2\pi i }{3\tau}}$, we conclude that the first fermionic spectral trace of local $\mathbb{P}^2$ can be expressed as
\allowdisplaybreaks{
\begin{align}
\log\Big[\mbox{Tr}(\rho_{\mathbb{P}^2})\Big] =&\notag - \log(3 \sqrt{3\tau}) -\frac{\pi i }{12} \tau +\frac{\pi i}{ 36} \tau^{-1} -\frac{\pi i }{4}\\*
&\notag -\frac{1}{2} \left(\mathcal{L}_1(\chi_{1,1};\tau)-\mathcal{L}_1(\chi_{1,1};-\frac{1}{\tau})\right) +\frac{1}{2}\left(\mathcal{L}_1(\chi_{1,1};3\tau)-\mathcal{L}_1(\chi_{1,1};-\frac{1}{3\tau})\right)\\*
&+ \frac{3}{2} \left(\tilde{\mathcal{L}}_1(\chi_{3,2};\tau) - i \sqrt{3} \mathcal{L}_1(\chi_{3,2};-\frac{1}{3\tau}) \right)\,.\label{eq:LogTr}
\end{align}}
Furthermore, thanks to the resurgent analysis of \cite{Dorigoni:2020oon} we know that $\mathcal{L}_1(\chi_{1,1};\tau)$, which is the integral of the quasi-modular form $G_2(\tau)$ or alternatively $\mathcal{L}_1(\chi_{1,1};q) = -\log[(q;q)_\infty]$, satisfies the $S$-transformation identity
\begin{equation}
\mathcal{L}_1(\chi_{1,1};\tau)  = \frac{\log (- i\tau)}{2} + \frac{\pi}{12} i ( \tau + \tau^{-1}) + \mathcal{L}_1(\chi_{1,1};-\frac{1}{\tau})\,, \label{eq:L1S}
\end{equation}
so that \eqref{eq:LogTr} further simplifies dramatically to
\begin{align}
\log\Big[\mbox{Tr}(\rho_{\mathbb{P}^2})\Big] =& - \frac{1}{2}\log(3^{\frac{5}{2}}\tau) -\frac{\pi i}{4} + \frac{3}{2} \left(\tilde{\mathcal{L}}_1(\chi_{3,2};\tau) - i \sqrt{3} \mathcal{L}_1(\chi_{3,2};-\frac{1}{3\tau}) \right)\,.\label{eq:LogTr2}
\end{align}
We see that the first fermionic spectral trace is precisely related to the $q$-series \eqref{eq:qser1}-\eqref{eq:qser2} for the case where $s=1$ and odd character $\chi_{3,2}$.
Note that in this present case the parity of $s$ is the same as that of the character and the $q$-series $\mathcal{L}_1(\chi_{3,2};q)$ and $\tilde{\mathcal{L}}_1(\chi_{3,2};q) $ do not originate from iterated integrals of twisted Eisenstein series.
Importantly, in the appropriate regimes the corresponding perturbative expansions give rise to the formal, asymptotic power series given in \eqref{eq:LtP} and \eqref{eq:AsyL}. 

The resurgent structures found in \cite{Rella:2022bwn,Fantini:2024snx} can be interpreted in terms of the general quantum-modular Fricke inversion \eqref{eq:TSQM}.
Specialising \eqref{eq:TSQM} to the case $s=1$ and $\chi_r = \chi_{3,2}$ we find:
\begin{equation}
\left( \begin{matrix} 
\mathcal{L}_{1}(\chi_{3,2};\tau) \\ 
{\tilde{\mathcal{L}}}_{1}(\chi_{3,2};\tau ) 
\end{matrix} \right) = \left(\begin{matrix} \mathcal{S}_-\left[ \mathcal{L}^{{\rm Pert}}_{1} \right](\chi_{3,2};\tau ) \vspace{0.1cm}\\ \mathcal{S}_-\left[ \tilde{\mathcal{L}}^{{\rm Pert}}_{1} \right](\chi_{3,2};\tau)
\end{matrix}\right) +  i  \left(\begin{matrix} 0 & \frac{1}{\sqrt{3}} \\ \sqrt{3} & 0\end{matrix}\right)\left( \begin{matrix} 
\mathcal{L}_{1}({\chi}_{3,2};-\frac{1}{3\tau}) \\ 
{\tilde{\mathcal{L}}}_{1}({\chi}_{3,2};-\frac{1}{3\tau} )\label{eq:TSQM32}
\end{matrix} \right) \,.
\end{equation}
The asymptotic expansion of \eqref{eq:LogTr2} for $\tau\to0$, i.e. $\hbar\to0$ in the variable of \cite{Rella:2022bwn,Fantini:2024snx}, can be obtained by substituting the second entry of the vector identity \eqref{eq:TSQM32} directly into \eqref{eq:LogTr2}. Similarly, for the asymptotic expansion of \eqref{eq:LogTr2} for $\tau\to  \infty$, i.e. $\hbar\to \infty$ in the variable of \cite{Rella:2022bwn,Fantini:2024snx}, we must use the first entry of \eqref{eq:TSQM32} upon substituting $\tau\to -1/(3\tau)$. We find
\begin{align}
\log\Big[\mbox{Tr}(\rho_{\mathbb{P}^2})\Big]  &\label{eq:tau0}\stackrel{\tau\to0}{=} - \frac{1}{2}\log(3^{\frac{5}{2}}\tau) -\frac{\pi i}{4} + \frac{3}{2} \mathcal{S}_-\left[ \tilde{\mathcal{L}}^{{\rm Pert}}_{1} \right](\chi_{3,2};\tau)\,, \\
\log\Big[\mbox{Tr}(\rho_{\mathbb{P}^2})\Big]  &\label{eq:tauInfty}\stackrel{\tau\to \infty}{=} - \frac{1}{2}\log(3^{\frac{5}{2}}\tau)-\frac{\pi i}{4} - \frac{3\sqrt{3}}{2} i\,  \mathcal{S}_-\left[ {\mathcal{L}}^{{\rm Pert}}_{1} \right](\chi_{3,2};-\frac{1}{3\tau})\,.
\end{align}

We then compute the formal, asymptotic power series \eqref{eq:LtP} and \eqref{eq:AsyL} using the explicit form for the Dirichlet $L$-values \eqref{eq:NegInt} to arrive at
\allowdisplaybreaks{
\begin{align}
&\label{eq:tau0PT}\frac{3}{2} \tilde{\mathcal{L}}^{{\rm Pert}}_{1} (\chi_{3,2};\tau)  = -\frac{1}{2} \log\left(\frac{16\pi^4 \tau}{\sqrt{3}}\right)+ 3\log \Gamma\left(\tfrac{1}{3}\right) + \frac{\pi i }{4} - 3 \sum_{n=1}^\infty \frac{B_{2n} B_{2n+1}(\frac{2}{3})}{2n(2n+1)!} (2\pi i \tau)^{2n} \,,\\*
&\label{eq:tauInftyPT}- \frac{3\sqrt{3}}{2} i\mathcal{L}^{{\rm Pert}}_{1} (\chi_{3,2};-\frac{1}{3\tau})  = -\left[\frac{\psi^{(1)}(\tfrac{1}{3}) - \psi^{(1)}(\tfrac{2}{3}) }{2\sqrt{3} }\right] \frac{3 \tau}{2\pi} + \frac{\pi i }{4} + i \sqrt{3} \sum_{n=1}^\infty \frac{B_{2n} B_{2n-1}(\frac{2}{3})}{ (2n-1) (2n)!} \left(\frac{2\pi }{i \tau}\right)^{2n-1} \,,
\end{align}}
where again $\psi^{(n)}(x)$ denotes the $n^{th}$ derivative of the digamma function $\psi(x) = \Gamma'(x)/\Gamma(x)$ and $B_k(x)$ denotes the $k^{th}$ Bernoulli polynomial.
Note that following the analysis carried out in section \ref{sec:Asy}, we find that the asymptotic expansion for $\tilde{\mathcal{L}}^{{\rm Pert}}_{1} (\chi_{3,2};\tau)$ near $\tau\to0$ does contain a logarithmic term as given in the general formula \eqref{eq:Asym2}.\footnote{It is worth noting that the logarithmic term $\log(\tau) = \log(\frac{\hbar}{2\pi})$ can be predicted from the Topological Strings/Spectral Theory (TS/ST) correspondence and it is in fact a generic feature of the perturbative expansion as $\hbar \to 0$ of the spectral traces of toric Calabi-Yau threefolds, see \cite{Rella:2022bwn}. We thank Claudia Rella for related discussions.}

Substituting \eqref{eq:tau0PT} and~\eqref{eq:tauInftyPT} respectively in~\eqref{eq:tau0} and~\eqref{eq:tauInfty} and changing variables back to~$\tau = \hbar/(2\pi)$ we retrieve exactly the asymptotic expansions derived in~\cite{Rella:2022bwn} in the limit~$\hbar\to0$ and~$\hbar\to \infty$.\footnote{For a precise match we note that quantity $V=[\psi^{(1)}\left(\frac{1}{3}\right) - \psi^{(1)}\left(\frac{2}{3}\right)]/(2\sqrt{3})$, here derived from the Riemann term in the asymptotic expansion \eqref{eq:AsyL} specialized to $s=1$ and character $\chi_{3,2}$, is also equal to $V = 2\,{\rm Im}( {\rm Li}_2(e^{\frac{ i \pi}{3}}))$ as defined in the reference.}.

\subsection{The spectral trace of local $\mathbb{P}^{m,n}$} 

Recently, \cite{Fantini:2025wap} analyzed a generalization of the $\mathbb{P}^2$ spectral trace \eqref{eq:TrFR} and considered the spectral trace ${\rm Tr}(\rho_{m,n})$  canonically associated with a local $\mathbb{P}^{m,n}$ geometry with $m,n\in \mathbb{N}$. 
This observable has been computed in \cite{Kashaev:2015kha} and it is once again given in terms of q-Pochhammer symbols \eqref{eq:qLog}, 
\begin{equation}
{\rm Tr}(\rho_{m,n}) =  \frac{1}{2 \sqrt{N \tau} \sin\left(\frac{\pi n }{N}\right)} e^{-\frac{\pi i }{12} \tau\left(\frac{m}{N}-\frac{N}{12}\right) +\frac{\pi i}{ 12 N} \tau^{-1} +\frac{\pi i }{4}} \exp\left[ G_{m,n}(q) +  F_{m,n}(\tilde{q}) \right]\,,\label{eq:rhonm}
\end{equation}
where $q=e^{2\pi i \tau}$ and $\tilde{q} = e^{-\frac{2\pi i }{ N\tau}}$ with $N= m+n+1$, and where we defined the functions:
\begin{align}
G_{m,n}(q)  &\coloneqq \log\left[\frac{(q^{N-m};q^N)_\infty (q^{N-1};q^N)_\infty}{(q^n;q^N)_\infty} \right]\,,\label{eq:G} \\
F_{m,n}(q) & \coloneqq \log\left[\frac{(\omega_N^n;q)_\infty }{(\omega_N^{N-m};q)_\infty (\omega_N^{N-1};q)_\infty } \right]\label{eq:F}\,,
\end{align}
with $\omega_N \coloneqq e^{\frac{2\pi i}{N}}$ denoting the $N^{th}$ root of unity.
Note that \eqref{eq:rhonm} is identical to \eqref{eq:TrFR2} for $n=m=1$, as expected given that the local $\mathbb{P}^{1,1}$ geometry coincides with $\mathbb{P}^2$.

In \cite{Fantini:2025wap}, the calculation of the all-order asymptotic expansion of \eqref{eq:rhonm} for $\tau\to \infty$ and $\tau\to0$ is aided by the auxiliary functions
\begin{align}
\mathfrak{g}(q) & \coloneqq \sum_{k=1}^N \chi_{N}(k) \log(q^k ; q^N)_\infty \label{eq:g}\,,\\
\mathfrak{f}(q) &\coloneqq \sum_{k =1}^N \chi_{N}(k) \log(\omega_N^k;q)_\infty\,,\label{eq:f}
\end{align}
where the authors considered the case where $\chi_N$ is a primitve odd character modulo $N$.
However, we can consider the case of $\chi_N$ a general primitive, not necessarily real, character modulo $N$ for which it is straightforward to see that
\begin{align}
\mathfrak{g}(q) & =  - \tilde{\mathcal{L}}_1(\chi_N;q) \label{eq:gq} \,,\\
\mathfrak{f}(q) &=  i^\kappa \epsilon(\chi_{N}) \sqrt{N} \mathcal{L}_1(\bar{\chi}_N;q) \label{eq:fq} \,,
\end{align}
where as always $\kappa \in \{0,1\}$ denotes the parity of the character, i.e. $\chi_N(-1) = (-1)^\kappa$ and $\epsilon(\chi_N)$ has been defined in \eqref{eq:GaussIT}.
We conclude that the resurgence properties of \eqref{eq:g} and \eqref{eq:f} are then captured by the Fricke involution \eqref{eq:TSQM}, where in particular for the case where $\chi_N$ is an even character these reduce to the Cheshire resurgence structure connected with the iterated twisted Eisenstein integral nature of the $q$-series considered.

We note that to discuss the resurgent properties of the spectral trace $\rm{Tr}(\rho_{m,n})$ and in particular of its building blocks \eqref{eq:G}-\eqref{eq:F}, one does in general encounter linear combinations of both resurgent and Cheshire-resurgent Lambert series with a character.
For the case $n=m=1$ corresponding to $N=3$ we have a unique odd primitive character $\chi_{3,2}$ and the resurgent structures can be read immediately from \eqref{eq:LogTr2}.
Similarly, for the case $N=4$ we have once again a unique non-principal character given by the odd primitive character $\chi_{4,3}$. We proceed as we did above for the case $N=3$, and consider first the case $m=2$ and $n=1$ and it is easy to show by expanding the logarithms that the building blocks \eqref{eq:G}-\eqref{eq:F} in the language of the present work are given by
\begin{align}
G_{2,1}(q) &= \tilde{\mathcal{L}}_1(\chi_{4,3};q) - \mathcal{L}_1(\chi_{1,1};q^2)+\mathcal{L}_1(\chi_{1,1};q^4)\,,\\
F_{2,1}(q) &= -\left(\log(2) + \frac{i \pi}{2}\right)-2i {\mathcal{L}}_1(\chi_{4,3};q) - \mathcal{L}_1(\chi_{1,1};q)+\mathcal{L}_1(\chi_{1,1};q^2)\,,
\end{align}
Using the modular property \eqref{eq:L1S} for $\mathcal{L}(\chi_{1,1};q)$ we then see that the combination appearing in \eqref{eq:rhonm} takes the neat form
\begin{equation}
G_{2,1}(q) + F_{2,1}(\tilde{q}) =\tilde{\mathcal{L}}_1(\chi_{4,3};q)-2i {\mathcal{L}}_1(\chi_{4,3};\tilde{q})-\frac{i \pi }{48 \tau }+\frac{i \pi  \tau }{6}- \frac{i \pi }{2}-\frac{\log (2)}{2}\,,\label{eq:21}
\end{equation}
with $\tilde{q} = e^{\frac{2\pi i }{4\tau}}$.
The weak/strong coupling resurgence properties are then fully captured by \eqref{eq:TSQM} specialized to the case $s=1$ and $\chi_r = \chi_{4,3}$.

Similarly, for the case $m=1$ and $n=2$ we find that the building blocks \eqref{eq:G}-\eqref{eq:F} can be written as
\allowdisplaybreaks{
\begin{align}
G_{1,2}(q) &= \tilde{\mathcal{L}}_1(\chi_{4,3};q) - \mathcal{L}_1(\chi_{1,1};q)+2\mathcal{L}_1(\chi_{1,1};q^2)-\mathcal{L}_1(\chi_{1,1};q^4)\,,\\
F_{1,2}(q) &= -  \frac{ i\pi}{2}-2i {\mathcal{L}}_1(\chi_{4,3};q) + \mathcal{L}_1(\chi_{1,1};q)-2\mathcal{L}_1(\chi_{1,1};q^2)+\mathcal{L}_1(\chi_{1,1};q^4)\,.
\end{align}}
Using again \eqref{eq:L1S}, we simplify the combination appearing in \eqref{eq:rhonm} to
\begin{equation}
G_{1,2}(q) + F_{1,2}(\tilde{q}) =\tilde{\mathcal{L}}_1(\chi_{4,3};q)-2i {\mathcal{L}}_1(\chi_{4,3};\tilde{q})-\frac{i \pi }{48 \tau }-\frac{i \pi  \tau }{12}-  \frac{i \pi}{2} \,,\label{eq:12}
\end{equation}
for which again the resurgence properties follow from \eqref{eq:TSQM}.
Substituting \eqref{eq:21} and \eqref{eq:12} in \eqref{eq:rhonm} for the respective values of $n,m$ we immediately see that ${\rm Tr}(\rho_{2,1}) = {\rm Tr}(\rho_{1,2})$ as it should.

As already stressed the cases $N=3$, given in  \eqref{eq:LogTr2}, and $N=4$ in \eqref{eq:21}-\eqref{eq:12}, are special in the sense that they each posses a unique non-principal character modulo $N$, namely the odd character $\chi_N(n) = \left( \frac{-N}{n}\right)$. In this simple scenario the weak/strong coupling resurgent properties of $\Tr(\rho_{m,n})$ are immediately encoded in the single building blocks \eqref{eq:gq}-\eqref{eq:fq} constructed out of the unique odd character $\chi_N$. However, for higher values of $N$ we have that \eqref{eq:G} and \eqref{eq:F} are in general given by linear combinations of the building blocks \eqref{eq:gq}-\eqref{eq:fq}, involving both odd \textit{and even} characters, thus displaying both resurgence and Cheshire-resurgence properties.

This can be exemplified neatly by considering the $N=5$ case, for which we have the primitive real character $\chi_{5,4}(n) = \left( \frac{5}{n}\right)$ as well as the two primitive and complex conjugated odd character $\chi_{5,2}(n) = \overline{\chi_{5,3}(n)}$ given below.
\begin{center}
\begin{tabular}{ |c|c|c|c|c|} 
\hline
 $n\,\rm{mod}\,5$ & 1 & 2 & 3 & 4 \\
\hline
\multirow{1}{4em}{$\chi_{5,2}(n)$} &1 & i &-i& -1\\
\hline
\end{tabular}
\end{center}

For simplicity, let us analyze the case $m=3$ and $n=1$, although similar results can be presented for different values of $n,m$ such that $m+n+1=5$. As with our previous analysis, we expand the logarithm in \eqref{eq:G} and \eqref{eq:F} to obtain
\allowdisplaybreaks{
\begin{align}
G_{3,1}(q) &\notag = \left(\frac{2+i}{4}\right) \tilde{\mathcal{L}}_1(\chi_{5,2},q)+\left(\frac{2-i}{4}\right) \tilde{\mathcal{L}}_1(\overline{\chi_{5,2}},q)+\frac{1}{4} \tilde{\mathcal{L}}_1(\chi_{5,4},q)\\*
&\phantom{=} -\frac{1}{4}\left( \mathcal{L}_1(\chi_{1,1};q) - \mathcal{L}_1(\chi_{1,1};q^5 \right))\,,\\
F_{3,1}(q) & \notag=  \frac{1}{2} \log \left( \frac{5-\sqrt{5}}{10}\right)- \frac{i \pi}{2} - i \epsilon(\chi_{5,2}) \sqrt{5} \left(\frac{2+i}{4}\right)  {\mathcal{L}}_1(\overline{\chi_{5,2}},q)\\*
&\phantom{=}-i \overline{\epsilon(\chi_{5,2})} \sqrt{5} \left(\frac{2-i}{4}\right) {\mathcal{L}}_1({\chi_{5,2}},q)- \frac{ \sqrt{5}}{4} {\mathcal{L}}_1(\chi_{5,4},q)-\frac{1}{4}\left( \mathcal{L}_1(\chi_{1,1};q) - \mathcal{L}_1(\chi_{1,1};q^5 \right))\,,
\end{align}}
where the Gauss sum defined in \eqref{eq:GaussIT} equals $\epsilon(\chi_{5,2})=e^{\frac{i}{4} \left[\pi - \arctan\left(\frac{4}{3}\right)\right]}$ and we remind the reader that since $\chi_{5,4}$ is a real primitive character we have $\epsilon(\chi_{5,4})=1$.
Just like before, we use the modular property \eqref{eq:L1S} to simplify the combination appearing in \eqref{eq:rhonm} and bring it to the neat form
\begin{align}
&G_{3,1}(q) + F_{3,1}(\tilde{q}) \notag = \\
&\notag \left(\frac{2+i}{4}\right) \Big[ \tilde{\mathcal{L}}_1(\chi_{5,2},q) - i \epsilon(\chi_{5,2}) \sqrt{5} {\mathcal{L}}_1(\overline{\chi_{5,2}},\tilde{q}) \Big]+\left(\frac{2-i}{4}\right) \Big[ \tilde{\mathcal{L}}_1(\overline{\chi_{5,2}},q) -i \overline{\epsilon(\chi_{5,2})} \sqrt{5}{\mathcal{L}}_1({\chi_{5,2}},\tilde{q}) \Big]\\
&\label{eq:GF31} -\frac{1}{4} \Big[ \tilde{\mathcal{L}}_1(\chi_{5,4},q) - \sqrt{5}{\mathcal{L}}_1(\chi_{5,4},\tilde{q})\Big] -\frac{i \pi }{60 \tau }+\frac{i \pi  \tau }{12}-\frac{i \pi}{2} -\frac{1}{8} \log \left( \frac{7-3 \sqrt{5}}{10}\right)\,,
\end{align}
with $\tilde{q} = e^{\frac{2\pi i }{5\tau}}$.

Once again, we see that the complete weak/strong coupling resurgence properties are fully captured by the Fricke involution \eqref{eq:TSQM}. However, we now notice that for $N\geq 5$ we find in general a linear combination of twisted Lambert series with $s=1$ and different characters $\chi_N$ modulo $N$. In particular, for $N=5$ we encounter in \eqref{eq:GF31}  both the twisted Lambert series with $s=1$ and the Fricke-paired odd complex conjugated characters $\chi_{5,2}$ and $\chi_{5,3} = \overline{\chi_{5,2}}$, as well as the iterated twisted Eisenstein case studied in \eqref{eq:L1chi5} and corresponding to the case $s=1$ with even character $\chi_{5,4}$. Although this second term does not produce a factorially divergent asymptotic expansion neither at weak nor at strong coupling, it does contribute in a crucial way to the non-perturbative effects. 

The reason why the weak/strong coupling expansions Stokes constants of \eqref{eq:rhonm} analyzed in \cite{Fantini:2025wap} do not necessarily give rise to $L$-series for $N=m+n+1\geq 5$ comes from the fact that the weak/strong coupling non-perturbative effects are in general given by sums of the building blocks \eqref{eq:gq}-\eqref{eq:fq} for both odd \textit{and} even characters modulo $N$, as elucidated in \eqref{eq:GF31}.

While in the present discussion we only encounter the twisted Lambert series $\mathcal{L}_1(\chi_N;q)$ and $\tilde{\mathcal{L}}_1(\chi_N;q)$, it would be extremely interesting to understand whether the quantum Fricke inversion \eqref{eq:TSQM}, or the more general \eqref{eq:TSgen}, play any role in the weak/strong coupling modular resurgent structure of other topological string observables also for different values of the parameters $s,s_1,s_2$.

\section*{Acknowledgments}
We thank Athanasios Bouganis, Veronica Fantini, Axel Kleinschmidt, Claudia Rella, and Oliver Schlotterer for helpful discussions and  useful comments on the draft.
We are grateful to the Galileo Galilei Institute for Theoretical
Physics for the hospitality and the INFN for partial support during the GGI programme “Resurgence
and Modularity in QFT and String Theory” where the present work started. DD is supported by
the Royal Society under the grants ICA$\backslash$R2$\backslash$242058 and IEC$\backslash$R3$\backslash$243103.

\appendix

\section{Dirichlet characters and their $L$-functions}
\label{sec:DirApp}

In this appendix we present some key properties of Dirichlet characters and their associated $L$-functions. We refer to \cite{Mont} for a more detailed account.

\subsection{Dirichlet characters}
\label{sec:DirCh}

Dirichlet characters provide an important class of arithmetic functions with particular relevance for number theory.
\\

\textbf{Definition}: A function $\chi:\mathbb{Z}\to \mathbb{C}$ is a Dirichlet character of modulus $r\in \mathbb{N}$ if $\forall \,a,b\in \mathbb{Z}$:
\begin{itemize}
\item $\chi(a)\chi(b) = \chi(a\cdot b)$, i.e. $\chi$ is a completely multiplicative function;
\item $ \chi(a) =  0$ when $\mbox{gcd}(a,r)>1$, while $\chi(a)\neq 0$  when $\mbox{gcd}(a,r)=1$;
\item $\chi(a+r) = \chi(a)$, i.e. $\chi$ is $r$-periodic.
\end{itemize}

Some key properties of Dirichlet characters are:
\begin{itemize}
\item[(i)] Since $\mbox{gcd}(1,r)=1$ we necessarily have $\chi(1)\neq0$, and from complete multiplicativity we deduce $\chi(1)^2 = \chi(1)\Rightarrow\chi(1)=1$;
\item[(ii)] Since $a^{\phi(r)} \equiv 1 (\mbox{mod}\,r)$ we have that $\chi(a)$ must be a $\phi(r)^{th}$ root of unity with $\phi$ denoting Euler totient function;
\item[(iii)] For every modulus $r$ we have $\phi(r)$ distinct Dirichlet characters which can be denoted by $\chi^{(\ell)}_{r}$ with $\ell\in \mathbb{N}$ and $1\leq \ell \leq \phi(r)$. We note that this labelling does not coincide with the more commonly used Conrey labels (see \cite{LMFDB}).
\end{itemize}

Since $\chi_r(1)=1$ we furthermore deduce that $\chi_r(-1)^2 =\chi_r(r-1)^2 = 1$, i.e. $\chi_r(-1) = \chi_r(r-1) = \pm1$. We write $\chi_r(-1) =\chi_{r}(r-1)= (-1)^\kappa$ with $\kappa\in \{0,1\}$ and we say that a character is \textit{even} if $\kappa=0$, i.e. $\chi_r(r-1) = \chi_r(-1) = +1$, otherwise the character is \textit{odd} if $\kappa=1$, i.e. $\chi_r(r-1) = \chi_r(-1) = -1$.

For every modulus $r$ we have a special character called the \textit{principal character} which we denote by $\chi_{r,1}$ using Conrey labels as on the LMFDB \cite{LMFDB} such that:
\begin{equation}
\chi_{r,1}(a) =
\begin{cases}
&0\,,\qquad \mbox{gcd}(a,r)>1\,,\\
&1 \,,\qquad \mbox{gcd}(a,r)=1\,.
\end{cases}
\end{equation}

As stated at point (iii), for every modulus $r$ we have $\phi(r)$ distinct Dirichlet characters. There is a useful labelling system introduced by Conrey to denote each of the $\phi(r)$ Dirichlet characters modulo $r$  and which makes it easy to recover many important properties of the Dirichlet character directly from its label.
These labels are based on the construction of an explicit isomorphism between the group of Dirichlet characters modulo $r$ and the multiplicative group $(\mathbb{Z}/r \mathbb{Z})^\times$.

We shall write a Dirichlet character modulo $r$ in \textit{Conrey labels} as $\chi_{r,n}$ and following \cite{LMFDB} we shall define the labels as follows.
Firstly, $\chi_{r,1}$ always denotes the principal character modulo $r$. 
Then we move to consider the case where $r$ has a single prime factor.
For prime powers $r = p^k$, the Conrey labels $\chi_{r,n}$ are defined as:
\begin{itemize}
\item For an odd prime $p$ we pick the least positive integer $\ell_p$ which is a primitive root\footnote{For all odd primes $p$ the multiplicative group $(\mathbb{Z} / p^k \mathbb{Z})^\times$ is cyclic for all $k\geq1$. Furthermore, any integer $\ell$ which generates $(\mathbb{Z}/p^2 \mathbb{Z})^\times$ will also generate $(\mathbb{Z} / p^k \mathbb{Z})^\times$ with $ k \geq 1$. The least of such $\ell$, which we shall denote by $\ell_p$, is called the primitive root modulo $p$.} for all $p^k$. For $n\equiv \ell_p^a $ and  $m\equiv \ell_p^b $ coprime to $p$ we define the character $\chi_{p^k,n}$ in Conrey labels as
\begin{equation}
\chi_{p^k,n}(m) = \exp\left( 2\pi i \frac{a b}{\phi(p^k)}\right)\,.
\end{equation}
\item For $p=2$ we start with $\chi_{2,1}$ which is the trivial principal character, while $\chi_{4,3}$ is the unique non-trivial character modulo $4$. When $r=2^k$ with $k>2$ the multiplicative group modulo $r$ is not cyclic since $(\mathbb{Z}/2^k\mathbb{Z})^\times \cong \mathbb{Z}_2 \times \mathbb{Z}_{2^{k-2}}$. Hence we pick $-1$ and $5$ as the generators of the multiplicative group modulo $r$ and we write $n\equiv e_a 5^a \,{\rm mod} \,r$ and $m\equiv e_b 5^b \,{\rm mod}\, r$ with $e_a,e_b\in \{\pm1\}$.
For $k>2$ we then define the Dirichlet character $\chi_{2^k,n}$ in Conrey labels
\begin{equation}
\chi_{2^k,n}(m) = \exp\left[ 2\pi \left( \frac{(e_a-1)(e_b-1)}{8} + \frac{ab}{2^{k-2}} \right) \right]\,.
\end{equation}
\end{itemize}

 Lastly, for general $r$, the Dirichlet character $\chi_{r,n}$ in Conrey labels is defined multiplicatively by $\chi_{r,n}(m) = \chi_{r_1,n}(m) \chi_{r_2,n}(m)$ for all coprime positive integers $r_1,r_2$ such that $r=r_1 r_2$. The Chinese remainder theorem ensures that in this way we do in fact define a Dirichlet character and that every Dirichlet character arises in this manner. In particular, each Dirichlet character $\chi_r$ modulo $r$ admits a unique factorization as a product of Dirichlet characters whose moduli are prime powers.

From property (ii) we see that in general $\chi_r(a) \in\mathbb{C}$, hence we say that a character is \textit{real} or \textit{quadratic} if $\chi_r(a)\in\mathbb{R}$, otherwise we say it is \textit{complex}. Note that for a real character we necessarily have $\chi_r(a)\in \{-1,1\}$.
Again from property (ii) we see that for every character $\chi^{(\ell)}_{r}$ we have $(\chi^{(\ell)}_{r})^{\phi(r)} = \chi_{r,1}$. We define the \textit{order} of a character $\chi_r$ as the smallest integer $n$ such that $\chi_r^n = \chi_{r,1}$. Real characters have order $2$.

We want furthermore to distinguish some sort of ``fundamental'' characters, i.e. characters of modulus $r$ which do not descend from some other character with modulus $D<r$. To this end we need to define:
\\

\textbf{Definition}: We say that $\chi_r$ has a quasi-period $D$ if $\chi_r(a) = \chi_r(b)$ for all $a,b$ coprime to $r$ such that $a\equiv b\,(\mbox{mod}\,D)$. The smallest $D$ quasi-period of $\chi_r$ is called the \textit{conductor} of $\chi_r$.
If the modulus and the conductor of a character are the same we say that the character is \textit{primitive}, otherwise we say that the character is imprimitive and it is induced by a character with smaller modulus, i.e. if $\chi_r$ is imprimitive we have:
\begin{equation}\label{eq:Induced}
\chi_{r}(a) =\begin{cases}
&0\,,\qquad \qquad\, \mbox{gcd}(a,r)>1\,,\\
&\chi_D(a) \,,\qquad  \mbox{gcd}(a,r)=1\,,\\
\end{cases}
\end{equation}
with $D$ equal to the conductor of $\chi_r$ and $\chi_D$ a primitive character modulo $D$.

Note that the principal character modulo $r$ is imprimitive and it is induced by $\chi_{1,1}$ which is the trivial character, i.e. $\chi_{1,1}(n) = 1$ for all $n\in \mathbb{Z}$.
A Dirichlet character is real if and only if it can be written as a  Kronecker symbol $\left( \frac{D}{\bullet}\right)$ for some integer $D\in \mathbb{Z}$.
Furthermore we have a complete characterization of all real and primitive characters which are in one-to-one correspondence with quadratic number fields $\mathbb{Q}(\sqrt{d})$ with $d$ a square-free integer which in turn are in one-to-one correspondence with \textit{fundamental discriminants}, i.e. $D\in \mathbb{Z}$ such that:
\begin{itemize}
\item $D=d$ if $d\equiv 1\,(\mbox{mod}\,4)$ and $d$ square-free or
\item $D= 4d$ if $d\equiv 2 $ or $3\,(\mbox{mod}\,4)$ and $d$ square-free.
\end{itemize}

For a given modulus $r\in \mathbb{N}$, a real and primitive characters modulo $r$ exists if and only if $ r = |D|$ with $D$ a fundamental discriminant.
We can be more specific and consider a square-free number $d = p_1 p_2 \cdots p_\ell$ with $p_1<p_2<...<p_\ell$ odd primes.
All real and primitive characters are given by
\begin{itemize}
\item For $r=d$ we either have $r \equiv 1 (\mbox{mod}\,4)$, in which case $r$ is a fundamental discriminant and hence $\chi_r(a) = \left( \frac{r}{a}\right)$ is the unique real primitive character modulo $r$, or $r \equiv 3 (\mbox{mod}\,4)$ in which case $-r$ is a fundamental discriminant and hence $\chi_r(a) = \left( \frac{-r}{a}\right)$ is the unique real primitive character modulo $r$;
\item for $r=4d$ we either have $d \equiv 1 (\mbox{mod}\,4)$ in which case $-r= 4(-d)$ is a fundamental discriminant and hence $\chi_r(a) = \left( \frac{-r}{a}\right)$  is the unique real primitive character modulo $r$, or $d \equiv 3 (\mbox{mod}\,4)$ in which case $r=4d$ is a fundamental discriminant and hence $\chi_r(a) = \left( \frac{r}{a}\right)$ is the unique real primitive character modulo $r$;
\item for $r=8d$ we have that $r=4(2d)$ and $-r=4(-2d)$ are both fundamental discriminants and hence there are exactly two real primitive characters modulo $r$ which are $\chi_r(a) = \left( \frac{r}{a}\right)$ and $\chi_r(a) = \left( \frac{-r}{a}\right)$;
\item for $r=2d$ or $r= 2^\alpha d$ with $\alpha>3$ or $r=d p_i$, i.e. for $r$ non-square free, there are no real primitive characters.
\end{itemize}

We note furthermore that if $r$ equals $D$ with $D$ a positive fundamental discriminant such as $1,5,8,12,...$ (see the OEIS sequence \href{https://oeis.org/A003658}{$A003658$}) the corresponding real primitive character $\chi_r =  \left( \frac{D}{\bullet}\right)$  is even. Conversely if $r$ equals $|D|$ with $D$ a negative fundamental discriminant, such as $-3,-4,-7,-8,...$ (see the OEIS sequence \href{https://oeis.org/A003657}{$A003657$}) the corresponding real primitive character $\chi_r =  \left( \frac{D}{\bullet}\right)$  is odd.

\subsection{Dirichlet $L$-functions}

Given a Dirichlet character $\chi$ we define its associated $L$-function as
\begin{equation}
L(\chi,s) \coloneqq \sum_{n=1}^\infty \frac{\chi(n)}{n^s}\,,\label{eq:Lseries}
\end{equation}
which converges for $\mbox{Re}(s)>0$ for a general character, while for $\chi$ a principal character it converges only for $\mbox{Re}(s)>1$.
To analyze further properties of Dirichlet $L$-functions we need to distinguish between primitive and imprimitive characters.

We start by assuming that $\chi=\chi_r$ is a primitive character modulo $r$ where $r>1$.
We define the completed Dirichlet $L$-function as
\begin{equation}
\xi(\chi_r,s) \coloneqq \Big(\frac{r}{\pi}\Big)^{\frac{s+\kappa}{2}} \Gamma\Big(\frac{s+\kappa}{2}\Big) L(\chi_r,s)\,,\label{eq:Lcompleted}
\end{equation}
where $\kappa\in \{0,1\}$ is defined via $\chi_r(-1)=(-1)^\kappa $, i.e. for even characters $\kappa=0$ while $\kappa=1$ for odd characters.
The completed $L$-function of a primitive character satisfies the functional equation (cfr. (10.19) of \cite{Mont})
\begin{equation}
\xi(\chi_r,s) = \epsilon(\chi_r) \xi(\bar{\chi}_r,1-s)\,,\label{eq:functional}
\end{equation}
where
\begin{equation}
\epsilon(\chi_r) \coloneqq \frac{\rho(\chi_r)}{i^\kappa \sqrt{r}} \,,\qquad\rho(\chi_r)\coloneqq \sum_{n=1}^r \chi_r(n) e^{\frac{2\pi i n}{r}}\,.\label{eq:Gauss}
\end{equation}
Alternatively, we can also use the functional equation (cfr. Corollary 10.9 of \cite{Mont})
\begin{equation}
L(\chi_r,s) = \epsilon(\chi_r)  2^s \pi^{s-1} r^{\frac{1}{2}-s} \Gamma(1-s) \sin\big(\pi\tfrac{(s+\kappa)}{2}\big) L(\bar{\chi}_r,1-s)\,.\label{eq:LFunct}
\end{equation}

Note that for the trivial character $\chi_r = \chi_{1,1}$, the Dirichlet series \eqref{eq:Lseries} reduces to the standard definition of the Riemann zeta function, $\zeta(s) = \sum_{n=1}^\infty n^{-s}$, and \eqref{eq:Lcompleted} produces the completed Riemann zeta function $\xi(s)=\pi^{-s/2}\Gamma(s/2)\zeta(s)$ which satisfies the well-known reflection identity $\xi(s) = \xi(1-s)$.

Furthermore we note that for $\chi_r$ a primitive character modulo $r$ we have:
\begin{itemize}
\item $|\epsilon(\chi_r)|=1$ and $\epsilon(\chi_r) \epsilon(\bar{\chi}_r) = 1$;
\item $\epsilon(\chi_r)=1$ for a real primitive character.
\end{itemize}

From the functional equation \eqref{eq:functional}-\eqref{eq:LFunct} we can easily deduce the location of the trivial zeroes of the Dirichlet $L$-functions:
\begin{equation}
L(\chi_r, -\kappa -2k) = 0\,,\qquad k\in\mathbb{N}\,,\label{eq:TZ}
\end{equation}
i.e. for even characters the associated $L$-functions vanish at negative even integers as well as at zero, while for odd characters the associated $L$-functions vanish at negative odd integers.
Using an integral contour argument, see \cite{Zuck}, it is possible to explicitly evaluate all Dirichlet $L$-values at negative integers for primitive characters which can be expressed in terms of Bernoulli polynomials ${\rm B}_k(x)$:
\begin{equation}
L(\chi_r, -k) = (-1)^k \frac{r^k}{k+1} \sum_{n=1}^r \chi_r(n) {\rm B}_{k+1}\left(1-\frac{n}{r}\right)\,,\label{eq:NegInt}
\end{equation}
valid for $k\in \mathbb{N}$. Note that while for real primitive characters these are all rational numbers, for complex primitive characters these are in general complex algebraic  numbers.

Combining \eqref{eq:NegInt} with the functional equation \eqref{eq:functional}-\eqref{eq:LFunct}, we can also read the ``easy'' non-critical primitive $L$-values, i.e. the analogues of even zeta values:
\begin{equation}
L(\chi_r,2k-\kappa) =\epsilon(\bar\chi_r) \frac{ (-1)^{k-1} 2^{2k-1-\kappa} \pi^{2k-\kappa} r^{-\frac{1}{2}}}{(2k-\kappa)!} \sum_{n=1}^r \bar{\chi}_r(n) {\rm B}_{2k-\kappa}\left(1-\frac{n}{r}\right)\,, \label{eq:PosInt}
\end{equation}
where again $\kappa=0$ for even characters, and $\kappa=1$ for odd ones. These are positive even $L$-values for even characters and positive odd $L$-values for odd characters.
If we specialize \eqref{eq:PosInt} to the case of real primitive characters for which $\epsilon(\chi_r) = 1$ and $\bar{\chi}_r = \chi_r$ we see that $L(\chi_r,2k-\kappa) $ is a rational multiple of $\pi^{2k-\kappa}$.

For $\chi_r$ an imprimitive character modulus  $r$ induced by the primitive character $\chi_D$ of modulus $D$ as in \eqref{eq:Induced} we have the relation between $L$-functions:
\begin{equation}
L(\chi_r,s) = L(\chi_D,s) \prod_{p\vert r} \Big( 1- p^{-s} \chi_D(p) \Big)\,,\label{eq:Lnonprim}
\end{equation}
where the product is over all prime factors $p$ of $r$. Note however that if $p$ divides $D$  we have $\chi_D(p)  =0$ hence we can restrict the product over all prime factors of $r$ coprime with $D$ and, by expanding the product, \eqref{eq:Lnonprim} is equivalent to
\begin{equation}
L(\chi_r,s) = L(\chi_D,s)\sum_{d|r}\mu(d)\frac{\chi_D(d)}{d^{s}}\,.
\end{equation}

By combining \eqref{eq:TZ} with \eqref{eq:Lnonprim} we see that also for an imprimitive character its $L$-values at negative even/odd integer vanish according to whether the character is even/odd. However, for the principal character $\chi_{r,1}$ its $L$-value at $s=0$ does not vanish.

\section{Twisted Eisenstein series}
\label{sec:Twist}

In this appendix we review some known properties regarding two families of \textit{twisted Eisenstein series}. 

Following \cite{Hida} we start by considering a primitive character $\chi_r$ modulo $r$ and define a first family of holomorphic Eisenstein series twisted by the character $\chi_r$ in terms of the lattice sum
\begin{equation}
E'_k(\chi_r;\tau) \coloneqq \sum_{(n,m)\neq(0,0)} \frac{\chi_r^{-1}(m)}{(r n \tau + m)^k}\,,\label{eq:EpLattice}
\end{equation}
with $k\in \mathbb{N}$. The series converges absolutely for $k>2$.
It is easy to show that the twisted Eisenstein series $E'_k(\chi_r;\tau)  \in \mathcal{M}_k (\Gamma_0(r),\chi_r)$, i.e. $E'_k(\chi_r;\tau)$ belongs to the vector space of modular forms of weight $k$ and character $\chi_r$ with respect to the congruence subgroup $\Gamma_0(r)$, that is
\begin{equation}
E'_k(\chi_r;\gamma\cdot \tau) = \chi_r(d) (c\tau+d)^k E'_k(\chi_r;\tau)\,,\qquad \gamma =\left(\begin{matrix} a & b \\ c & d \end{matrix}\right)\in \Gamma_0(r)\,,
\end{equation}
where as usual
\begin{equation}
\gamma\cdot \tau \coloneqq \frac{a\tau+b}{c\tau+d}\,,\qquad {\rm with}\quad \gamma =\left(\begin{matrix} a & b \\ c & d \end{matrix}\right)\,.
\end{equation}
The congruence group $\Gamma_0(r)$ is defined as $\Gamma_0(r) \coloneqq \{ \left(\begin{smallmatrix} a & b \\ c & d \end{smallmatrix}\right) \in {\rm SL}(2,\mathbb{Z})\,\vert \, c \equiv 0 \,{\rm mod} \,r\}$.

It is convenient to use a slightly different normalization for this first type of twisted Eisenstein series and define
\begin{align}
E_k(\chi_r ; \tau) &\notag \coloneqq \left( 2 r^{-k} \rho(\chi^{-1}_r) \frac{(-2\pi i )^k}{(k-1)!}\right) ^{-1} E'_k(\chi_r ; \tau)\\
&\phantom{:} = \frac{L(\chi_r,1-k)}{2} + \sum_{n=1}^\infty \sigma_{k-1,\chi_r}(n) q^n\,,\label{eq:TwEisen1}
\end{align}
where $q=e^{2\pi i \tau}$ and we note that the inverse character trivially satisfies $\chi^{-1}_r(n) = \bar{\chi}_r(n)$. The Gauss sum $\rho(\chi_r)$ appearing in the above equation has been defined in \eqref{eq:Gauss}, while the twisted divisor sigma function $\sigma_{s,\chi_r}(n)$ is given by
\begin{equation}
\sigma_{s,\chi_r}(n) : = \sum_{d|n} \chi_r(d) d^s\,. \label{eq:Twsigma1App}
\end{equation}
Note that the above expression is also valid for $k=1$ and $k=2$ in the case of a non-trivial character.

We now define a second family of twisted Eisenstein series via the lattice sum
\begin{align}
G'_k(\chi_r;\tau) &\coloneqq \sum_{(m,n)\neq(0,0)} \frac{\chi_r(n)}{(n\tau+m)^k} \label{eq:GpLattice}\,,
\end{align}
It is again fairly straightforward to show that 
$G'_k(\chi_r;\tau)  \in \mathcal{M}_k (\Gamma_0(r);\chi_r)$, i.e.
\begin{equation}
G'_k(\chi_r;\gamma\cdot \tau) = \chi_r(d) (c\tau+d)^k G'_k(\chi_r;\tau)\,,\qquad \forall\,\gamma =\left(\begin{matrix} a & b \\ c & d \end{matrix}\right)\in \Gamma_0(r)\,,
\end{equation}
where we note that in the case of a non-trivial character \eqref{eq:GpLattice} does define a modular form also for $k=1$ and $k=2$.
We find it convenient to choose a slightly different normalization for this second twisted Eisenstein series and define
\begin{equation}
G_k(\chi_r;\tau)  \coloneqq \left(  \frac{2(-2\pi i )^k}{(k-1)!}\right)^{-1} G'_k(\chi_r;\tau)   = \delta_{k,1} \frac{L(\chi_{r};0)}{2} + \sum_{n=1}^\infty \sigma'_{k-1,\chi_r}(n) q^n\,,\label{eq:TwEisen2}
\end{equation}
where we defined a second flavour of twisted divisor sigma function
\begin{equation}
\sigma'_{s,\chi_r}(n) : = \sum_{d|n} \chi_r\left(\frac{n}{d}\right) d^s\,. \label{eq:Twsigma2App}
\end{equation}

Note for $k=1$ there is no difference between the twisted divisor functions \eqref{eq:Twsigma1}-\eqref{eq:Twsigma2} so that trivially 
\begin{equation}
G_1(\chi_r;\tau) =E_1(\chi_r;\tau)\,,
\end{equation}
 as we can see from  \eqref{eq:TwEisen1}-\eqref{eq:TwEisen2}.
Furthermore for the trivial character $\chi_{r}(n)= \chi_{1,1}(n)=1$ and when $k\in \mathbb{N}^{\geq4}$ is an even integer, we have that both Eisenstein series \eqref{eq:EpLattice}-\eqref{eq:GpLattice} reduce to the standard holomorphic Eisenstein series,
\begin{align}
 G'_{k}(\tau) &\label{eq:holoGlattice}=  \frac{2(-2\pi i )^k}{(k-1)!} G_{k}(\tau) \coloneqq \sum_{(n,m)\neq (0,0)} \frac{1}{(n\tau+m)^k}\,,\\
 G_{k}(\tau)  &=  -\frac{B_{k}}{2k} + \sum_{N=1}^\infty \sigma_{k-1}(N) q^N\,,\label{eq:holoG}
\end{align}
where $B_k$ denotes the $k^{th}$ Bernoulli number.
 \vspace{0.2cm}
 
 \textbf{Remark.}
It is important to note that if we change summation variables in the lattice sums \eqref{eq:EpLattice}-\eqref{eq:GpLattice}  to $(m,n) \to (-m,-n)$, we deduce that both twisted Eisenstein are non-zero only for 
\begin{equation}
 k \equiv\kappa \, \rm{mod}\,2\,,\label{eq:WeightKappa}
\end{equation}
where again $\chi_r(-1) = (-1)^\kappa$ with $\kappa\in\{0,1\}$ according to the character being even/odd.
That is for an even, respectively odd, character we can only consider twisted Eisenstein series with even, respectively odd, weight $k$. \vspace{0.2cm}

Let us now define the slash action via
\begin{equation}
f(\tau) \vert_k \gamma \coloneqq {\rm det}(\gamma)^{k-1} (c\tau+d)^{-k} f(\gamma\cdot \tau)\,, \qquad \gamma =\left(\begin{matrix} a & b \\ c & d \end{matrix}\right)\in M_2(\mathbb{R})\,.
\end{equation}
A key relation between $E_k(\chi_r;\tau)$ and $G_k(\chi_r;\tau) $ is obtained by considering the action of Fricke involution $\hat{S}_r =\left(\begin{matrix} 0 & -1 \\ r & 0 \end{matrix}\right)$ which reads
\begin{equation}
E_k(\chi_r;\tau) \Big\vert_k \hat{S}_r = \chi_r(-1) \,r^{k-2} \rho(\chi_r) G_k(\bar{\chi}_r;\tau)\,.\label{eq:Fricke1}
\end{equation}
Alternatively, we combine equation \eqref{eq:Gauss} with the fact that $\chi_r$ is a primitive character and with the weight/parity condition \eqref{eq:WeightKappa} to rewrite \eqref{eq:Fricke1} as
\begin{equation}
G_k(\chi_r ;\tau) = i^ \kappa \epsilon(\chi_r) r^{\frac{1}{2}-k} \tau^{-k} E_k(\bar{\chi}_r;-\frac{1}{r\tau})\,.\label{eq:FrickeApp}
\end{equation}
Throughout this paper we use interchangeably as variables either the modular
parameter $\tau\in \mathbb{C}$ such that ${\rm Im}(\tau)>0$ or the variable $q= e^{2\pi i \tau}\in \mathbb{C}$ with $|q|<1$.
By slight abuse of notation we write
\begin{equation}
E_k(\chi_r ; \tau) = E_k(\chi_r;q) \,, \qquad G_k(\chi_r ; \tau) = G_k(\chi_r;q) \,, \qquad {\rm with}\quad  q= e^{2\pi i \tau}\,.
\end{equation}

Note that the lattice sums \eqref{eq:EpLattice}-\eqref{eq:GpLattice} can be seen as special cases of a more general lattice sum involving two Dirichlet character.
Following \cite{Miyake} we define 
\begin{equation}
G_k'(\chi_{r_1},\chi_{r_2};\tau) \coloneqq \sum_{(n,m)\neq (0,0)} \frac{\chi_{r_1}(n) \chi_{r_2}^{-1}(m)}{(nr_2 \tau+ m)^k}\,,\label{eq:E2chi}
\end{equation}
where $k\in \mathbb{N}$ and $\chi_{r_1},\chi_{r_2}$ are two primitive characters modulo $r_1$ and $r_2$ respectively.
Note that this lattice sum vanishes unless
\begin{equation}
k \equiv \kappa_1 +\kappa_2 \text{ mod }2\,,
\end{equation}
where $\chi_{r_1}(-1) = (-1)^{\kappa_1}$ and $\chi_{r_2}(-1) = (-1)^{\kappa_2}$. Furthermore, we note that when $\chi_{r_1} =\chi_{r_2} = \chi_{1,1}$ the lattice sum \eqref{eq:E2chi} is non-trivial only for $k$ even and it is absolutely convergent only when $k\geq4$.

From the definitions \eqref{eq:EpLattice}-\eqref{eq:GpLattice} and \eqref{eq:holoGlattice}, it is immediate to see that when either one or both of the two characters in \eqref{eq:E2chi} become the trivial character $\chi_{1,1}$ we retrieve the Eisenstein series previously discussed, i.e.
\begin{equation}
G_k'(\chi_{r},\chi_{1,1};\tau) = G_k'(\chi_r;\tau)\,,\qquad G_k'(\chi_{1,1},\chi_{r};\tau) = E_k'(\chi_r;\tau)\,,\qquad G_k'(\chi_{1,1},\chi_{1,1};\tau) =G'_k(\tau)\,.
\label{eq:G2chiDict}
\end{equation}
Furthermore, it is an easy exercise to show that $G_k'(\chi_{r_1},\chi_{r_2};\tau)\in \mathcal{M}_k(\Gamma_0(r_1r_2),\chi_{r_1}\chi_{r_2})$ and it satisfies the Fricke inversion
\begin{equation}
G_k'\left(\chi_{r},\chi_{r_2}; -\frac{1}{r_1 r_2\tau}\right) = (-1)^{\kappa_1} (r_1 \tau)^k G'_k(\bar{\chi}_{r_2},\bar{\chi}_{r_1};\tau)\,.\label{eq:FrickeGp2chi}
\end{equation}

Rather than using the definition \eqref{eq:E2chi}, we introduce the convenient normalization
\begin{equation}
G_k(\chi_{r_1},\chi_{r_2};\tau) = \left( 2 r_2^{-k} \rho(\chi^{-1}_{r_2}) \frac{(-2\pi i )^k}{(k-1)!}\right)^{-1} G'_k(\chi_{r_1},\chi_{r_2};\tau)\,,
\end{equation}
for which we can write
\begin{equation}
G_k(\chi_{r_1},\chi_{r_2};\tau) = A_k(\chi_{r_1},\chi_{r_2})+ \sum_{n=1}^\infty  \sigma_{k-1}(\chi_{r_1},\chi_{r_2};n)\,q^n\,,\label{eq:qserG2chiApp}
\end{equation}
where we defined the doubly twisted divisor sigma function
\begin{equation}
\sigma_s(\chi_{r_1},\chi_{r_2};n) \coloneqq \sum_{d|n} \chi_{r_1}\left(\frac{n}{d}\right) \chi_{r_2}(d) d^{s}\,.\label{eq:div2chiApp}
\end{equation}
We note that the doubly twisted divisor sigma function satisfies the identity
\begin{equation}
\sigma_s(\chi_{r_1},\chi_{r_2};n)  = n^{s} \sigma_{-s}(\chi_{r_2},\chi_{r_1};n)\,. \label{eq:div2chiId}
\end{equation}
The constant term in \eqref{eq:qserG2chiApp} depends on the nature of the characters considered
\begin{equation}
A_k(\chi_{r_1},\chi_{r_2}) = \begin{cases}
\frac{L(\chi_{r_2},1-k)}{2} \,,\qquad \,\,\,\,\chi_{r_1} = \chi_{1,1}\,,\\
\delta_{k,1} \frac{L(\chi_{r_1},0)}{2} \,,\qquad \,\chi_{r_2} = \chi_{1,1}\,,\\
0\,,\qquad \qquad\qquad\,\, {\rm otherwise}\,.
\end{cases}
\end{equation}

We conclude by noting that the Fricke inversion \eqref{eq:FrickeGp2chi} in the new normalization yields 
\begin{equation}
G_k(\chi_{r_1},\chi_{r_2}; \tau) = i^{\kappa_1 +\kappa_2} \epsilon(\chi_{r_1})\epsilon(\chi_{r_2}) r_1^{\frac{1}{2}-k} r_2^{-\frac{1}{2}}  \tau^{-k} G_k\left(\bar{\chi}_{r_2},\bar{\chi}_{r_2}; -\frac{1}{r_1 r_2\tau}\right)\,,\label{eq:FrickeG2chiApp}
\end{equation}
where we used \eqref{eq:Gauss} to rewrite the Gauss sums.
From \eqref{eq:G2chiDict}, it is immediate to see that by specialising \eqref{eq:FrickeG2chiApp} to the case $\chi_{r_2} = \chi_{1,1}$, for which we have $\epsilon(\chi_{1,1})=1$, we find precisely the Fricke inversion \eqref{eq:FrickeApp} relating the two families of twisted Eisenstein series previously discussed.

\section{Asymptotic expansion for a single character}
\label{sec:AsyApp}

In this appendix we show how to extract the asymptotic expansion of \eqref{eq:Repeated1} in the case when one of the two characters reduces to the trivial case using the analysis of \cite{Zagier}.
Given the symmetry \eqref{eq:LamGen}, without loss of generality we can assume that $\chi_{r_2}=\chi_{1,1}$ and $\chi_{r_1}=\chi_r$ is non-principal so that \eqref{eq:Repeated1} becomes
\begin{equation}
\Xi_{s_1,s_2}(\chi_{r},\chi_{1,1} ; y) =  y^{s_2} \sum_{m=1}^\infty \frac{1}{(my)^{s_2}} \Phi_{s_1}(\chi_{r} ; e^{-2\pi my})\,,\label{eq:Repeated1App}
\end{equation}
where the auxiliary function $\Phi_s(\chi_r;q)$ has been defined in \eqref{eq:Phi}.

Following \cite{Zagier}, we see that to extract the asymptotic expansion of \eqref{eq:Repeated1App} we need considering the asymptotic expansion for $y\to 0^+$ of a function defined via the series
\begin{equation}
F(y) = \sum_{m=0}^\infty f((m+a)y)\,,\label{eq:AsymDef}
\end{equation}
where $f$ is $C^\infty$ near the origin and $a>0$.
We assume that $f(y)$ is a smooth function for $y>0$ with all its derivatives of rapid decay at infinity and that it is has an asymptotic expansion around $y=0$ of the form $f(y) \sim \sum_{n=0}^\infty b_n y^n$. Under these conditions, the perturbative asymptotic expansion of \eqref{eq:Asym} around $y=0$ is given by
\begin{equation}
F(y) \sim\frac{I_f}{y}+ \sum_{n= 0}^\infty b_n \zeta(-n,a) y^n\,.\label{eq:Asym}
\end{equation}

The first term in this asymptotic expansion is the Riemann integral term
\begin{equation}
I_f \coloneqq \int_0^\infty f(y) \,{\rm d}y\,,
\end{equation}
while the infinite series over Hurwitz zeta functions arises from the Taylor expansion around the origin
\begin{equation}
f((m+a)y) = \sum_{n= 0}^\infty b_n (m+a)^n y^n \,,
\end{equation}
which would naively give
\begin{equation}
\sum_{m= 0}^\infty f((m+a)y) \stackrel{\mbox{\scriptsize{naive}}}{=}\sum_{n=0}^\infty \sum_{m=0}^\infty b_n (m+a)^n y^n = \sum_{n=0}^\infty b_n  \zeta(-n,a) y^n\,.\label{eq:naive}
\end{equation}
By combining the Riemann integral term with the naive expansion we find the perturbative expansion \eqref{eq:Asym}. Importantly \eqref{eq:Asym} does not account for exponentially suppressed corrections as $y\to 0^+$.

For later convenience we also consider the case where $f(y)$ has a singularity at the origin of the form $f(y) = b_{-1} y^{-1}+O(y^0)$ for which the naive contribution \eqref{eq:naive} becomes singular. 
This case can be easily cured by adding and subtracting to $f(y)$ an auxiliary function, such as $b_{-1} e^{-y}/y$, for which we can compute \eqref{eq:Asym} exactly. In this case the net result is that the asymptotic expansion of \eqref{eq:AsymDef} becomes
\begin{equation}
F(y)  \sim \frac{1}{y}\left( I^\star_f    - b_{-1} \log y  \right)+ \sum_{n= 0}^\infty b_n \zeta(-n,a) y^n\,,\label{eq:Asym2}
\end{equation}
where the Riemann contribution is now given by the regulated term
\begin{equation}
I^\star_f \coloneqq \int_0^\infty \Big[ f(y) - \frac{e^{-y}}{y} \Big] {\rm d}y\,.
\end{equation}
Note that more singular terms in the asymptotic expansion of $f(y)$ near $y=0$, i.e. terms of the form $f(y)\sim b_{\lambda} y^\lambda$ with ${\Re}(\lambda)<-1$, are not particularly interesting since they can be subtracted off from $f(y)$ given that their contribution to \eqref{eq:Asym} is absolutely convergent.

Going back to the case of interest  \eqref{eq:Repeated1App}, we then need to understand the behaviour of \eqref{eq:Phi} as $y \to 0^+$.
This can be easily done by making use of its Mellin-Barnes integral representation
\begin{equation}
\Phi_{s}(\chi_r ; y)  = \int_{\gamma-i \infty}^{\gamma+i\infty} \Gamma(t) L(\chi_r, t+s) (2\pi y)^{-t} \,\frac{{\rm d}t}{2\pi i}\,,\label{eq:MellinApp}
\end{equation}
where $\gamma\in \mathbb{R}$ is chosen such that $\gamma >0$ for $\chi_r$ non-primitive.
With this choice of integration contour we can replace the $L$-function in terms of its convergent $L$-series~\eqref{eq:Lseries} and then use the dominated convergence theorem to exchange the series with the contour integral to reproduce~\eqref{eq:Phi}.

To analyze the small-$y$ expansion of \eqref{eq:Phi}, we need to close the contour of integration in \eqref{eq:MellinApp} to the left semi-half plane ${\rm Re}(t)<0$ and collect the residues coming from the poles of the integrand originating solely from the gamma function $\Gamma(t)$ given that $L(\chi_r,s)$ is an entire function of $s$ for $\chi_r$ not-principal.
We then have that $\Phi_{s}(\chi_r ; y) = L(\chi_r,s) +O(y)$ as $y\to 0^+$ exactly as one might have naively deduced by expanding the exponential directly in \eqref{eq:Phi}. Note in particular that for $s=-2n$ with $n\in \mathbb{N}$ and $\chi_r$ even we have $L(\chi_r,s) =0$, see \eqref{eq:TZ}, hence $\Phi_{s}(\chi_r ; y) = O(y)$ as $y\to 0^+$.

We can now move to analyze the asymptotic expansion of \eqref{eq:Repeated1App} using Zagier's method. 
Let us start first with the case $(s_1,s_2)=(s,0)$ as it is simpler.
For the Riemann integral contribution to  \eqref{eq:Repeated1App} when $(s_1,s_2)=(s,0)$ we have 
\begin{equation}
I_f (\chi_r;s)= \int_0^\infty \Phi_s(\chi_r,y) \,{\rm d}y = \frac{L(\chi_r,s+1)}{2\pi} \,,\label{eq:RiemannApp}
\end{equation}
while for the naive part of the asymptotic expansion we use
\begin{align}
\sum_{m=1}^\infty \Phi_s(\chi_r,my) &\notag \stackrel{\mbox{\scriptsize{naive}}}{=} \sum_{m=1}^\infty \sum_{n=1}^\infty \frac{\chi_r(n)}{n^s} e^{-2\pi n m y} = \sum_{k=0}^\infty \frac{(-2\pi y)^k}{k!}  \sum_{n=1}^\infty\frac{\chi_r(n)}{n^{s-k}}\sum_{m=1}^\infty m^k\\
& = \sum_{k=0}^\infty \frac{(-2\pi y)^k}{k!} L(\chi_r,s-k)\zeta(-k) \,.\label{eq:naive1}
\end{align} 
We note in particular that the Riemann term \eqref{eq:RiemannApp} can be obtained as the limit for $k\to -1$ of the naive expansion, i.e.
\begin{equation}
\lim_{k\to -1} \frac{(-2\pi y)^k}{k!}L(\chi_r,s-k) \zeta(-k)  = \frac{L(\chi_r,s+1)}{2\pi y}\,,
\end{equation}
so that the asymptotic expansion of \eqref{eq:LamChi} for $y\to 0^+$ is given by
\begin{equation}
 \mathcal{L}_s(\chi_r ; y) = \Xi_{s,0}(\chi_r,\chi_{1,1};y) \sim \mathcal{L}^{\rm Pert}_s(\chi_r ; y) =  \sum_{k=-1}^\infty \frac{(-2\pi y)^k}{k!}L(\chi_r,s-k) \zeta(-k) \,.\label{eq:AsyLApp}
\end{equation}

For the asymptotic expansion as $y\to 0^+$ of \eqref{eq:Repeated1App} with general parameters $(s_1,s_2)$ the story is analogous. However, a minor complications arises from the fact that for ${\rm Re}(s_2)\geq 1$ the expansion of the summand in \eqref{eq:Repeated1App} contains finitely many power behaved terms of the form $(my)^\lambda$ with ${\rm Re}(\lambda) \leq -1$. As explained at the beginning of this section, all such terms with ${\rm Re}(\lambda) < -1$ can be straightforwardly subtracted off given that their contribution to \eqref{eq:Repeated1App}  is absolutely convergent. However, if the expansion of the summand near $y=0$ contains a term proportional to $y^{-1}$, which is only possible for $s_2\in \mathbb{N}$, we find a logarithmic term upon summation as in \eqref{eq:Asym2}.
Both the contributions originating from absolutely convergent terms as well as a possible $\log(y)$ term can be obtained via a proper analytic continuation of \eqref{eq:Asym}.

Assuming that $\mbox{Re}(s_2)<1$, the Riemann term is now given by
\begin{equation}
I_f (\chi_r;s_1,s_2)= \int_0^\infty \frac{1}{y^{s_2}} \Phi_{s_1}(\chi_r,y) \,{\rm d}y = (2\pi)^{s_2-1} \Gamma(1-s_2) L(\chi_r,s_1+1-s_2)\,.\label{eq:RiemannApp2}
\end{equation}
Similarly, the perturbative corrections originating from the naive contribution to \eqref{eq:Repeated1App} are now
\begin{align}
& \notag y^{s_2} \sum_{m=1}^\infty \frac{1}{(my)^{s_2}}\Phi_{s_1}(\chi_r,my)  \stackrel{\mbox{\scriptsize{naive}}}{=} \sum_{m=1}^\infty \sum_{n=1}^\infty \frac{\chi_r(n)}{n^{s_1}} \frac{1}{m^{s_2}} e^{-2\pi n m y} \\
&= \sum_{k=0}^\infty \frac{(-2\pi y)^k}{k!}  \sum_{n=1}^\infty\frac{\chi_r(n)}{n^{s_1-k}} \sum_{m=1}^\infty \frac{1}{m^{s_2-k}}= \sum_{k=0}^\infty \frac{(-2\pi y)^k}{k!}  L(\chi_r,s_1-k) \zeta(s_2-k)\,.\label{eq:naive2}
\end{align} 

We then deduce that the general perturbative asymptotic expansion near $y\to 0^+$ of \eqref{eq:Repeated1App} is given by
\begin{align}
&\notag\Xi_{s_1,s_2}(\chi_r,\chi_{1,1};y) \sim  \Xi^{\rm Pert}_{s_1,s_2}(\chi_r,\chi_{1,1};y)  \\
&= (2\pi y)^{s_2-1}\Gamma(1-s_2) L(\chi_r,s_1+1-s_2)+ \sum_{k=0}^\infty \frac{(-2\pi y)^k}{k!}  L(\chi_r,s_1-k)\zeta(s_2-k)\,.\label{eq:XiPApp}
\end{align} 
We note that for general values of the parameters $(s_1,s_2)$, this is a factorially divergent asymptotic expansion.

As mentioned above, the case $s_2\in \mathbb{N}$ deserves some additional care since both the first term in \eqref{eq:XiPApp} as well as the summand with $k=s_2-1$ are divergent. However, their combined contribution is finite. This can be seen by considering the analytic continuation in $s_2$ as given in \eqref{eq:XiPApp} and then taking the limit $s_2\to m\in \mathbb{N}$ for which we have
\begin{align}
&\notag \lim_{s_2\to m} \left[ (2\pi y)^{s_2-1}\Gamma(1-s_2) L(\chi_r,s_1+1-s_2){+}\frac{(-2\pi y)^{m-1}}{(m-1)!} L(\chi_r,s_1+1-m)\zeta(1+s_2-m) \right]\\
&\label{eq:s2IntApp}= \left[ L'(\chi_r , s_1+1-m) + L(\chi_r , s_1+1-m)\left( \gamma+\psi(m) - \log(2\pi y)\right) \right]  \frac{(-2\pi y)^{m-1}}{(m-1)!}  \,,
\end{align}
where $\gamma$ is the Euler-Mascheroni constant, $\psi(x) = \Gamma'(x)/\Gamma(x)$ is the digamma function and $L'(\chi_r,s) = d L(\chi_r,s)/ds$. This result is identical to what can be derived using the particular expansion \eqref{eq:Asym2}.

Specialising \eqref{eq:XiPApp} to $(s_1,s_2) = (0,s)$, we find the asymptotic perturbative expansion of $\tilde{\mathcal{L}}_s(\chi_r;y)$,
\begin{align}
&\notag \tilde{\mathcal{L}}_s(\chi_r;y) = \Xi_{0,s}(\chi_{r},\chi_{1,1};y) \sim \tilde{\mathcal{L}}^{\rm Pert}_s(\chi_r;y)\\
& =(2\pi y)^{s-1}\Gamma(1-s) L(\chi_r,1-s)+ \sum_{k=0}^\infty \frac{(-2\pi y)^k}{k!}L(\chi_r,-k) \zeta(s-k) \,.\label{eq:LtPApp}
\end{align}
The same caveat applies here when considering the case $s\in \mathbb{N}$ where in the above expression we must consider the limit \eqref{eq:s2IntApp}.
Note that the asymptotic expansions \eqref{eq:AsyLApp}-\eqref{eq:XiPApp} and \eqref{eq:LtPApp} can all be derived as special cases of the general expression \eqref{eq:XiAsyP} presented in the main text upon substituting for the characters $\chi_{r_1}=\chi_r$ and $\chi_{r_2}=\chi_{1,1}$.

\end{document}